\documentclass[onefignum,onetabnum]{siamonline220329}

\usepackage{braket,amsfonts}
\usepackage[position=top]{subfig}
\usepackage{amssymb}
\usepackage{amsbsy}
\usepackage{stmaryrd}

\usepackage{algorithm}
\usepackage{algpseudocode}

\usepackage{pgfplots}
\usepackage{multicol}
\usepackage{hhline}
\usepackage{multirow}
\usepackage{enumerate}
\usepackage{enumitem}
\usepackage{caption}
\usepackage{float}
\usepackage{amssymb}
\usepackage{pifont} 
\usepackage{color}
\usepackage{chngcntr}
\counterwithin{table}{section}
\usepackage{amsbsy, mathrsfs, mathdots}
\newsiamthm{claim}{Claim}
\newsiamremark{remark}{Remark}
\newsiamremark{hypothesis}{Hypothesis}
\crefname{hypothesis}{Hypothesis}{Hypotheses}

\usepackage{bbm}

\usepackage[utf8]{inputenc} 
\usepackage[T1]{fontenc}    
\usepackage{lmodern}

\usepackage{booktabs}       
\usepackage{amsfonts}       
\usepackage{nicefrac}       
\usepackage{microtype}      
\usepackage{yhmath}         
\usepackage{titlesec}
\usepackage{mwe}
\usepackage{cleveref}

\usepackage{multirow}
\usepackage{hhline}
\usepackage{color}
\usepackage{epsfig}
\usepackage{subfig}
\usepackage{graphicx}   
\usepackage{xcolor}   

\usepackage{bm}
\usepackage{fixmath}

\usepackage{comment}
\usepackage[numbers,square,comma, sort&compress]{natbib}
\usepackage{epstopdf}

\usepackage{physics}
\usepackage{setspace}

\graphicspath{{eps/}{img/}{video/}}

\DeclareMathOperator*{\minimize}{\mathrm{minimize}}
\DeclareMathOperator*{\subject}{\mathrm{subject~to~}}

\def\R{\mathbb{R}}

\def\BL{\bm{L}}

\def\BS{\bm{S}}
\def\BC{\bm{C}}
\def\BU{\bm{U}}
\def\BR{\bm{R}}
\def\BV{\bm{V}}

\def\BX{\bm{X}}
\def\BA{\bm{A}}

\def\BQ{\bm{Q}}

\def\cE{\mathcal{E}}
\def\cR{\mathcal{R}}
\def\cX{\mathcal{X}}
\def\cL{\mathcal{L}}
\def\cS{\mathcal{S}}

\def\cO{\mathcal{O}}
\def\fro{\mathrm{F}}
\newcommand{\tens}[1]{\mathcal{#1}}

\newcommand{\alg}{RTCUR}
\newcommand{\algf}{\alg-F}
\newcommand{\algr}{\alg-R}

\newcommand{\supp}{{\rm supp\,}}

\definecolor{OliveGreen}{rgb}{0.29,0.5,0.3}
\definecolor{Green1}{rgb}{0.3,0.7,0.3}

\definecolor{dncolor}{rgb}{0.85, 0.15, 0.65}

\crefname{Definition}{Definition}{Definitions}

\usepackage{amsopn}

\usepackage{xspace}
\usepackage{bold-extra}
\usepackage[most]{tcolorbox}

\colorlet{texcscolor}{blue!50!black}
\colorlet{texemcolor}{red!70!black}
\colorlet{texpreamble}{red!70!black}
\colorlet{codebackground}{black!25!white!25}

\lstdefinestyle{siamlatex}{%
	style=tcblatex,
	texcsstyle=*\color{texcscolor},
	texcsstyle=[2]\color{texemcolor},
	keywordstyle=[2]\color{texemcolor},
	moretexcs={cref,Cref,maketitle,mathcal,text,headers,email,url},
}

\tcbset{%
	colframe=black!75!white!75,
	coltitle=white,
	colback=codebackground, 
	colbacklower=white, 
	fonttitle=\bfseries,
	arc=0pt,outer arc=0pt,
	top=1pt,bottom=1pt,left=1mm,right=1mm,middle=1mm,boxsep=1mm,
	leftrule=0.3mm,rightrule=0.3mm,toprule=0.3mm,bottomrule=0.3mm,
	listing options={style=siamlatex}
}

\newtcblisting[use counter=example]{example}[2][]{%
	title={Example~\thetcbcounter: #2},#1}

\newtcbinputlisting[use counter=example]{\examplefile}[3][]{%
	title={Example~\thetcbcounter: #2},listing file={#3},#1}

	\title{Robust Tensor CUR Decompositions: Rapid Low-Tucker-Rank\\Tensor Recovery with Sparse Corruptions\thanks{Received by the editors XXXX, 2023; accepted for publication (in revised form) XXXX XX, 2023; published electronically XXXX XX, 2023. 
			This paper is an extension of work originally presented in ICCVW \cite{cai2021fast}. 
			\funding{The authors contributed equally. This work was partially supported by NSF DMS 2011140, NSF DMS 2108479, NSF DMS 2304489, and AMS Simons Travel Grant.}}}
	
	\author{
		HanQin Cai\thanks{Department of Statistics and Data Science and Department of Computer Science, University of Central Florida, Orlando, FL 32816, USA. (email: \href{mailto:hqcai@ucf.edu}{hqcai@ucf.edu})}
		\and
		Zehan Chao\thanks{Department of  Mathematics, University of California  Los Angeles, Los Angeles, CA 90095, USA. (email:   \href{mailto:zchao@math.ucla.edu}{zchao@math.ucla.edu}) }
		\and
		Longxiu Huang\thanks{Department of Computational Mathematics, Science and Engineering and Department of Mathematics, Michigan State University, East Lansing, MI 48823, USA. (email:  \href{mailto:huangl3@msu.edu}{huangl3@msu.edu})}
		\and Deanna Needell\thanks{Department of  Mathematics, University of California Los Angeles,  Los Angeles,  CA 90095, USA. (email: \href{mailto:deanna@math.ucla.edu}{deanna@math.ucla.edu})} 
	}
	
	\headers{Robust Tensor CUR Decompositions}{H.Q.~Cai, Z.~Chao, L.~Huang, and D.~Needell }

\ifpdf
\hypersetup{pdftitle={Robust Tensor CUR Decompositions: Rapid Low-Tucker-Rank Tensor Recovery with Sparse Corruptions}}
\fi

\begin{document}
	\maketitle

	\begin{abstract}
		We study the tensor robust principal component analysis (TRPCA) problem, a tensorial extension of matrix robust principal component analysis (RPCA), that aims to split the given tensor into an underlying low-rank component and a sparse outlier component. This work proposes a fast algorithm, called \textbf{R}obust \textbf{T}ensor \textbf{CUR} Decompositions (\alg), for large-scale non-convex TRPCA problems under the Tucker rank setting. \alg\ is developed within a framework of alternating projections that projects between the set of low-rank tensors and the set of sparse tensors. We utilize the recently developed tensor CUR decomposition to substantially reduce the computational complexity in each projection. In addition, we develop four variants of \alg\ for different application settings. We demonstrate the effectiveness and computational advantages of \alg\ against state-of-the-art methods on both synthetic and real-world datasets.
	\end{abstract}
	

	\section{Introduction}\label{sec:intro}
	
	In our real world, high-dimensional data, such as images, videos, and DNA microarrays, often reside approximately on low-dimensional manifolds  \cite{Lin2017}. This association between low-dimensional manifolds and high-dimensional data has led mathematicians to develop various dimension reduction methods, such as principal component analysis (PCA) \cite{abdi2010principal} and non-negative matrix factorization (NMF) \cite{lee1999learning} under the low-rank assumption. It becomes increasingly important to study the low-rank structures for data in many fields such as  image processing \cite{li2018tucker, liu2022characterizing, chao2021hosvd}, video processing \cite{zhou2017tensor, lu2016tensor}, text analysis \cite{luo2015subgraph, chambua2018tensor}, and recommendation system \cite{zheng2016topic, song2017based}. 
	In reality, many higher-order data are naturally represented by tensors \cite{LU2011survey, zhou2013tensor, liu2012tensor}, which are  
	higher-order extensions of matrices. The tensor-related tasks such as tensor compression and tensor recovery usually involve finding a low-rank  
	structure from a given tensor. 
	
	As in the matrix setting, but to an even greater extent, principal component analysis is one of the most widely used methods for such dimension reduction tasks.  
	However, standard PCA is over-sensitive to extreme outliers \cite{candes2007sparsity}. To overcome this weakness, robust principal component analysis (RPCA) has been proposed to tolerate the sparse outliers in data analysis \cite{candes2011robust}. In particular, RPCA aims to reconstruct a low-rank matrix $\BL^\star$ and a sparse outlier matrix $\BS^\star$ from the corrupted observation 
	\begin{equation} \label{eq:rpca}
		\BX=\BL^\star+\BS^\star.
	\end{equation} 
	Existing studies on RPCA seek to find $\BL^\star$ and $\BS^\star$ by solving the following non-convex problem \cite{candes2011robust}, and we use $\BL$ and $\BS$ to denote the outcomes:
	\begin{equation}
		\begin{split}
			\minimize_{\BL,\BS} ~~&\|\BX - \BL - \BS\|_\fro \\ 
			\subject ~&~ \BL\textnormal{ is low-rank and } \BS \textnormal{ is sparse}.
			\label{eq:PCA}
		\end{split}
	\end{equation} 
	The term \textit{low-rank} refers to the constraint that the rank of $\BL$ is much smaller than its size, and the term \textit{sparse} refers to the restriction on the number of non-zero entries in $\BS$ (for example, one possible restriction is to allow each column and row of $\BS$ contain at most $10\%$ non-zero entries). This RPCA model has been widely studied \cite{netrapalli2014non-convex,bouwmans2018applications,cai2019accelerated,cai2021accelerated,cai2020rapid,cai2021robust,cai2021learned,hamm2022RieCUR,cai2022structured} and applied to many applications, e.g., 
	face modeling  \cite{wright2008robust}, 
	feature identification \cite{hu2019dstpca}, 
	and video background subtraction \cite{li2004statistical}.
	However, the original RPCA method can only handle $2$-mode arrays (i.e., matrices), while real-world data are often more naturally represented by higher-dimensional arrays (i.e., tensors). For instance, in the application of video background subtraction, a color video is automatically a $4$-mode tensor (height, width, frame, and color). To apply RPCA to tensor data, one has to unfold the original tensor into a matrix along some specific mode(s). 
	Although the upper bound of the unfolded matrix rank depends on the original tensor rank, the exact rank of the unfolded matrix remains unclear in some tensor rank settings. 
	In addition, we seek methods that utilize the structural information of tensor rather than ignoring the information. Therefore, it is important to generalize the standard RPCA to tensor settings. This task is called Tensor Robust Principal Component Analysis (TRPCA) \cite{lu2019tensor}. Moving from matrix PCA to the tensor setting could be challenging because some standard results known in the matrix case may not be generalized to tensors smoothly \cite{chen2009tensor}. For example, the rank of a matrix is well and uniquely defined, but researchers have proposed several definitions of tensor rank, such as Tucker rank \cite{tucker1966}, CP rank \cite{carroll1970analysis}, and Tubal rank \cite{zhang2014novel}.

	\subsection{Notation and Definitions}\label{SEC:Prelim}
	
	A tensor is a multi-dimensional array and its number of    dimensions 
	is called the \emph{order} or \emph{mode}. The space of real tensors of order $n$ and of size $(d_1,\cdots, d_n)$ is denoted as $\mathbb{R}^{d_1\times  \cdots\times d_n}$. 
	In this section, we first bring in the tensor-related notation and review some basic tensor properties, which will be used throughout the rest of the paper. 
	We denote tensors, matrices, vectors, and scalars in different typeface for clarity. More specifically, calligraphic capital letters (e.g., $\cX$) are used for tensors,  bold capital letters (e.g., $\BX$) are used for matrices, bold lower case letters (e.g., $\bm{x}$) for vectors, and regular letters (e.g., $x$) for scalars. We use $\BX(I,:)$ and $\BX(:,J)$ to denote the row and column submatrices of $\BX$ with index sets $I$ and $J$, respectively. $\mathcal{X}(I_1,\cdots,I_n)$ denotes the subtensor of $\mathcal{X}$ with index sets $I_k$ at $k$-th mode. A single element in a tensor is indexed as $\mathcal{X}_{i_1,\cdots, i_n}$. Moreover,
	\begin{equation*}
		\left\|\mathcal{X}\right\|_{\infty}=\max_{i_1,\cdots,i_n}|\mathcal{X}_{i_1,\cdots,i_n}| \text{ and } \left\|\mathcal{X}\right\|_\fro=\sqrt{\sum_{i_1,\cdots, i_n}\mathcal{X}_{i_1,\cdots,i_n}^2}
	\end{equation*}
	denote the max magnitude and Frobenius norm of a tensor, respectively. $\BX^\dagger$ denotes the Moore-Penrose pseudo-inverse of a matrix.
	The set of the first $d$ natural numbers is denoted by $[d]:=\{1,\cdots,d\}$.

	\begin{definition}[\textbf{Tensor matricization/unfolding}]
		An $n$-mode tensor $\mathcal{X}$ can be matricized, or reshaped into a matrix, in $n$ ways by unfolding it along each of the $n$ modes.  
		The mode-$k$ matricization/unfolding of tensor $\mathcal{X}\in\mathbb{R}^{d_1\times  \cdots\times d_n}$ is the matrix denoted by
		\begin{equation}
			\mathcal{X}_{(k)}\in\mathbb{R}^{d_k\times\prod_{j\neq k}d_j}
		\end{equation}
		whose columns are composed of all the vectors obtained from $\mathcal{X}$ by fixing all indices except for the $k$-th dimension.  The mapping $\mathcal{X}\mapsto \mathcal{X}_{(k)}$ is called the mode-$k$ unfolding operator.
	\end{definition}
	
	\begin{definition}[\textbf{Mode-$k$ product}]~
		Let $\mathcal{X}\in\mathbb{R}^{d_1\times \cdots\times d_n}$ and $\BA\in\mathbb{R}^{J\times d_k}$. The $k$-th mode multiplication between $\mathcal{X}$ and $\BA$ is denoted by $\mathcal{Y}=\mathcal{X}\times_k \BA$, with 
		\begin{equation}
			\mathcal{Y}_{i_1,\cdots,i_{k-1},j,i_{k+1},\cdots,i_{n}}=\sum_{s=1}^{d_k}\mathcal{X}_{i_1,\cdots,i_{k-1},s,i_{k+1},\cdots,i_{n}}\BA_{j,s}.
		\end{equation}
		This can be written as a matrix product by noting that $\mathcal{Y}_{(k)}=\BA \mathcal{X}_{(k)}$. If we have multiple tensor matrix products from different modes, we  use the notation $\mathcal{X}\times_{i=t}^{s} \BA_i$ to denote the product $\mathcal{X}\times_{t}\BA_{t}\times_{t+1}\cdots\times_{s}\BA_{s}$. We also use `tensor-matrix product' to name this operation throughout our paper. 
	\end{definition}
	
	In tensor analysis, Tucker rank \cite{HF1928} (also known as multilinear rank) is one of the most essential tensor ranks related to subspace estimation, compression, and dimensionality reduction \cite{rabanser2017introduction}.
	
	\begin{definition}[\textbf{Tensor Tucker rank and Tucker decomposition}] 
		The Tucker decomposition of tensor $\mathcal{X}$ is defined as an approximation of a  core tensor $\mathcal{C}$
		multiplied by $n$ factor matrices $\BA_k$ 
		(whose columns are usually orthonormal)  along each mode, such that
		\begin{equation}\label{eqn:Tucker_D}
			\mathcal{X}\approx\mathcal{C}\times_{i=1}^n \BA_i.
		\end{equation}
		If \eqref{eqn:Tucker_D} becomes an equation and $\mathcal{C}\in\R^{r_1\times\cdots\times r_n}$, 
		then we say this decomposition is an exact Tucker decomposition of $\mathcal{X}$.
	\end{definition}
	
	Note that higher-order singular value decomposition (HOSVD) \cite{de2000best} is a specific orthogonal Tucker decomposition which is popularly used in the literature.

	\subsection{Related Work: Tensor CUR Decompostions} \label{sec:2cur}
\begin{figure}[th]
		\centering
		\begin{minipage}{0.45\linewidth}
			\centering
			\includegraphics[width=\linewidth]{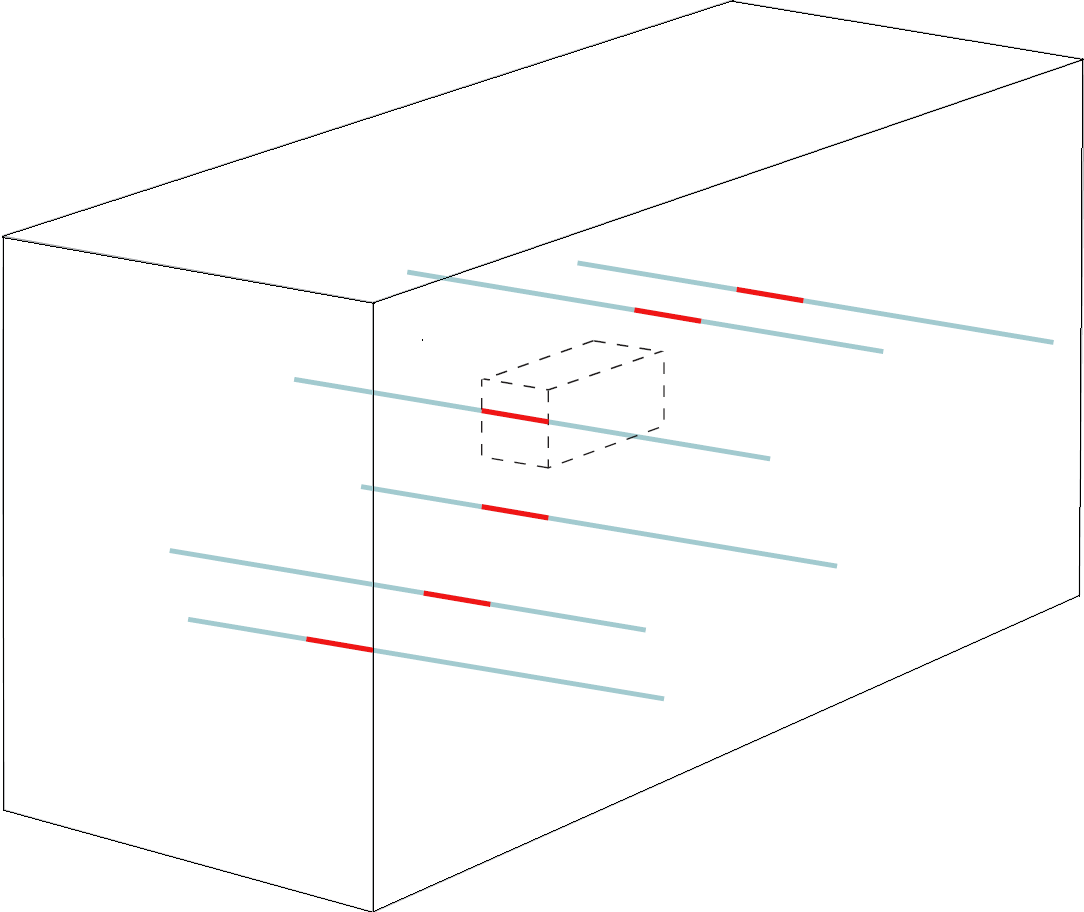}
			\caption{(\cite{cai2021mode}). 
				Illustration of the Fiber CUR Decomposition of \Cref{thm: CUR Char} in which $J_i$ is not necessarily related to $I_i$.  The lines correspond to rows of $\BC_2$, and red indices correspond to rows of $\BU_2$.  Note that the lines may (but do not have to) pass through the core subtensor $\cR$ outlined by dotted lines.  For the figure's clarity, we do not show fibers in $\BC_1$ and $\BC_3$.}
			\label{FIG:TensorCURFiber}
		\end{minipage}
		\hfill
		\begin{minipage}{0.45\linewidth}
			\centering\includegraphics[width=\linewidth]{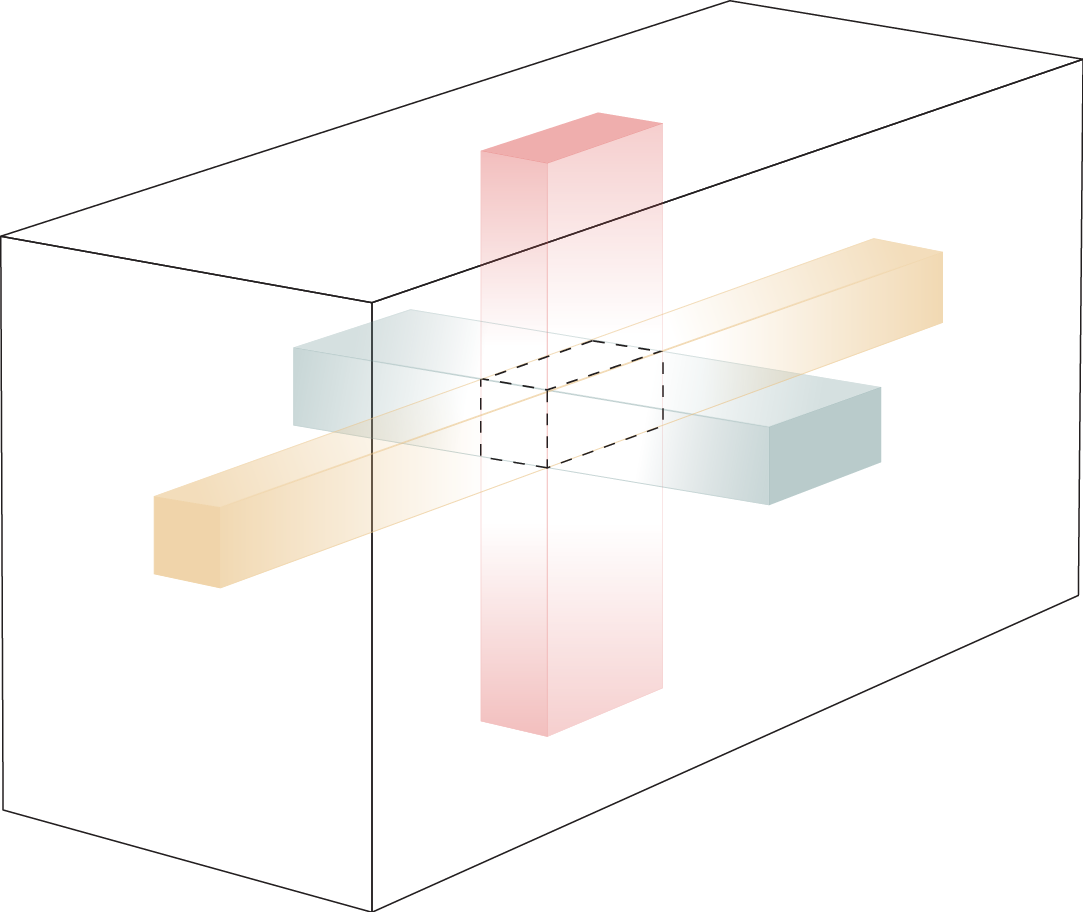}
			\caption{(\cite{cai2021mode}). Illustration of Chidori CUR decomposition of a 3-mode tensor in the case when the indices $I_i$ are each an interval and $J_i=\otimes_{j\neq i} I_j$ (see \Cref{thm: CUR Char}). The matrix $\BC_1$ is obtained by unfolding the red subtensor along mode 1, $\BC_2$ by unfolding the green subtensor along mode 2, and $\BC_3$ by unfolding the yellow subtensor along mode 3.  The dotted line shows the boundaries of $\mathcal{R}$. In this case $\BU_i=\mathcal{R}_{(i)}$ for all $i$.
			}\label{FIG:TensorCUR}
		\end{minipage}
		\vspace{-0.1in}
	\end{figure}

	Researchers have been actively studying CUR decompositions for matrices in recent years \cite{Goreinov,HH2020}. 
	For a matrix $\BX\in\mathbb{R}^{d_1\times d_2}$, let $\BC$ be  a submatrix consisting a subset of columns of $\BX$ with column indices $J$ , $\BR$ be a submatrix consisting a subset of rows of $\BX$ with row indices $I$, and $\BU=\BX(I,J)$. The theory of CUR decompositions states that $\BX=\BC\BU^\dagger \BR$ if $\rank(\BU)=\rank(\BX)$. 
	The first extension of CUR decompositions to tensors involved a single-mode unfolding of 3-mode tensors \cite{MMD2008}. 
	Later, \cite{caiafa2010generalizing} proposed 
	a different variant of tensor CUR that accounts for all modes. Recently, \cite{cai2021mode} dubbed these decompositions with more descriptive monikers, namely Fiber and Chidori CUR decompositions. In this paper, we will employ both  Fiber CUR decomposition and Chidori CUR decomposition (see \Cref{FIG:TensorCUR,FIG:TensorCURFiber} for illustration) to accelerate an essential step in the proposed algorithm. 
	We state the Fiber CUR and Chidori CUR decomposition characterizations below for the reader's convenience.
	
	\begin{theorem}[{\cite[Theorem~3.3]{cai2021mode}}]\label{thm: CUR Char}
		Let $\mathcal{A}\in\mathbb{R}^{d_1\times\cdots\times d_n}$ with Tucker rank $(r_1,\dots,r_n)$. Let $I_i\subseteq [d_i]$ and $J_i\subseteq[\prod_{j\neq i}d_j]$. Set $\mathcal{R}=\mathcal{A}(I_1,\cdots,I_n)$, $\BC_i=\mathcal{A}_{(i)}(:,J_i)$ and $\BU_i=\BC_i(I_i,:)$. 
		Then the following statements are equivalent:
		\begin{enumerate}[label=(\roman*),leftmargin=.5in]
			\item  \label{CUR Char:item2} $\mathcal{A}=\mathcal{R}\times_{i=1}^{n}(\BC_i\BU_i^\dagger)$,
			\item  \label{CUR Char:item1} $\rank(\BU_i)=r_i$,
			\item  \label{CUR Char:item3} $\rank(\BC_i)=r_i$ for all $i$ and the Tucker rank of $\mathcal{R}$ is $(r_1,\cdots,r_n)$.
		\end{enumerate}
	\end{theorem}
	\begin{remark}
		In particular, when $J_i$  
		are sampled independently from $I_i$, \Cref{thm: CUR Char}\ref{CUR Char:item2} is called Fiber CUR decomposition. When $J_i=\otimes_{j\neq i} I_j$, \Cref{thm: CUR Char}\ref{CUR Char:item2} is called  Chidori CUR decomposition.
		\label{rmk:fiberchidori}
	\end{remark}

	In addition, according to \cite[Corollary~5.2]{hamm2020stability}, if one uniformly samples indices $I_i$ and $J_i$ with size $|I_i|= \cO(r_i\log(d_i))$ and $|J_i|=\cO\left(r_i\log(\prod_{j\neq i}d_j)\right)$, then  $\rank(\BU_i)=r_i$ holds for all $i$ with high probability under some mild assumptions.
	Thus, the tensor CUR decomposition holds and its computational complexity is dominated by computing the pseudo-inverse of $\BU_i$.  Given the dimension of $\BU_i$, the computational complexity of the pseudo-inverse of $\BU_i$ with Fiber sampling is $\cO\left((n-1)r^3\log^2d\right)$, thus  Fiber CUR decomposition costs $\cO\left(nr^3\log^2d\right)$.\footnote{For notation simplicity, we assume the tensor has the same $d$ and $r$ along each mode when we discuss complexities. All $\log$ operators used in this paper stand for natural logarithms.} The Chidori CUR decomposition has a slightly larger $|J_i|$, which is $\prod_{j\neq i}r_i\log(d_i) = \cO((r\log d)^{n-1})$, thus the decomposition cost $\cO\left(r^{n+1}\log^nd\right)$. By contrast, the computational complexity of HOSVD is at least $\cO(rd^n)$.

	\subsection{Related Work: Tensor Robust Principal Component Analysis}
	There is a long list of studies on RPCA \cite{vaswani2018static} and low-rank tensor approximation \cite{LU2011survey}, so we refer readers to those two review articles for the aforementioned topics and focus on TRPCA works in this section. Consider a given tensor $\cX$ that can represent a hypergraph network or a multi-dimensional observation \cite{cai2022generalized};
	the general assumption of TRPCA is that $\cX$ can be decomposed as the sum of two tensors:
	\begin{equation}
		\cX =\cL^\star+\cS^\star,
	\end{equation}
	where $\cL^\star\in\mathbb{R}^{d_1\times\cdots\times  d_n}$ is the underlying low-rank tensor and $\cS^\star\in\mathbb{R}^{d_1\times\cdots\times  d_n}$ is the underlying sparse tensor. Compared to the exact low-rank tensor models, TRPCA model contains an additional sparse tensor $\cS$, which accounts for potential model outliers and hence more stable with sparse noise.
	Different from the well-defined matrix rank, there exist various definitions of tensor decompositions that lead to various versions of tensor rank, and lead to different versions of robust tensor decompositions. For example, \cite{liu2018improved,lu2019tensor,mohammadpour2022randomized} formulate TRPCA as a convex optimization model  based on the tubal rank \cite{zhang2014novel}. 
	
	Based on the Tucker rank, we aim to solve the non-convex optimization problem in this work:
	\begin{equation}\label{eq:objective}
		\begin{split}
			\minimize_{\cL,\cS} ~~&\|\cX-\cL-\cS\|_\fro \cr
			\subject ~&~\cL \textnormal{ is low-Tucker-rank and } \cS \textnormal{ is sparse}.
		\end{split}
	\end{equation}
	
	Researchers have developed different optimization methods to solve \eqref{eq:objective} \cite{dong2022fast,cai2022generalized,sofuoglu2018two,hu2020robust}. For example, the work in \cite{cai2022generalized} integrated the Riemannian gradient descent (RGD) and gradient pruning methods to develop a linearly convergent algorithm for \eqref{eq:objective}. This RGD algorithm will also serve as a guideline approach in our experiments. However, one of the major challenges in solving the Tucker rank based TRPCA problem is the high computational cost for computing the Tucker decomposition. If $\cL^\star$ is rank-$(r_1,\cdots,r_n)$, the existing methods, e.g., \cite{huang2014provable,gu2014robust,sofuoglu2018two,hu2020robust,cai2022generalized}, have computational complexity at least $\cO(nd^n r)$---they are thus computationally challenging in large-scale problems. Thus, it is necessary to develop a highly efficient TRPCA algorithm for time-intensive applications.

	\subsection{Contributions}
	In this work, we consider the TRPCA problem under the Tucker rank setting. Our main contributions are three-fold:
	\begin{enumerate}
		\item We provide  theoretical evidence supporting the generality of Tensor Robust Principal Component Analysis (TRPCA) model over the  matrix robust PCA  model obtained from the unfolded tensor (see \Cref{sec:cha-SOT}). That is, TRPCA requires a much weaker sparsity condition on the outlier component. Our theoretical finds will be empirically verified later in the numerical section.
		\item We propose a novel non-convex approach, coined \textbf{R}obust \textbf{T}ensor \textbf{CUR} Decompositions (\alg), for large-scale\footnote{In our context, `large-scale' refers to large $d$.} TRPCA problems (see \Cref{sec:proposed approach}). \alg\ uses a framework of alternating projections and employs a novel mode-wise tensor decomposition \cite{cai2021mode} for fast low-rank tensor approximation. We present four variants of \alg\ with different sampling strategies (see \Cref{sec:4vars} for the details about sampling strategies). The computational complexity of \alg\ is as low as $\cO(n^2d r^2\log^2d)$ or $\cO(nd r^n\log^nd)$ flops, depending on the sampling strategy, for an input $n$-mode tensor of size\footnote{For notation simplicity, we assume the tensor has the same $d$ and $r$ along each mode when we discuss complexities. All $\log$ operators used in this paper stand for natural logarithms.} $\mathbb{R}^{d\times  \cdots\times d}$ with Tucker rank ($r,\dots,r)$. Both computational complexities are substantially lower than the state-of-the-art TRPCA methods. For instance, two state-of-the-art methods \cite{lu2019tensor,lu2016tensor} based on tensor singular value decomposition have computational costs at least $\cO(nd^nr)$ flops.
		\item We verify the empirical advantages of \alg\ with synthetic datasets and three real-world applications (see \Cref{sec:experiment}), including robust face modeling, video background subtraction, and network clustering. We show that \alg\ has not only speed efficiency but also superior robustness compared to the state-of-the-arts. 
		In particular, we provide certain outlier patterns that can be detected by \alg\ but fails all matrix based methods (see \Cref{sec:face}). This further verifies our theoretical finds in \Cref{sec:cha-SOT}. 
	\end{enumerate}
	
	\section{Characterizations of Sparse Outlier Tensor} \label{sec:cha-SOT}
	Tensor RPCA problem can be also solved by using the RPCA method if we unfold the tensor into a matrix along a certain mode. To solve the tensor-to-matrix RPCA problem successfully, the unfolded outlier tensor must satisfy $\alpha$-M-sparsity, a commonly used assumption in matrix RPCA problem. We state the definition of $\alpha$-M-sparsity for matrix as follows. 
	\begin{definition}[$\alpha$-M-sparsity for matrix] \label{def:matrix sparse}
		$S\in\R^{d_1\times d_2}$ is $\alpha$-M-sparse if
		\[  
		\|S\bm{e}_i\|_0\leq \alpha d_1 \quad \textnormal{and} \quad \|\bm{e}_j^\top S\|_0\leq \alpha d_2
		\]
		for all $i=1,\cdots,d_2$ and $j=1,\cdots,d_1$.
	\end{definition}
	However, when solving TRPCA directly with tensor-based methods, it is more natural to generalize $\alpha$-M-sparsity to tensor setting. We consider the following sparsity condition for outlier tensor \cite{cai2022generalized}.
	{
		\begin{definition}[$\alpha$-T-sparsity for tensor] \label{def:tensor sparse}
			A tensor $\cS\in\mathbb{R}^{d_1\times\cdots\times d_n}$ is  $\alpha$-T-sparse if
			$$
			\|\mathcal{S}\times_j \bm{e}_{k_j,j}^\top\|_0\leq \alpha \prod_{i\neq j}d_i, 
			$$
			for all $k_j=1,\cdots,d_j$, where $ \{\bm{e}_{k_j,j}\}_{k_j=1}^{d_j}$ is the standard basis of $\R^{d_j}$. Thus, $\mathcal{S}\times_j \bm{e}_{k_j,j}^\top$ corresponds to the $k_j$-th slice of $\cS$ along mode $j$.
		\end{definition}
	}
	\begin{figure}[h]
		\centering
		\includegraphics[width=0.67\textwidth]{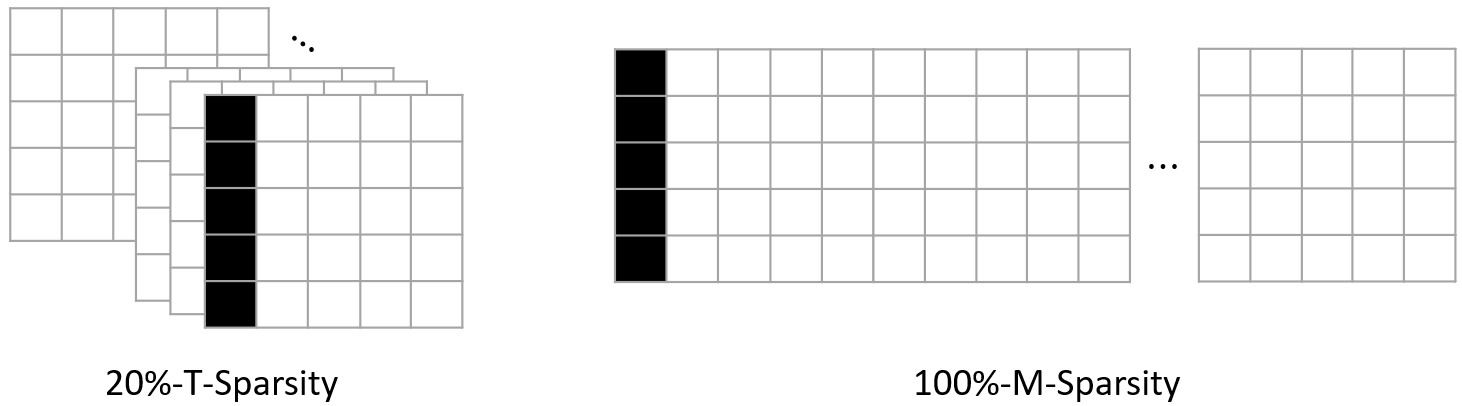}
		\caption{T-Sparsity vs.~M-Sparsity. A black box represents an outlier entry and a white box represents a good entry. The right-hand-side matrix is unfolded from the left-hand-side tensor.} \label{fig:TvsM}
		\vspace{-0.1in}
	\end{figure}
	We emphasize that \Cref{def:matrix sparse,def:tensor sparse} are not equivalent. In fact, from \Cref{fig:TvsM}, one can see that unfolding an outlier tensor along a certain mode can result in much worse sparsity in the unfolded matrix, with exact same outlier pattern. Note one can argue that unfolding alone a different mode may result in a not-worse M-sparsity. Although this is true, in practice the user usually does not have prior knowledge of outlier patterns, and thus cannot determine the most robust mode to unfold alone. Hence, we claim that \Cref{def:tensor sparse} is a much weaker condition than \Cref{def:matrix sparse}. This is further verified by the following theorem. 
	\begin{theorem}\label{asparse}
		Suppose 
		that an $n$-order tensor  $\mathcal{S} \in{\mathbb{R}^{d\times\dots\times d}}$ is generated according to the Bernoulli distribution with expectation $\frac{\alpha}{2}$, i.e.,  $\mathcal{S}_{i_1,\cdots,i_n}\sim\mathtt{Ber}(\frac{\alpha}{2})$, and $\alpha d>\frac{2n\log d}{\log 4-1}$.
		Then, $\cS$ is $\alpha$-T-sparse with probability at least 
		Moreover, by unfolding $\cS$ into a matrix $M\in\R^{d^{k}\times d^{n-k}}$ along  some  particular modes with $k\in [1,n-1]$, $M$ is $\alpha$-M-sparse with probability at least $1-d^{k-nd^{n-k-1}}-d^{n-k-nd^{k-1}}$. 
	\end{theorem}

	\begin{remark}\label{rmk: t v.s. m}
		It is evident that for large-scale tensor, i.e., $d\gg 1$ and $n>2$, we have $1-nd^{n-2}\ll \max\{k-nd^{n-k-1},n-k-nd^{k-1}\}$. Thus, the condition of $\alpha$-T-sparsity for $\cS$ can be satisfied with much better probability than that of $\alpha$-M-sparsity for the unfolders of $\cS$. 
	\end{remark}
	\begin{proof}[Proof of \Cref{asparse}]
		Since $\mathcal{S}_{i_1,\cdots,i_n}\sim \mathtt{Ber}(\frac{\alpha}{2})$, we have $\mathbb{E}(\mathcal{S}_{i_1,\dots,i_n})=\frac{\alpha}{2}$. Let's first consider the sparsity of one slice of $\mathcal{S}$: we use $\mathcal{S}_{(j)[k]}$ to denote the $k$-th slice along mode $j$.
		According to the multiplicative Chernoff bound, we have 
		\begin{equation}
			\mathbb{P}\left(\sum_{i_1,\cdots,i_{j-1},i_{j+1},\cdots,i_n} [\mathcal{S}_{(j)[k]}]_{i_1,\cdots,i_{j-1},i_{j+1},\cdots,i_n} \geq \alpha d^{n-1} \right) < \Big(\frac{e}{4}\Big)^{\frac{\alpha d^{n-1}}{2}}.
		\end{equation}
		Taking the sparsity of all slices along all modes into account,  then the probability  that $\mathcal{S}$ is $\alpha$-sparsity can be bounded by:
		\begin{align}
			\nonumber
			\mathbb{P}(\mathcal{S} \text{ is } \alpha\textnormal{-T-sparse}) & \geq \bigg(1-\Big(\frac{e}{4}\Big)^{\frac{\alpha d^{n-1}}{2}}\bigg)^{n d}  \geq
			1-nd\Big(\frac{e}{4}\Big)^{\frac{\alpha d^{n-1}}{2}}.
		\end{align}
		When 		$\alpha d>\frac{2n\log d}{\log4-1}$,
		we thus have:
		\begin{align}
			\mathbb{P}(\mathcal{S} \text{ is } \alpha\textnormal{-T-sparse}) \geq 1- nd^{1-nd^{n-2}}.
		\end{align}
		Similarly, if we unfolder $\mathcal{S}$ into a matrix $M\in\R^{d^{k}\times d^{n-k}}$, we have that
		\begin{align}
			\mathbb{P}(M \text{ is } \alpha\textnormal{-M-sparse}) \geq 1-d^{k-nd^{n-k-1}}-d^{n-k-nd^{k-1}}
		\end{align}
		provided that $ \frac{2n\log d}{\log4-1}$.
	\end{proof}
	
	\begin{remark} \label{rmk:uncertainty principle}
		The purpose of the sparsity condition is to ensure a well-defined robust tensor PCA problem, i.e., the low-rank and sparse tensors are separable. By matrix rank-sparsity uncertainty principle \cite{chandrasekaran2011rank}, $\alpha$-M-sparsity ensures that the matrix problem unfolded from tensor is well-defined. However, with the newly developed tensor Tucker-rank-sparsity uncertainty principle \cite[Proposition 2]{xia2021statistically}, the more relaxed $\alpha$-T-sparsity is enough for well-defined robust tensor PCA problems. Moreover, to solve the problems with $\alpha$-T-sparsity right, an algorithm based directly on tensor structures is needed, like the one that will be proposed in the next section. The unfolded matrix based algorithms will still require the more restricted $\alpha$-M-sparsity condition, thus more likely to fail.   
	\end{remark}

	\section{Proposed Approach} \label{sec:proposed approach}
	In this section, we propose an efficient approach, called \textbf{R}obust \textbf{T}ensor \textbf{CUR} Decompositions (\alg), for the non-convex TRPCA problem \eqref{eq:objective}. \alg\ is developed in a framework of alternating projections: (I) First,  we project $\cX-\cL^{(k)}$ onto the space of sparse tensors to update the estimate of outliers (i.e., $\cS^{(k+1)}$); (II) then we project the less corrupted data $\cX-\cS^{(k+1)}$ onto the space of low-Tucker-rank tensors to update the estimate (i.e., $\cL^{(k+1)}$). In our algorithm, the key to acceleration is using the tensor CUR decomposition for inexact low-Tucker-rank tensor approximation in Step (II), which is proved 
	to be much more efficient than the standard HOSVD \cite{cai2021mode}, in terms of computational complexity. Consequently, in Step~(I), this inexact approximation allows us to estimate only the outliers in the smaller subtensors and submatrices involved in the tensor CUR decomposition. 
	\alg\ is summarized in \Cref{alg:TRCUR}. Notice that there are two variants of tensor CUR decompositions which will result in different $J_i$ (see \Cref{rmk:fiberchidori}), but the steps of \Cref{alg:TRCUR} will remain the same. Therefore, we will not distinguish the two decomposition methods in \Cref{sec:updateS,sec:updateL} when discussing the details of Step~(I) and (II). We will then show the computational complexity for \Cref{alg:TRCUR} with both Fiber and Chidori CUR decompositions in \Cref{sec:cc}.
	\subsection{Step~(I): Update Sparse Component $\cS$} 
	\label{sec:updateS}
	We consider the simple yet effective hard thresholding operator $\mathrm{HT}_\zeta$ for outlier estimation. The operator is defined as:
	
	\begin{equation}
		(\mathrm{HT}_{\zeta}(\cX))_{i_1,\cdots,i_n} =
		\begin{cases}
			\cX_{i_1,\cdots,i_n}, & \quad|\cX_{i_1,\cdots,i_n}| >\zeta;\\
			0,  & \quad\mbox{otherwise.}
		\end{cases}
	\end{equation}
	
	\begin{algorithm}[t]
		\caption{\textbf{R}obust \textbf{T}ensor  \textbf{CUR} Decompositions (\alg)}
		\label{alg:TRCUR}
		\begin{algorithmic}[1]
			\State \textbf{Input: }{$\cX=\cL^\star+\cS^\star\in\mathbb{R}^{d_1\times \cdots \times d_n}$: observed tensor; $(r_1, \cdots, r_n)$: underlying Tucker rank of $\cL^\star$; $\varepsilon$: targeted precision; $\zeta^{(0)},\gamma$: thresholding parameters; $\{|I_i|\}_{i=1}^n,\{|J_i|\}_{i=1}^n$: cardinalities for sample indices.} \qquad\quad {\color{OliveGreen}// $J_i$ is defined differently for different sampling strategies. See \Cref{sec:4vars} for details about $J_i$ and sampling strategies.}
			\State \textbf{Initialization:} $\cL^{(0)} = \bm{0},\cS^{(0)} = \bm{0},k=0$ 
			\State Uniformly sample the indices $\{I_i\}_{i=1}^n, \{J_i\}_{i=1}^n$ 
			\While {$e^{(k)} > \varepsilon$} \qquad\qquad\qquad {\color{OliveGreen}// $e^{(k)}$ is defined in \eqref{eq:rel_err}}
			\State 
			{\color{OliveGreen} (Optional)} Resample the indices $\{I_i\}_{i=1}^{n}, \{J_i\}_{i=1}^{n}$ 
			\State {\color{OliveGreen} // Step (I): Updating $\cS$} 
			\State $\zeta^{(k+1)} = \gamma\cdot \zeta^{(k)}$ 
			\State $\cS^{(k+1)} = \mathrm{HT}_{\zeta^{(k+1)}}(\tens{X}-\tens{L}^{(k)})$  
			
			\State {\color{OliveGreen} // Step (II): Updating $\cL$} 
			\State $\cR^{(k+1)}= (\cX-\cS^{(k+1)})(I_1,\cdots,I_n)$ 
			\For{$i = 1,\cdots,n$}
			\State $\BC_i^{(k+1)} = (\tens{X}-\cS^{(k+1)})_{(i)}(:,J_i)$ 
			\State $\BU_i^{(k+1)} = \mathrm{SVD}_{r_i}(\BC_i^{(k+1)}(I_i,:))$ 
			\EndFor
			\State $\cL^{(k+1)} = \cR^{(k+1)}\times_{i=1}^{n} \BC_i^{(k+1)}\left(\BU_{i}^{(k+1)}\right)^{\dagger}$ 
			\State $k = k+1$ 
			\EndWhile
			\State \textbf{Output: }{$\cR^{(k)}, \BC_i^{(k)},\BU_i^{(k)}$ for $i=1,\cdots,n$: the estimates of  the tensor   CUR decomposition of $\cL_{\star}$.
			}
		\end{algorithmic}
	\end{algorithm}

	As shown in \cite{cai2019accelerated,cai2020rapid,netrapalli2014non-convex}, with a properly chosen thresholding value, $\mathrm{HT}_{\zeta}$ is effectively a projection operator onto the support of $\cS^\star$. More specifically, we update
	\begin{equation} \label{eq:S=HT(X-L)}
		\tens{S}^{(k+1)} = \mathrm{HT}_{\zeta^{(k+1)}}(\cX-\cL^{(k)}).
	\end{equation}
	If $\zeta^{(k+1)}=\|\cL^\star-\cL^{(k)}\|_\infty$ 
	is chosen, then we have $\supp(\cS^{(k+1)})\subseteq\supp(\cS^\star)$ and $\|\cS^\star-\cS^{(k+1)}\|_\infty\leq 2\|\cL^\star-\cL^{(k)}\|_\infty$. Empirically, we find that iteratively decaying thresholding values
	\begin{equation}
		\zeta^{(k+1)} = \gamma \cdot \zeta^{(k)}
	\end{equation}
	provide superb performance with carefully tuned $\gamma$ and $\zeta^{(0)}$. Note that a favorable choice of $\zeta^{(0)}$ is $\|\cL^\star\|_\infty$, which can be easily estimated in many applications. The decay factor $\gamma\in(0,1)$ should be tuned according to the level of difficulty of the TRPCA problem, e.g., those problems with higher rank, more dense outliers, or large condition numbers are considered to be harder.
	For successful reconstruction of $\cL^\star$ and $\cS^\star$, the harder problems require larger $\gamma$. When applying \alg\ on both synthetic and real-world data, we observe that $\gamma\in[0.6,0.9]$ generally performs well. Since real-world data normally leads to more difficult problems, we fix $\gamma=0.7$ for the  synthetic experiment and $\gamma=0.8$ for the real-world data studies in \Cref{sec:experiment}. 
	
	\subsection{Step~(II): Update Low-Tucker-rank Component $\cL$}
	\label{sec:updateL}
	SVD is the most popular method for low-rank approximation under matrix settings since SVD gives the best rank-$r$ approximation of given matrix $\BX$, both with respect to the operator norm and to the Frobenius norm \cite{bergqvist2010higher}. Similarly, HOSVD has been the standard method for low-Tucker-rank approximation under tensor settings in many works \cite{de2000best,liu2012tensor,bergqvist2010higher,rajwade2012image}. However, the computational complexity of HOSVD is at least $\cO(rd^n)$; hence computing HOSVD is very expensive when the problem scale is large. 
	As highlighted in \cite[Sections~3.2 and 3.3]{cai2021mode}, tensor CUR decomposition can serve as an effective low-Tucker-rank approximation method, even with perturbations.  As such, we employ tensor CUR decomposition for accelerated inexact low-Tucker-rank tensor approximations. Namely, we update the estimate of the  low-Tucker-rank component $\cL$ by setting
	\begin{equation} \label{eq:L=CUR}
		\cL^{(k+1)} = \cR^{(k+1)}\times_{i=1}^{n} \BC_i^{(k+1)}\left(\BU_{i}^{(k+1)}\right)^{\dagger},
	\end{equation}
	where 
	\begin{equation} \label{eq:R_and_C_and_U}
		\begin{split}
			\cR^{(k+1)} &= (\cX-\cS^{(k+1)})(I_1,\cdots,I_n), \cr
			\BC_i^{(k+1)} &= (\cX-\cS^{(k+1)})_{(i)}(:,J_i), \cr
			\BU_i^{(k+1)} &= \mathrm{SVD}_{r_i}(\BC_i^{(k+1)}(I_i,:)).
		\end{split}
	\end{equation}

	\subsection{Computational Complexities} \label{sec:cc}
	As mentioned in \Cref{sec:2cur}, the complexity for computing a tensor CUR decomposition is much lower than HOSVD, and the dominating steps in \alg\ are the hard thresholding operator and the tensor/matrix multiplications. 
	For both Fiber and Chidori CUR decompositions, only the sampled subtensors and submatrices are required when computing \eqref{eq:R_and_C_and_U}. Thus, we merely need to estimate the outliers on these subtensors and submatrices, and \eqref{eq:S=HT(X-L)} should not be fully executed. Instead, we only compute
	\begin{equation} \label{eq:partial S update}
		\begin{split}
			\cS^{(k+1)}(I_1,\cdots,I_n) &= \mathrm{HT}_{\zeta^{(k+1)}}((\cX-\cL^{(k)})(I_1,\cdots,I_n)), \cr
			\cS^{(k+1)}_{(i)}(:,J_i) &= \mathrm{HT}_{\zeta^{(k+1)}}((\cX-\cL^{(k)})_{(i)}(:,J_i))
		\end{split}
	\end{equation}
	for all $i$. 
	Not only can we save the computational complexity on hard thresholding but also, much smaller subtensors of $\cL^{(k)}$ need to be formed in  \eqref{eq:partial S update}. We can form the required subtensors from the saved tensor CUR components, which is much cheaper than forming and saving the whole $\cL^{(k)}$.
	
	In particular,  for $\cX\in\mathbb{R}^{d\times \cdots\times d}$, $r_1=\cdots=r_n=r$ and $|I_1|=\cdots=|I_n|=\cO(r\log d)$, computing $\cL^{(k)}(I_1,\cdots,I_n)$ requires $n$ tensor-matrix product operations so the complexity for computing $\cL^{(k)}(I_1,\cdots,I_n)$ is $\cO(n(r\log d)^{n+1})$ flops for both Fiber and Chidori CUR decompositions. The complexity for computing $\cL^{(k)}_{(i)}(:,J_i)$ with Fiber CUR is different from the complexity of computing $\cL^{(k)}_{(i)}(:,J_i)$ with Chidori CUR. With Fiber CUR, we compute each fiber in $\cL^{(k)}_{(i)}(:,J_i)$ independently and each fiber takes $n$ tensor-matrix product operations. The first $n-1$ operations transform the $n$-mode core tensor $\cL^{(k)}(I_1,\cdots,I_n)$ into a $1$-mode tensor, which is a vector of length $\cO(r\log d)$, and the last operation transforms this vector into another vector of length $d$. Since there are $J_i = \cO(nr\log d)$ fibers in total, the complexity for computing $\cL^{(k)}_{(i)}(:,J_i)$ with Fiber CUR decomposition is $\cO(nr\log d((r\log d)^n+d r\log d))$ flops. With Chidori CUR, we compute $\cL^{(k)}_{(i)}(:,J_i)$ as a complete unit using $n$ tensor-matrix product operations. The first $n-1$ operations on the core tensor do not change its size, and the last operation changes the size of the $i$th mode to $d$. Therefore the complexity for computing $\cL^{(k)}_{(i)}(:,J_i)$ with Chidori CUR decomposition is $\cO(n(r\log d)^{n+1}+d(r\log d)^n)$ flops.
	
	Moreover, for time-saving purposes, we may avoid {computing} the Frobenius norm of the full tensor when computing the relative error for the stopping criterion. In \alg, we adjust the relative error formula to be  
	\begin{equation} \label{eq:rel_err}
		e^{(k)}=  \frac{\|\cE^{(k)}(I_1,\cdots,I_n)\|_\fro+\sum_{i=1}^n\|\cE^{(k)}_{(i)}(:,J_i)\|_\fro}{\|\cX(I_1,\cdots,I_n)\|_\fro+\sum_{i=1}^n\|\cX_{(i)}(:,J_i)\|_\fro},
	\end{equation}
	where $\cE^{(k)}=\cX-\cL^{(k)}-\cS^{(k)}$, 
	so that it does not use any extra subtensor or fiber but only those we already have. 
	We hereby summarize the computational complexity for each step from \Cref{alg:TRCUR} in \Cref{tab:complexity}.
	
	If we assume that the tensor size $d$ is comparable with or greater than $\cO((r\log d)^{n-1})$, we can conclude that the total computational complexity is $\cO(n^2 d r^2\log^2d$) flops for \alg\ with Fiber CUR decomposition and $\cO(nd r^n\log^nd$) flops for \alg\ with Chidori CUR decomposition. Otherwise, the computational complexity would be $\cO(n^2 r^{n+1}\log^{n+1}d$) flops for \alg\ with Fiber CUR, and the complexity for \alg\ with Chidori CUR remains unchanged. For all tensors tested in \Cref{sec:experiment}, the first case holds. Therefore, in \Cref{tab:complexity}, we highlighted $\cO(n^2 d r^2\log^2d$) and $\cO(nd r^n\log^nd$) as the dominating terms.
	
	\begin{table}[h]
		\caption{computational complexity for each step from \Cref{alg:TRCUR}. The complexity for computing $\cS(I_1,\cdots,I_n)$ and $\cR$ are the same as their size; the complexity for $\BC_i\BU_i^\dagger$ is introduced in \Cref{sec:2cur}; the complexity for computing $\cL$ and error term is introduced in \Cref{sec:cc}. The dominating terms are highlighted in bold.}
		\label{tab:complexity}
		\centering
		\vspace{-0.1in}
		\resizebox{\textwidth}{!}{  \begin{tabular}{c|cc}\toprule
				\textsc{Computational Complexity} & \textsc{Fiber Sampling} & \textsc{Chidori Sampling}\\\midrule
				Sparse subtensor $\cS(I_1,\cdots,I_n)$ or $\cR$ & $\cO(r^n\log^nd)$ & $\cO(r^n\log^nd)$\\
				All $\cS_{(i)}(:,J_i)$ or $\BC_i$ for $n$ modes & $\cO(n^2rd\log d )$ & $\cO(nr^{n-1}d\log^{n-1}d)$\\
				All $\BU_i^\dagger$ for $n$ modes & $\cO(n^2r^3\log^2d)$ & $\cO(nr^{n+1} \log^nd)$\\
				All $\BC_i\BU_i^\dagger$ for $n$ modes & $\cO(\pmb{n^2r^2d \log^2d})$ & $\cO(\pmb{nr^nd \log^nd})$\\
				Low-rank subtensor $\cL(I_1,\cdots,I_n)$ & $\cO(nr^{n+1}\log^{n+1}d)$ & $\cO(nr^{n+1}\log^{n+1}d)$\\
				All $\cL_{(i)}(:,J_i)$ for $n$ modes & 
				$\cO(n^2r^{n+1}\log^{n+1}d+\pmb{n^2r^2d \log^2d})$
				& $\cO(n^2r^{n+1}\log^{n+1}d+\pmb{nr^{n}d\log^nd})$\\
				Error term $\cE^{(k)}(I_1,\cdots,I_n)$ and $\cE^{(k)}_{(i)}(:,J_i)$ & $\cO(r^n\log^nd+n^2dr\log d)$ & $\cO(dr^{n-1}\log^{n-1}(d))$
				\\\bottomrule
			\end{tabular}
		}
	\end{table}

	\subsection{Four Variants of \alg} \label{sec:4vars}
	In \Cref{sec:2cur}, we discussed two versions of tensor CUR decomposition: Fiber CUR decomposition and Chidori CUR decomposition. Each of the decomposition methods could derive two slightly different \alg\ algorithms depending on if we fix sample indices through all iterations (see \Cref{alg:TRCUR}).  
	As a result of this, we obtain four variants in total. We give different suffixes for each variant of \alg\ algorithm: \alg-FF, \alg-FC, \alg-RF, and \alg-RC. We will showcase experimental results for all variants in \Cref{sec:experiment}. The first letter in the suffix indicate whether we fix the sample indices through all iterations: `F' stands for `fix', where the variant uses fixed sample indices through all iterations; `R' stands for `resampling',  where the variant resamples $\{I_i\}_{i=1}^n$ and $\{J_i\}_{i=1}^n$ in each iteration. The second letter indicate which type of CUR decomposition we use in \alg. `F' represents that \alg\ is derived from Fiber CUR decomposition and `C' stands for Chidori CUR. For Fiber CUR, the amount of fibers to be sampled refers to \cite[Corollary 5.2]{hamm2020stability}, i.e., $|I_i|= \upsilon r_i\log(d_i)$, $|J_i|=\upsilon r_i\log(\prod_{j\neq i}d_j)$ and $J_i$ is sampled independently from $I_i$. Here, $\upsilon$ denotes the sampling constant, a hyper-parameter that will be tuned in the experiments. For Chidori CUR, $|I_i|= \upsilon r_i\log(d_i)$ and $J_i=\otimes_{j\neq i} I_j$.  Of these four variants, \alg-FF requires minimal data accessibility and runs slightly faster than other variants. The resampling variants access more data and take some extra computing; for example, the denominator of \eqref{eq:rel_err} has to be recomputed per iteration. However, accessing more redundant data means resampling variants have better chances of correcting any ``unlucky'' sampling over the iterations. Thus, we expect resampling variants to have superior outlier tolerance than fixed sampling variants. And the fixed sampling variants have an efficiency advantage over the resampling variants under specific conditions (e.g., when re-accessing the data is expansive).

	The difference between Chidori variants and Fiber variants has similar properties: if we choose the same $\upsilon$ and let $|I_i|=  \upsilon r_i\log(d_i)$ for both Chidori and Fiber CUR described in \Cref{sec:2cur}, the Chidori variants generally access more tensor entries compared to the Fiber variants. Therefore, with the same sampling constant $\upsilon$, Chidori variants requires more computing time in each iteration. Nevertheless, Chidori variants can tolerate more dense outliers with these extra entries than Fiber variants.  
	We will further investigate their computational efficiency and practical performance in \Cref{sec:experiment}.
	\begin{remark}
		The tensor CUR decomposition represented in \Cref{thm: CUR Char}\ref{CUR Char:item2} is also in Tucker decomposition form. We can efficiently convert the tensor CUR decomposition to HOSVD with \Cref{alg:CUR to HOSVD} \cite{cai2021mode}. 
		\begin{algorithm}[t!]
			\caption{Conversion from CUR to HOSVD}
			\label{alg:CUR to HOSVD}
			\begin{algorithmic}[1]
				\State \textbf{Input: }{$\cR, \BC_i, \BU_i$ : CUR decomposition of the tensor $\tens{A}$}
				\State $\left[\BQ_{i}, \BR_{i}\right]=\operatorname{qr}\left(\BC_{i} \BU_{i}^{\dagger}\right) \text { for } i=1, \cdots, n$ 
				\State $\mathcal{T}_{1}=\mathcal{R} \times_{1} \BR_{1}   \times_{2} \cdots  \times_{n}\BR_{n}$ \State $\text {Compute HOSVD of } \mathcal{T}_{1} \text { to find } \mathcal{T}_{1}=\mathcal{T} \times_{1} \BV_{1} \times{ }_{2} \cdots \times_{n} \BV_{n}$
				
				\State \textbf{Output: }{$\llbracket\mathcal{T} ; \BQ_{1} \BV_{1}, \cdots, \BQ_{n} \BV_{n}\rrbracket$: HOSVD decomposition of $\tens{A}$}
			\end{algorithmic}
		\end{algorithm}
		In contrast, converting HOSVD to a tensor CUR decomposition is not as straightforward.
	\end{remark}

	\section{Numerical Experiments}
	\label{sec:experiment}
	\begin{figure}[ht]
		\centering
		\includegraphics[width=0.24\linewidth]{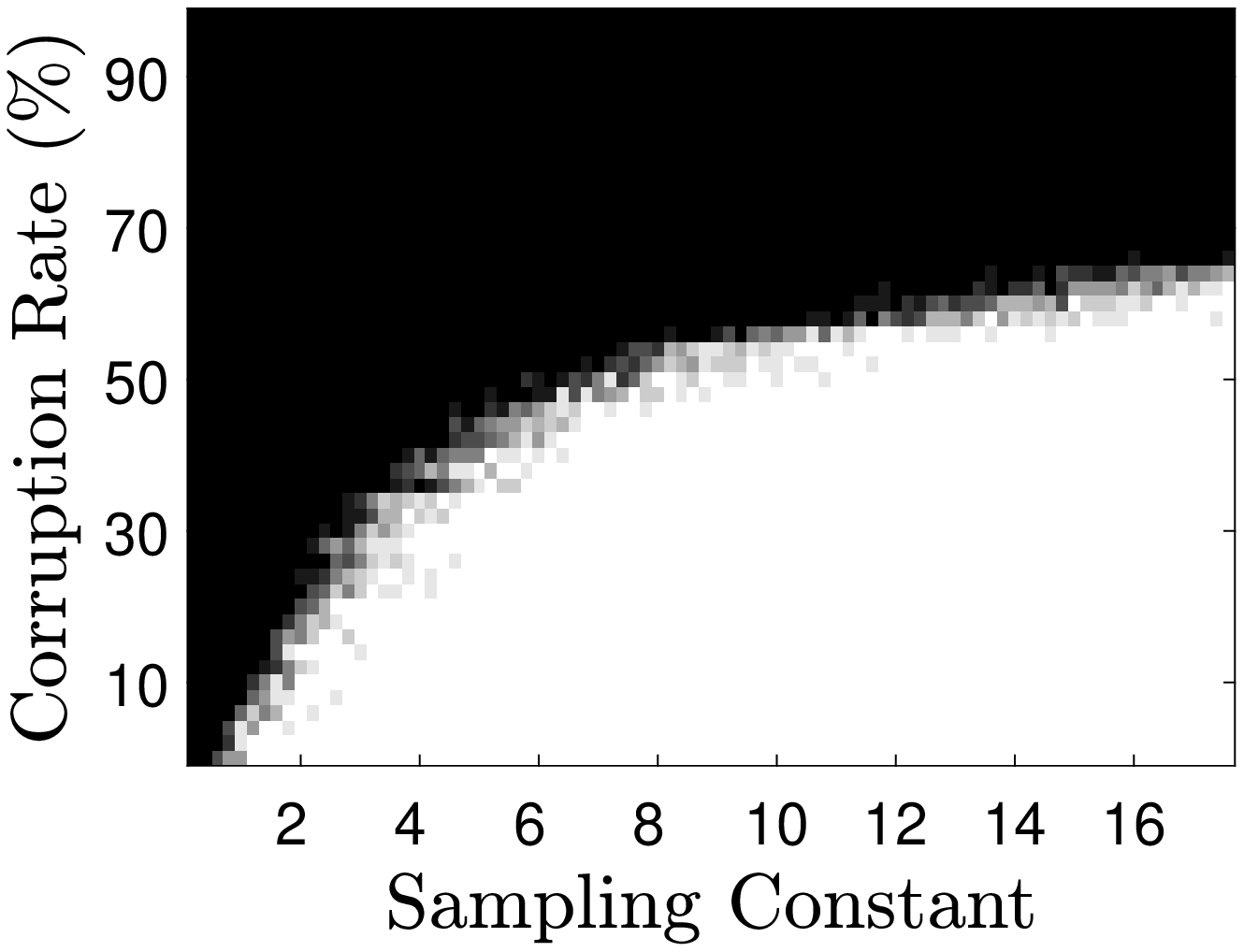}
		\hfill
		\includegraphics[width=0.24\linewidth]{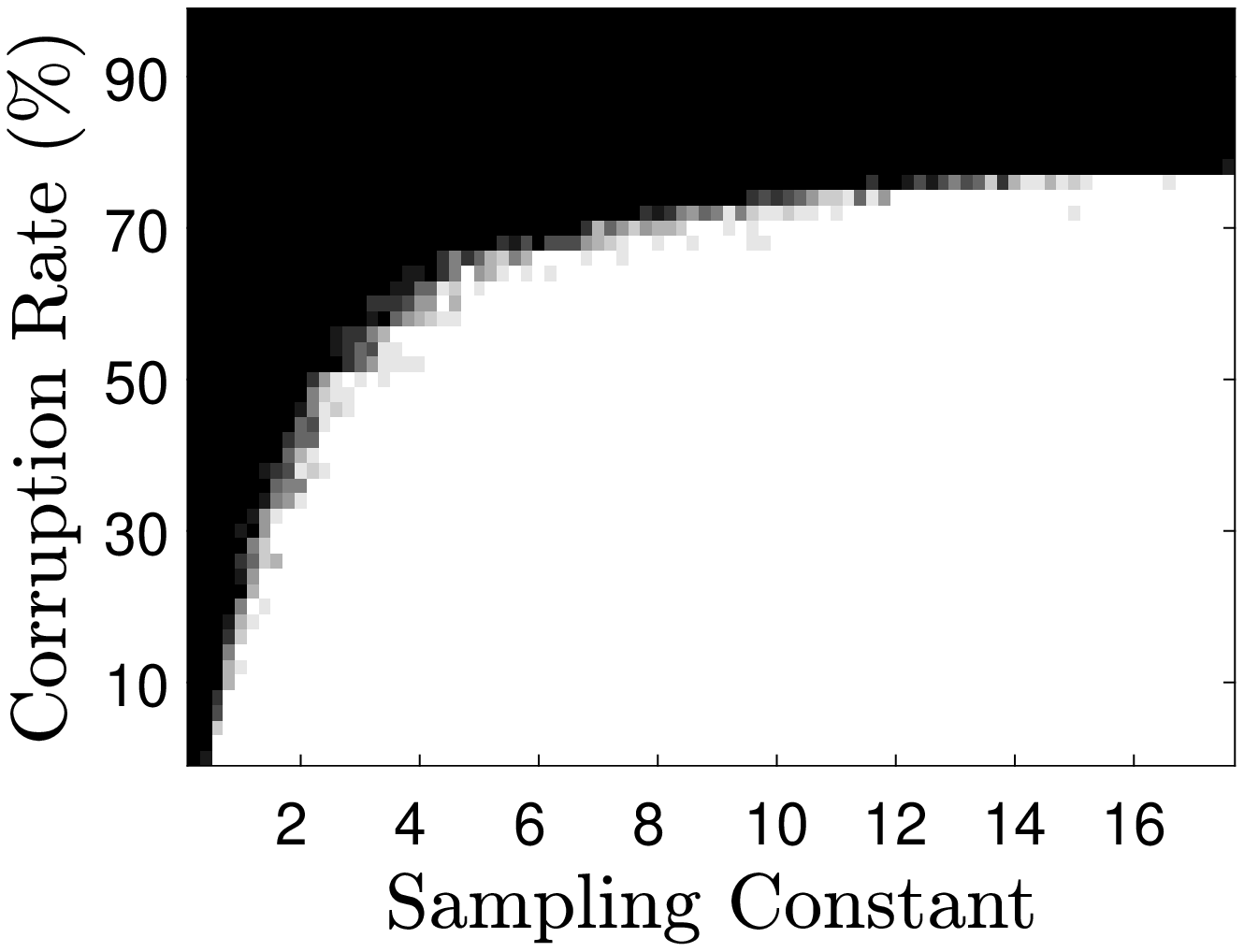}
		\hfill
		\includegraphics[width=0.24\linewidth]{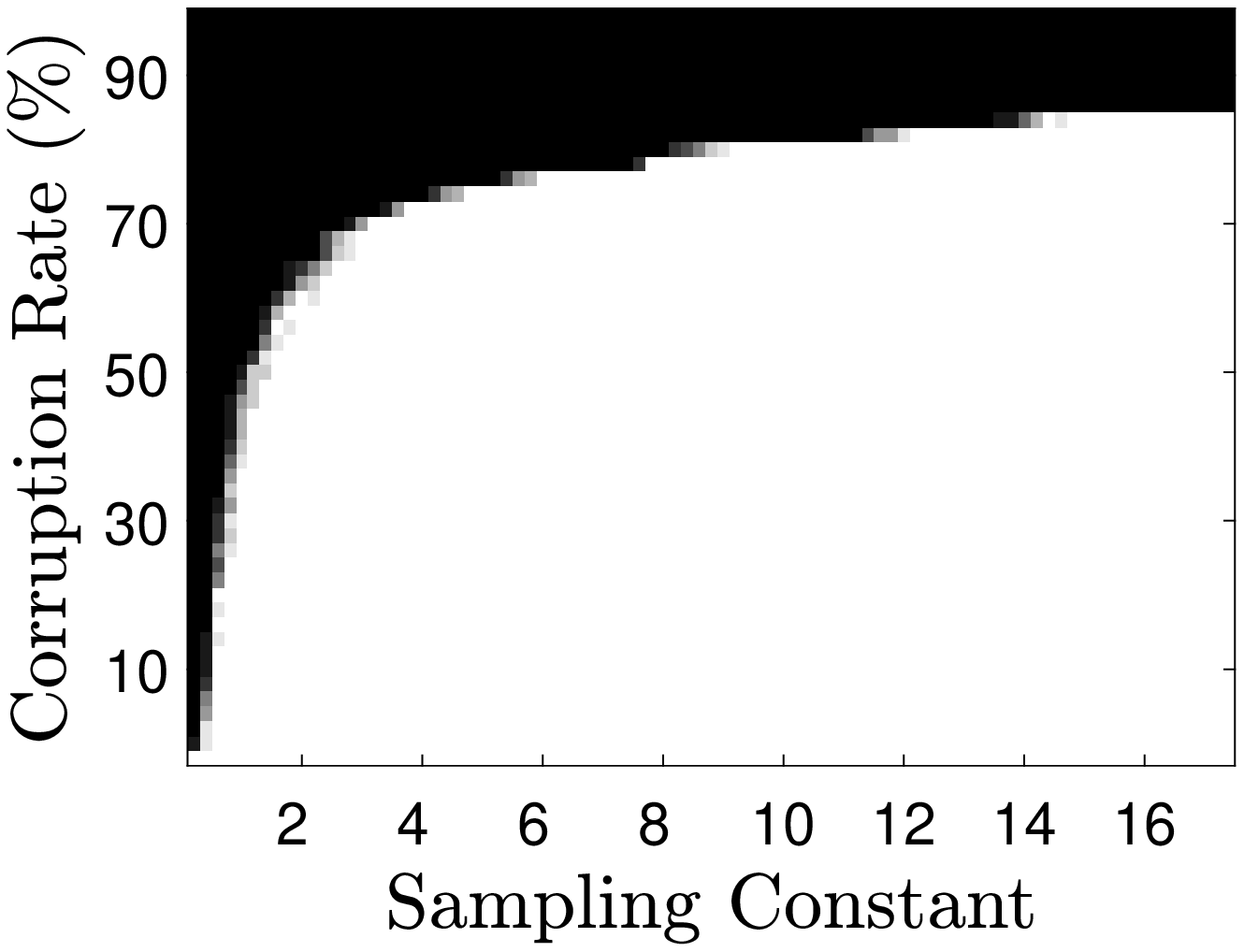}
		\hfill
		\includegraphics[width=0.24\linewidth]{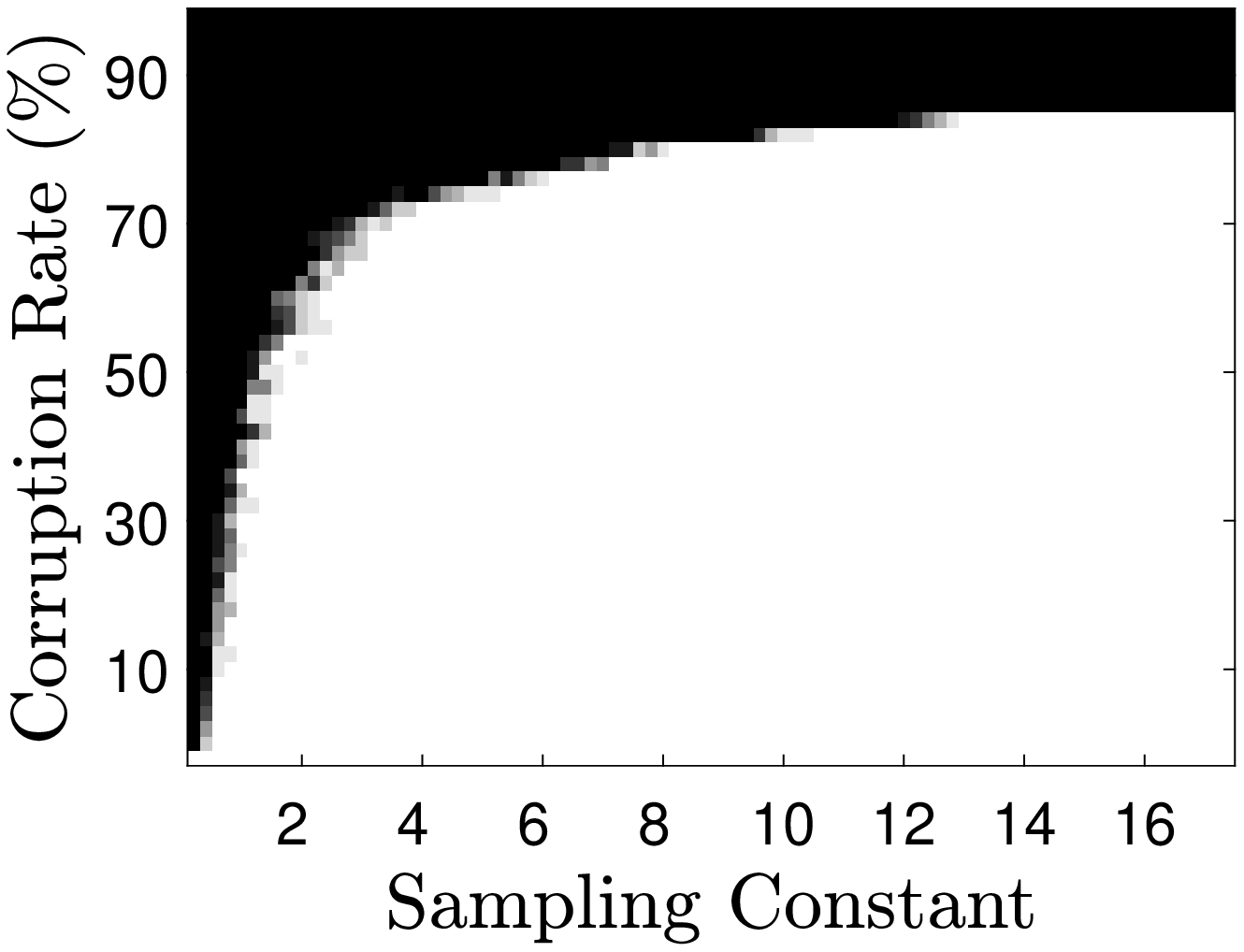}\\\vspace{0.1in}
		
		\includegraphics[width=0.24\linewidth]{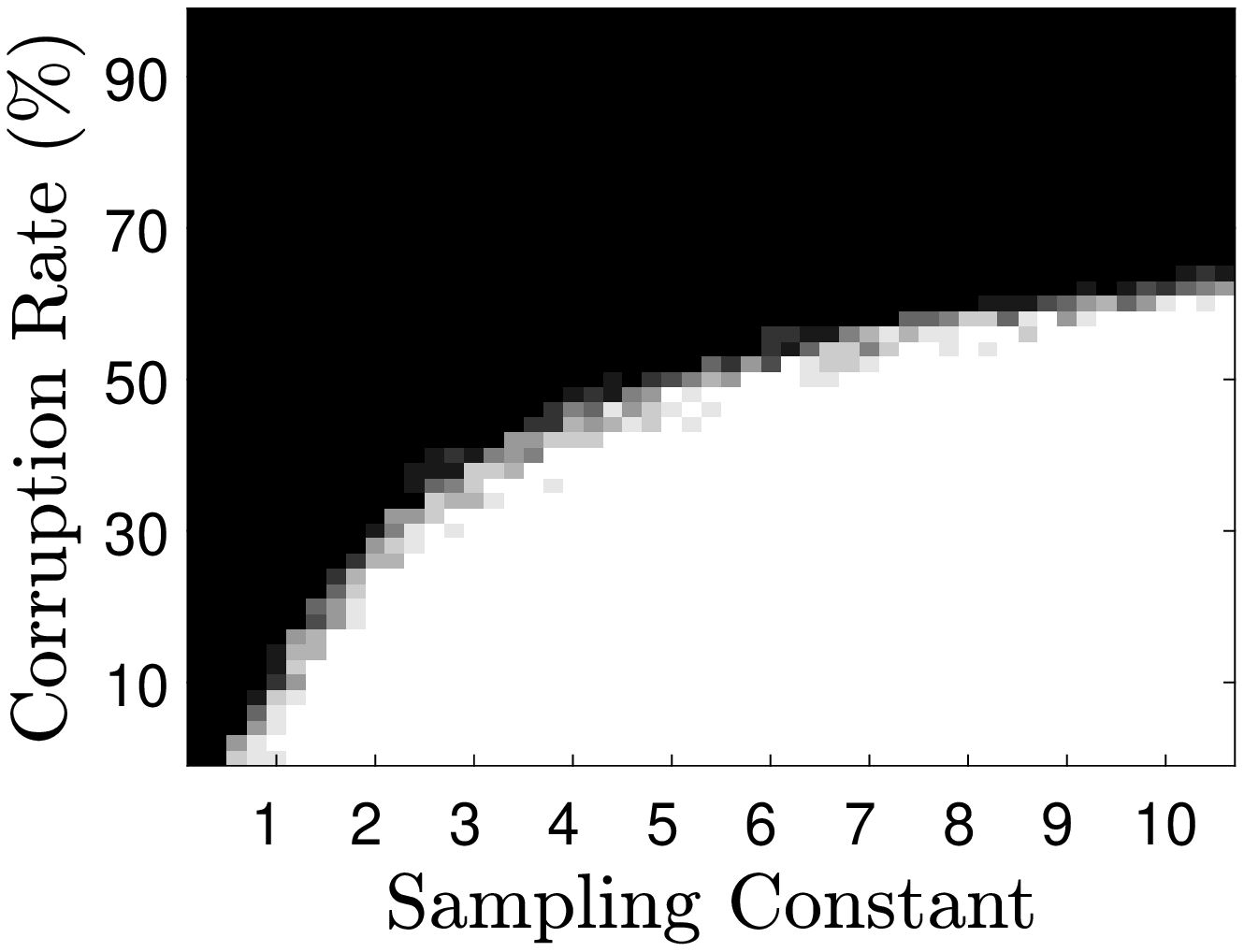}
		\hfill
		\includegraphics[width=0.24\linewidth]{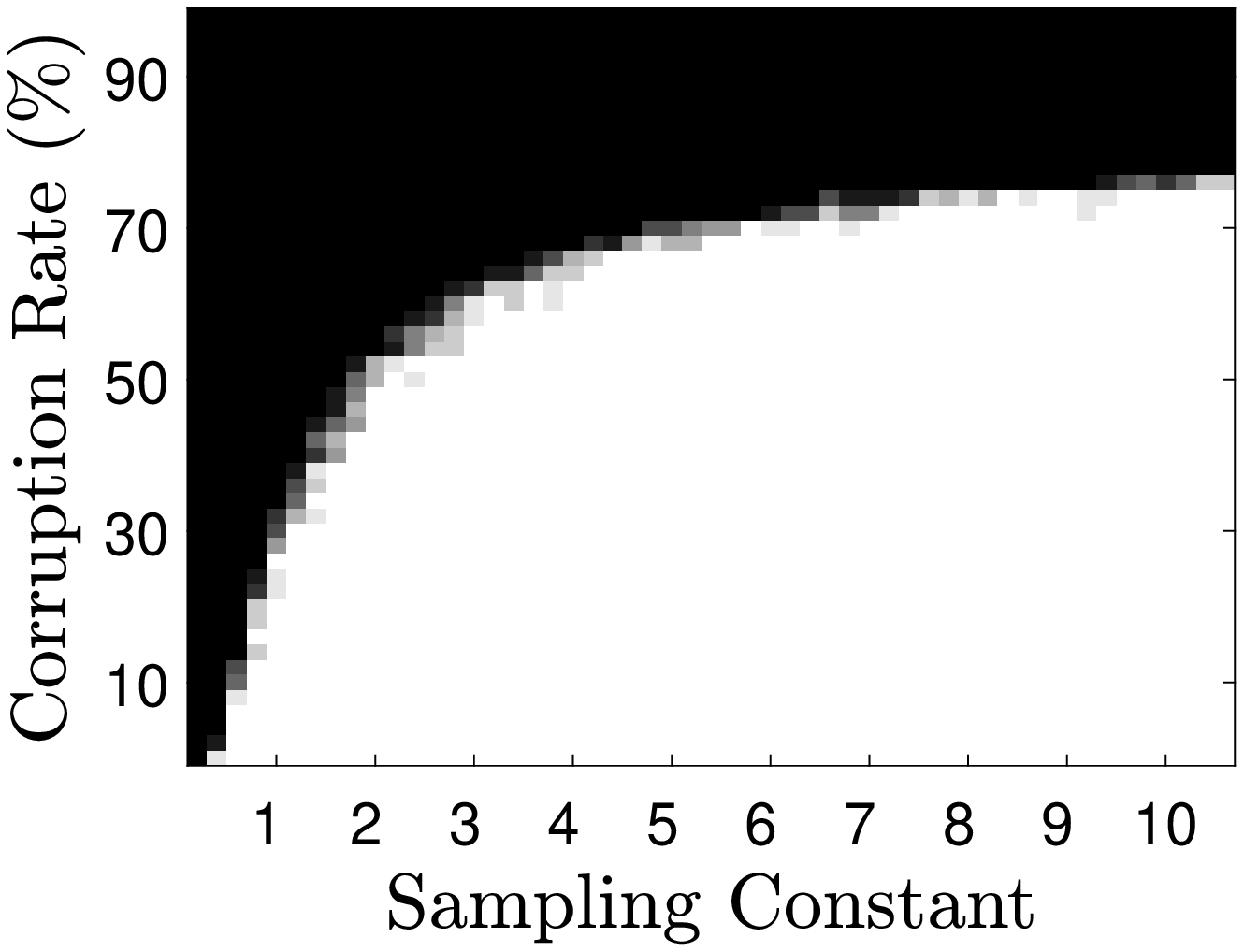}
		\hfill
		\includegraphics[width=0.24\linewidth]{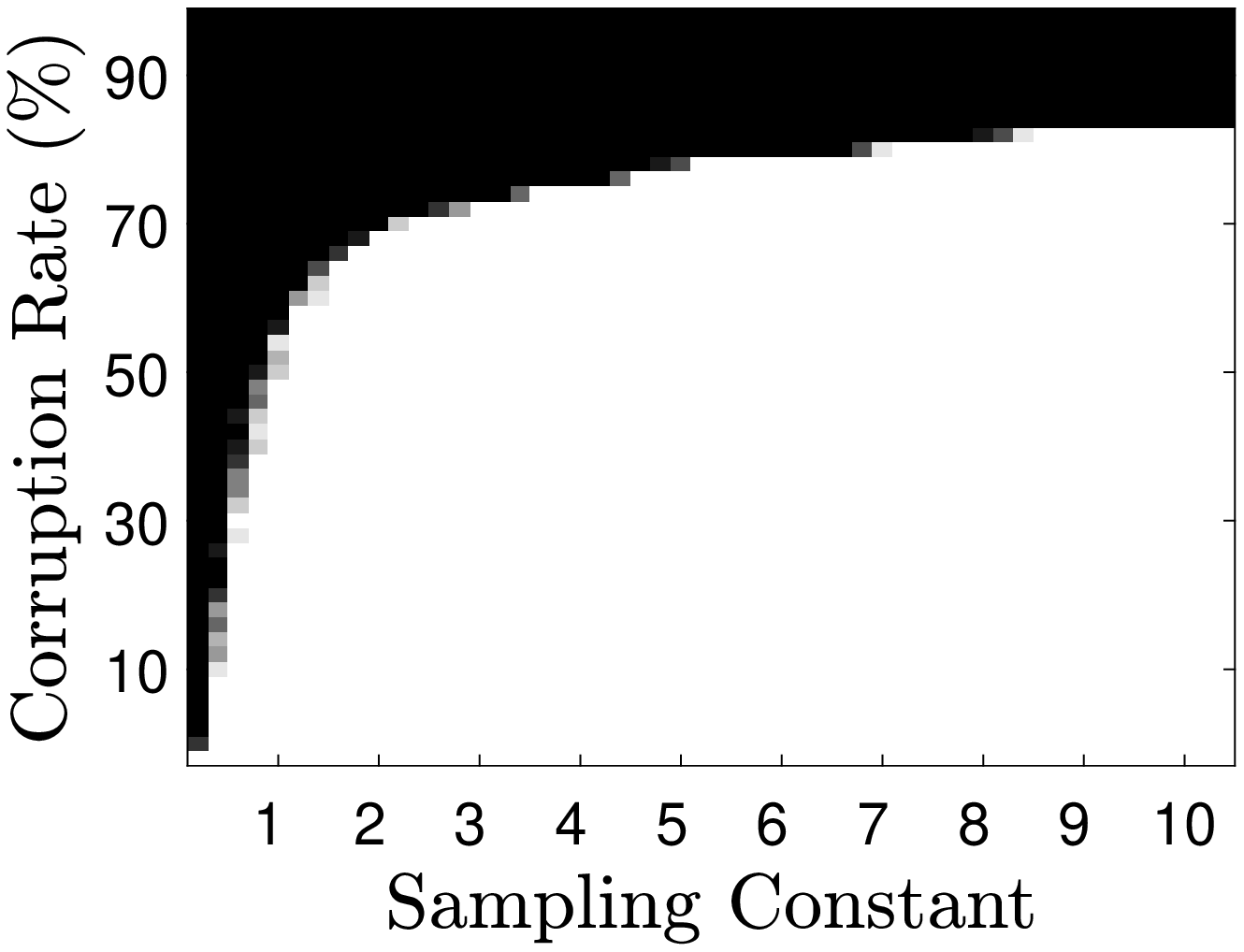}
		\hfill
		\includegraphics[width=0.24\linewidth]{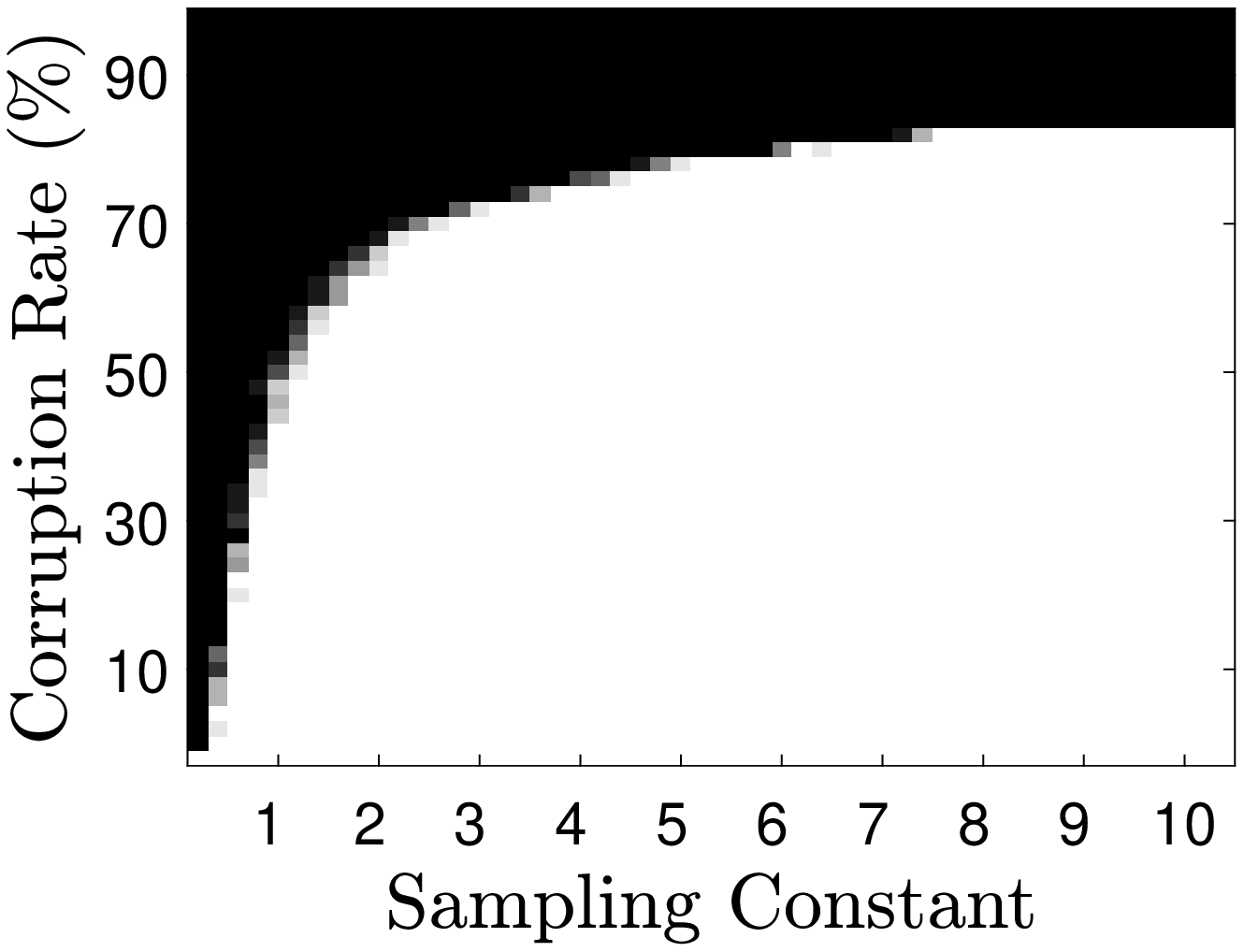}\\\vspace{0.1in}
		
		\includegraphics[width=0.24\linewidth]{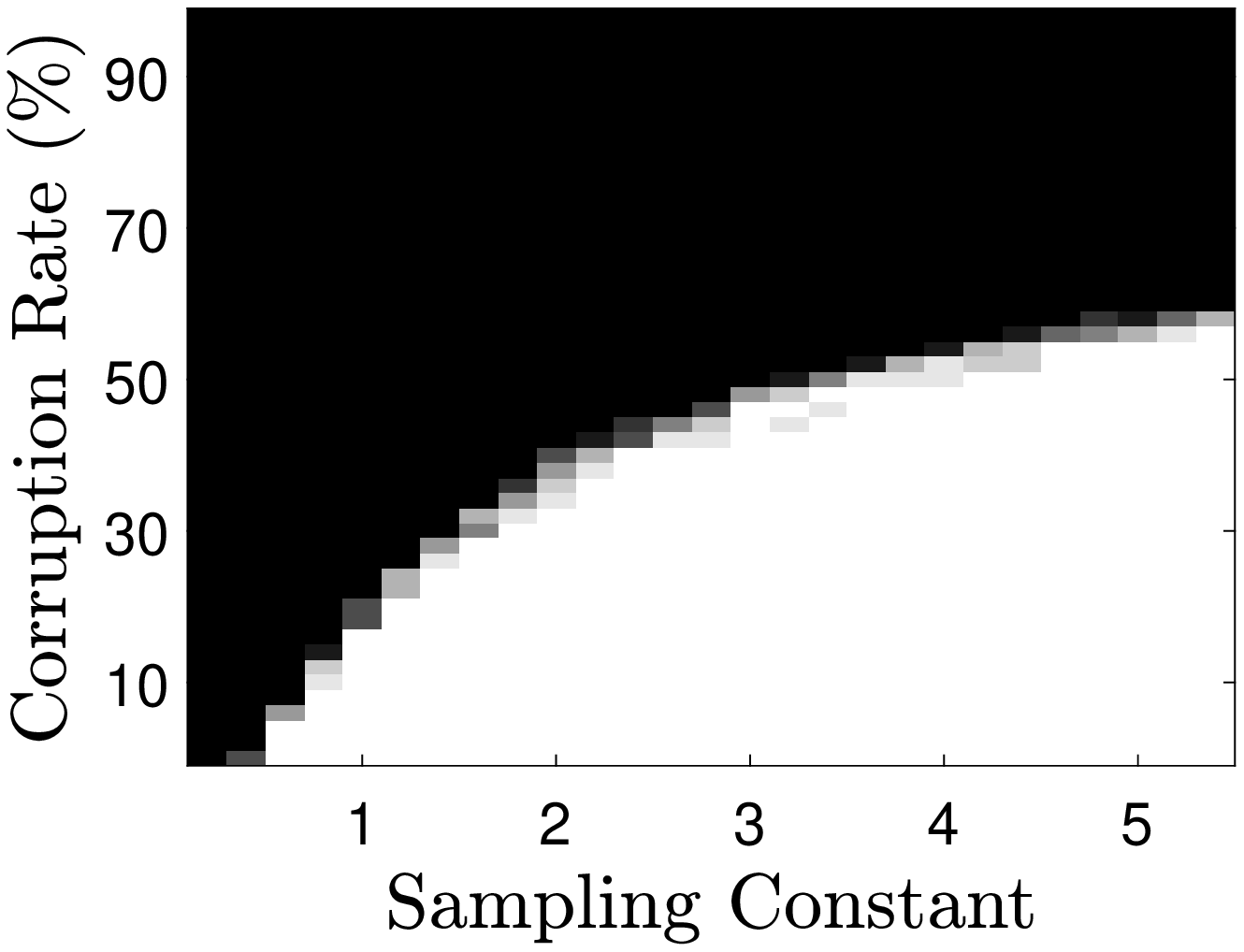}
		\hfill
		\includegraphics[width=0.24\linewidth]{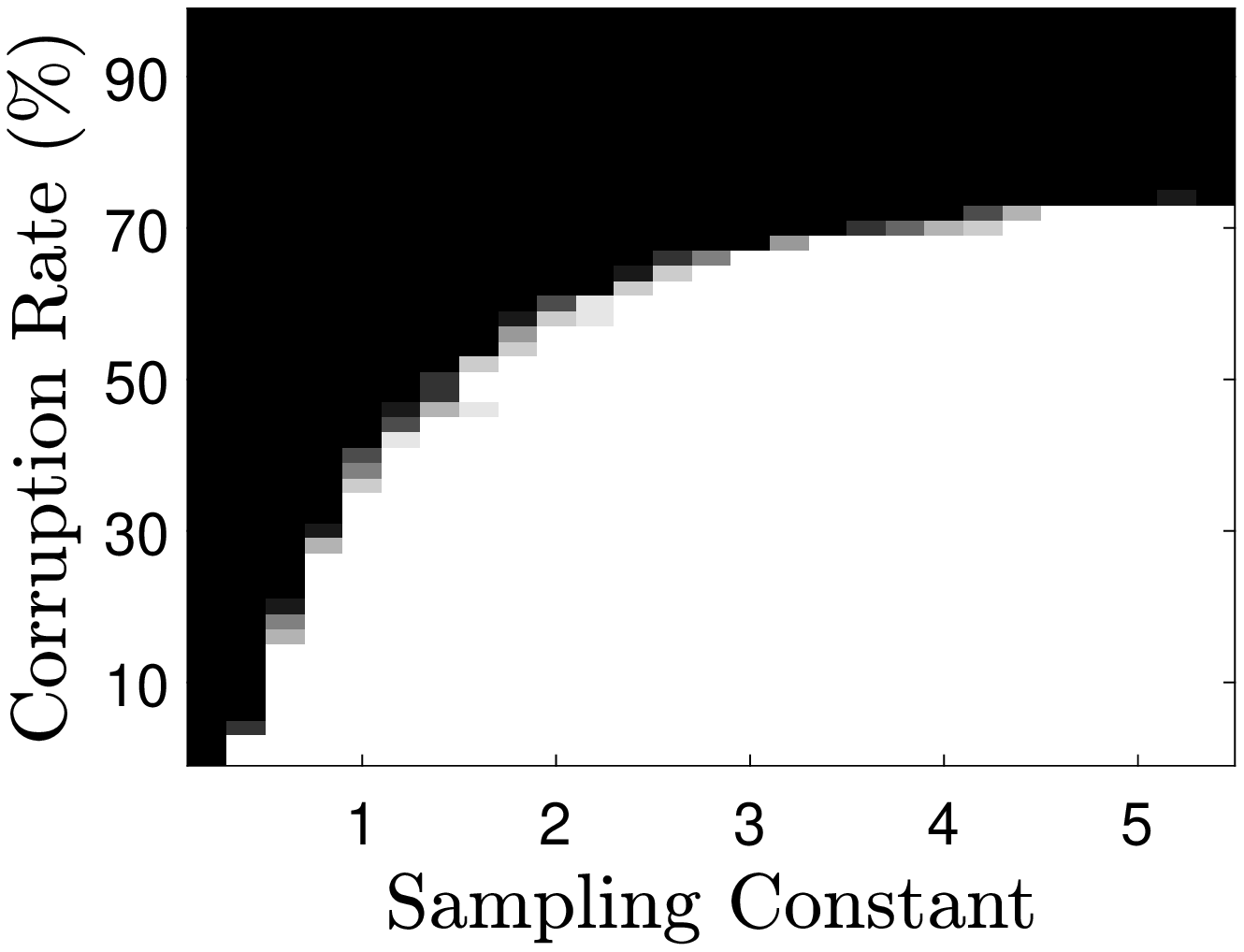}
		\hfill
		\includegraphics[width=0.24\linewidth]{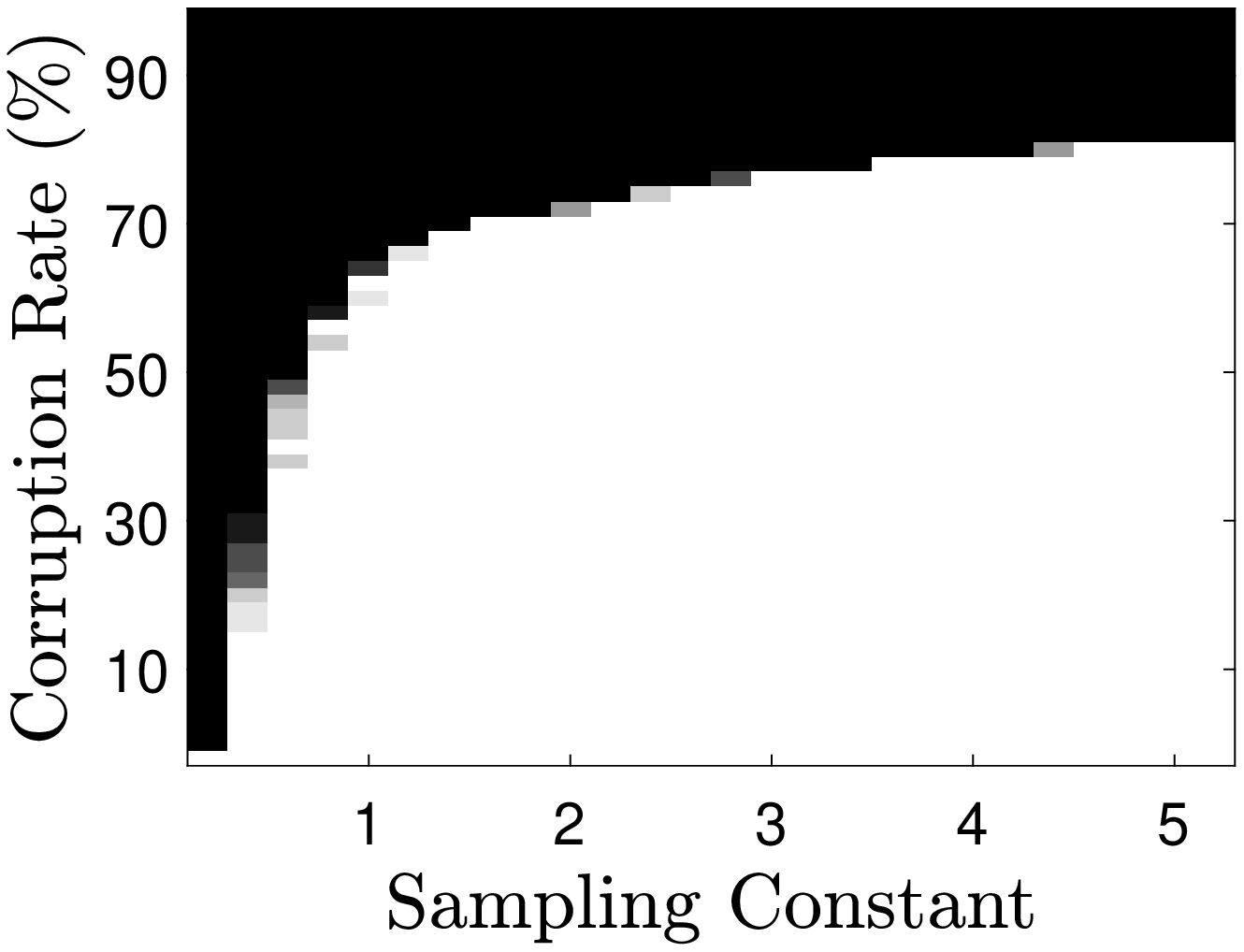}
		\hfill
		\includegraphics[width=0.24\linewidth]{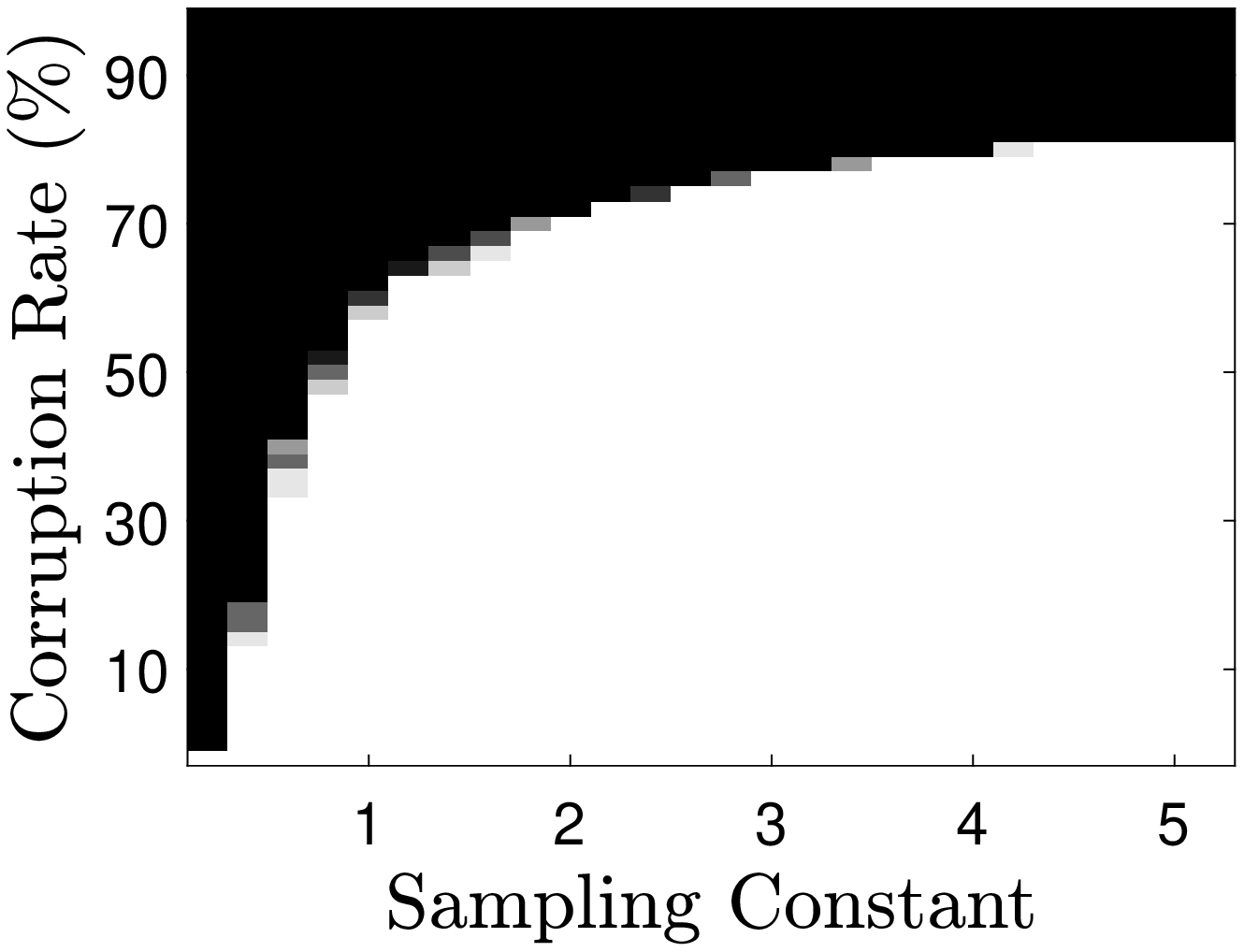}\\
		\caption{Empirical phase transition in corruption rate $\alpha$ and sampling constant $\upsilon$. Left to Right: RTCUR-FF, RTCUR-RF, RTCUR-FC, RTCUR-RC, \textbf{Top}: $r = 3$. \textbf{Middle}: $r = 5$. \textbf{Bottom}: $r = 10$.} \label{FIG:phase}
		\vspace{-0.1in}
	\end{figure}

	This section presents a set of numerical experiments that compare the empirical performance of \alg\ with several state-of-the-art robust matrix/tensor PCA algorithms, including  Riemannian gradient descent (RGD) \cite{cai2022generalized}, alternating direction method of multipliers (ADMM) \cite{lu2019tensor}, accelerated alternating projections (AAP) \cite{cai2019accelerated}, and iterative robust CUR (IRCUR) \cite{cai2020rapid}. RGD and ADMM are designed for the TRPCA task, while AAP and IRCUR are designed for the traditional matrix RPCA task. Note that RGD is Tucker-rank based and ADMM is tubal-rank based.
	
	In each subsection, we evaluate the performance of all four proposed variants, \alg-FF,  \alg-RF, \alg-FC, and \alg-RC. However, for the network clustering experiment, we only use  fixed sampling (\alg-FF and \alg-FC). As the coauthorship network is highly sparse, the resampling variants may diminish the core tensor with a  high  probability.
	
	The rest of this section is structured as follows. In \Cref{sec:synthetic}, we present two synthetic experiments. Specifically, \Cref{sec:phase} examines the empirical relationship between outlier tolerance and sample size for \alg, while \Cref{sec:speed_comparison} demonstrates the speed advantage of \alg\ over state-of-the-art. In \Cref{sec:video,sec:face}, we apply \alg\ to two real-world problems, namely  face modeling and color video background subtraction. In \Cref{sec:network}, we apply \alg\ on network clustering applications and analyze the obtained results.
	
	We obtain the codes for all compared algorithms from the authors' websites and hand-tune the parameters for their best performance. For \alg, we sample $|I_i|=\upsilon r_i\log(d_i)$ (and $|J_i|=\upsilon r_i \log(\prod_{j\neq i}d_j)$ for Fiber variants) for all $i$, and $\upsilon$ is called the \textit{sampling constant} through this section. All the tests are executed from Matlab R2020a on an Ubuntu workstation with an Intel i9-9940X CPU and 128GB RAM. The relevant codes are available at  
	\href{https://github.com/huangl3/RTCUR}{https://github.com/huangl3/RTCUR}.

	\subsection{Synthetic Examples} \label{sec:synthetic}
	
	\begin{figure}[th]
		\centering
		\subfloat{\includegraphics[width = 0.45\linewidth]{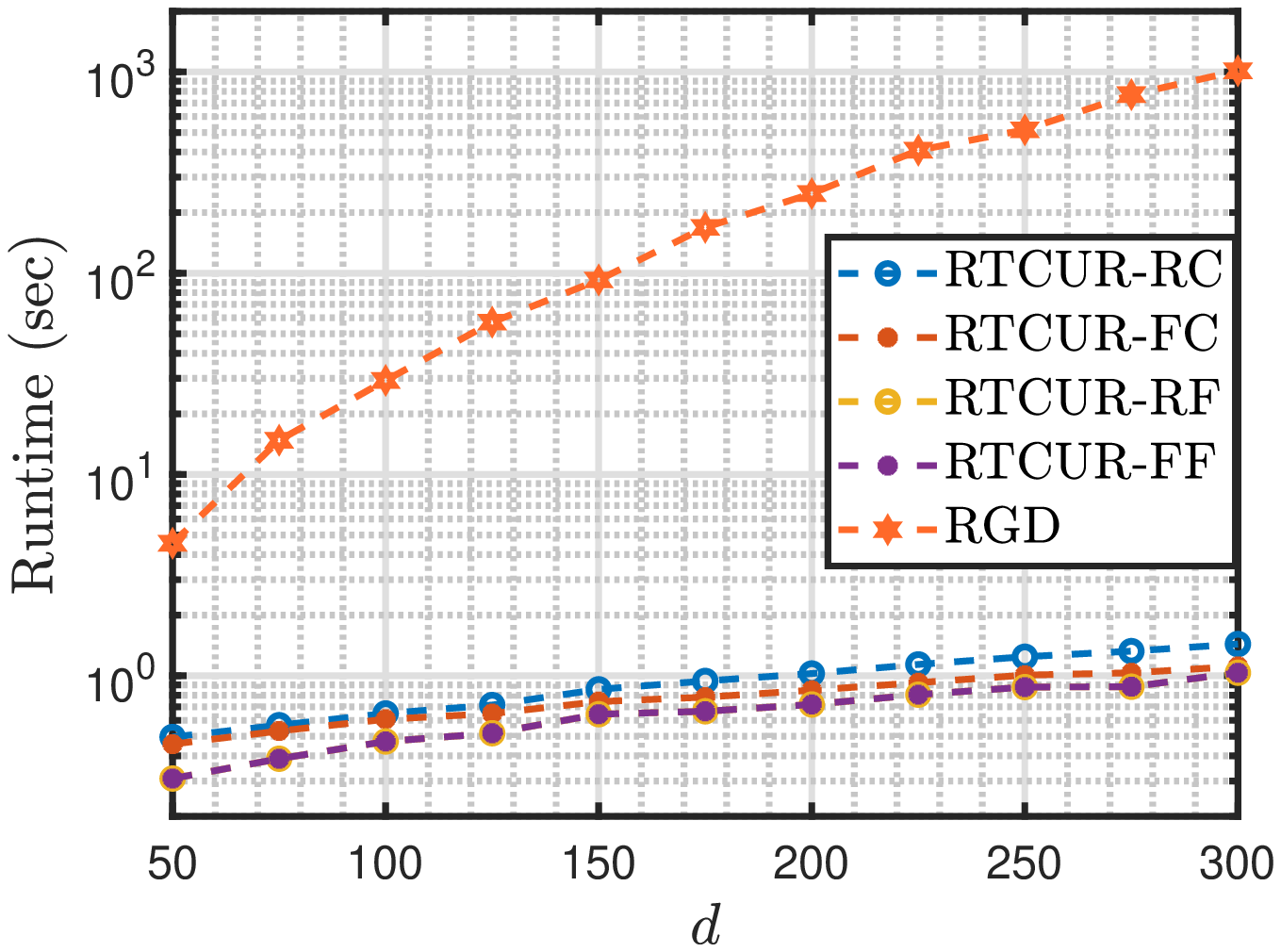}} 
		\subfloat{\includegraphics[width = 0.45\linewidth]{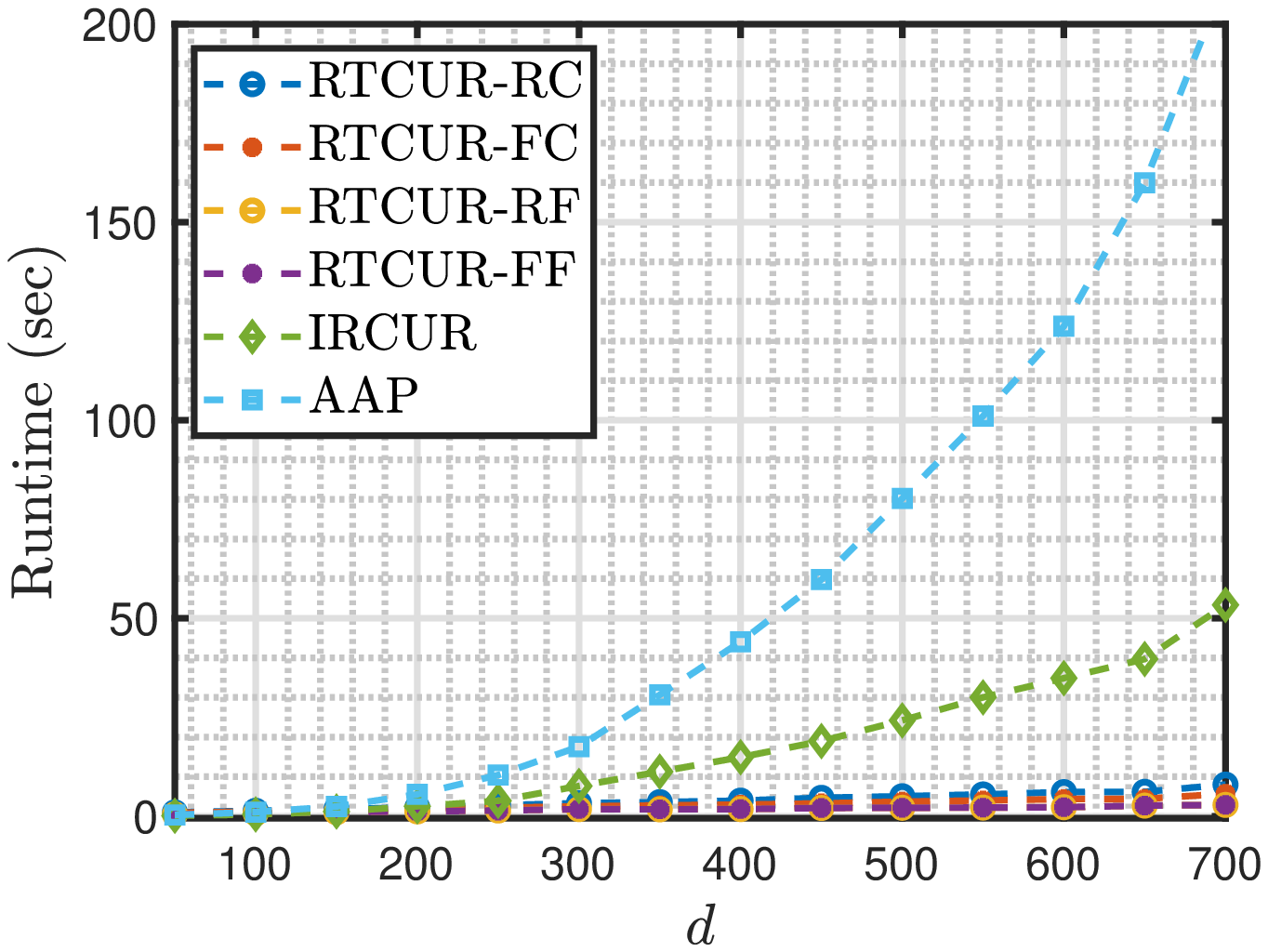}}
		\caption{Runtime vs.~dimension comparison among variants of \alg, RGD, IRCUR and AAP on tensors with size $d\times d\times d$~ and Tucker rank $(3,3,3)$. The RGD method proceeds relatively slowly for larger tensors, so we only test the RGD runtime for tensors with a size smaller than 300 for each mode.}
		\label{fig:CUR_time}
		\vspace{-0.1in}
	\end{figure}
	
	For the synthetic experiments, we use $d:=d_1=\cdots=d_n$ and $r:=r_1=\cdots=r_n$. The observed tensor $\cX$ is composed as $\cX=\cL^\star+\cS^\star$. 
	To generate $n$-mode $\cL^\star\in\mathbb{R}^{d\times\cdots\times d}$ with Tucker rank $(r,\cdots,r)$, we take $\cL^\star=\mathcal{Y}\times_1 Y_1 \times_2 \cdots \times_n Y_n$ where $\mathcal{Y}\in\mathbb{R}^{r\times\cdots\times r}$ and $\{Y_i\in\mathbb{R}^{d\times r}\}_{i=1}^n$ are Gaussian random tensor and matrices with standard normal entries. To generate the sparse outlier tensor $\cS^\star$, we uniformly sample $\lfloor\alpha d^n\rfloor$ entries to be the support of $\cS^\star$ and the values of the non-zero entries are uniformly sampled from the interval $[-\mathbb{E}(|\cL^\star_{i_1,\cdots,i_n}|),\mathbb{E}(|\cL^\star_{i_1,\cdots,i_n}|)]$.
	
	\subsubsection{Phase Transition} \label{sec:phase}
	We study the empirical relation between the outlier corruption rate $\alpha$ and sampling constant $\upsilon$ for all four variants of \alg\ using $300\times300\times300$ (i.e., $n=3$ and $d=300$) problems with Tucker rank $(r,r,r)$, where $r=3,5,$ or $10$. We set the thresholding parameters  to  $\zeta^{(0)}=\|\cL\|_\infty$ and $\gamma=0.7$, and use the stopping condition   $e^{(k)}<10^{-5}$. A test example is considered successfully solved if $\|\cL^\star-\cL^{(k)}\|_\fro/\|\cL^\star\|_\fro\leq 10^{-3}$.  For each pair of $\alpha$ and $\upsilon$, we generate $10$ test examples. 
	
	We summarize the experimental results in \Cref{FIG:phase}, where a white pixel means all $10$ test examples are successfully solved under the corresponding problem parameter setting, and a black pixel means all $10$ test cases fail. 
	Upon observation,  one has that the Chidori variants, \alg-FC, and \alg-RC, can recover the low-rank tensor with higher outlier rate than the Fiber variants with the same sampling constant $\upsilon$. This expected behavior can be attributed to the fact that Chidori sampling accesses more data from the tensor with the same $\upsilon$, leading to more stable performance compared to Fiber sampling. Furthermore, \alg-RF also outperforms \alg-FF, as resampling the fibers allows access to more tensor entries throughout the process. This phenomenon holds true for the Chidori variants as well, where \alg-FC and \alg-RC exhibit very similar phase transitions. Additionally, we observe that smaller values of $r$ tolerate more outliers, as larger values of $r$ make the TRPCA task more complex. Increasing $\upsilon$ improves outlier tolerance, but at the cost of larger subtensors and more fibers to be sampled, resulting in longer computational time

	\begin{figure}[t]
		\centering  \subfloat{\includegraphics[width = .45\linewidth]{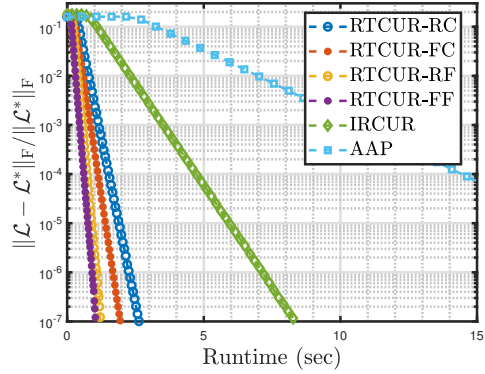}}
		\caption{
			Runtime vs.~relative error comparison among \algf, \algr, AAP, and IRCUR on tensor with size $500\times 500\times 500$~ and Tucker rank ($3,3,3$).} \label{fig:CUR_time2}
		\vspace{-0.1in}
	\end{figure}

	\subsubsection{Computational Efficiency} \label{sec:speed_comparison}
	In this section, we present a comparative analysis of the computational efficiency of \alg\ and the state-of-the-art tensor/matrix RPCA algorithms mentioned in \Cref{sec:experiment}. To apply the matrix RPCA algorithms, we first unfold an $n$-mode $d\times \cdots \times d$ tensor to a $d\times d^{n-1}$ matrix and then solve the matrix RPCA problem by setting the matrix rank as $r$. It is worth noting that we do not apply ADMM to this fixed-Tucker-rank experiment, as the target tensor rank for ADMM is tubal-rank instead of Tucker-rank \cite{lu2019tensor}.
	For all tests, we set the underlying low-Tucker-rank component as rank $[3,3,3]$ and add $20\%$ corruption. We set parameters $\upsilon=3$, $\zeta^{(0)}=\|\cL\|_\infty$, $\gamma = 0.7$ for all four variants of \alg. The reported runtime is averaged over $20$ trials.
	
	In \Cref{fig:CUR_time}, we investigate the $3$-mode TRPCA problem with varying dimension $d$ and compare the total runtime (all methods halt when relative error $e^{(k)}<10^{-5}$). It is evident that all variants of \alg\ are significantly faster than the compared algorithms when $d$ is large. Next, we evaluate the convergence behavior of the tested algorithms in \Cref{fig:CUR_time2}. We exclude RGD as well because running RGD on this experiment is too expensive. We find all the tested algorithms converge linearly, and variants of \alg\ run the fastest. Moreover, as discussed in \Cref{sec:4vars}, the Fiber sampling has a running time advantage over Chidori sampling with the same sampling constant, and fixed sampling runs slightly faster than resampling in all tests.
	
	\subsection{Robust Face Modeling} \label{sec:face}
	
	\begin{figure}[h]
		\centering
		\includegraphics[width=0.78 \linewidth]{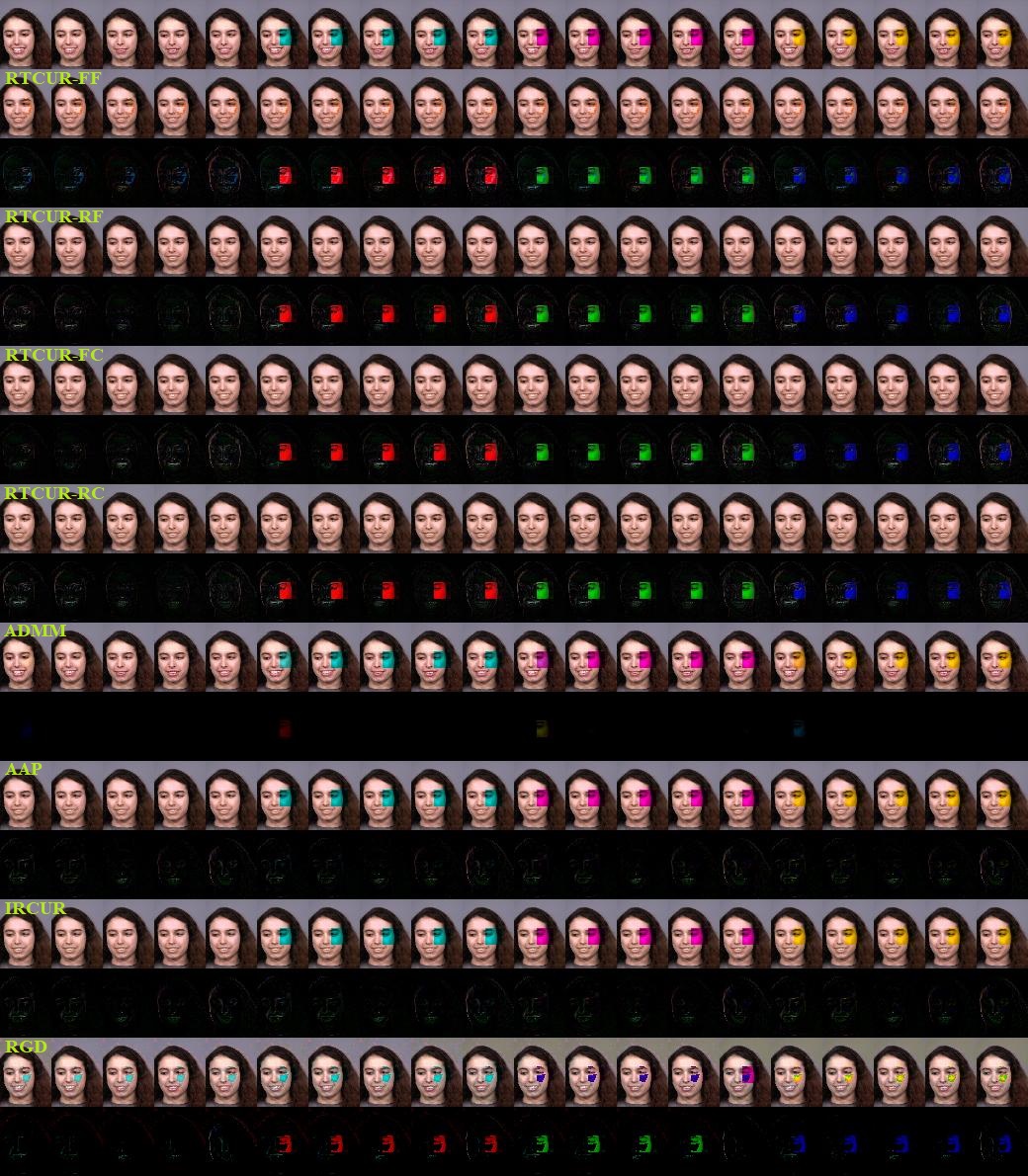}
		\caption{
			Visual results for robust face modeling. The \textbf{top row} contains the corrupted faces, the \textbf{second and third rows} are the recovered faces and detected outliers outputted by \alg-FF; the \textbf{fourth and fifth rows} are results from \alg-RF; the \textbf{sixth and seventh rows} are results from \alg-FC; the \textbf{eighth and ninth} rows are results from \alg-RC; the \textbf{tenth and eleventh rows} are results from ADMM; the \textbf{twelfth and thirteenth rows} are results from AAP; the \textbf{fourteenth and fifteenth rows} are results from IRCUR; the \textbf{sixteenth and seventeenth rows} are results from RGD.} \label{fig: face2}
		\vspace{-0.1in}
	\end{figure}
	
	In this section, we apply the four variants of \alg\ and compare them with the aforementioned tensor/matrix RPCA algorithms on the robust face modeling task using data from the UT Dallas database \cite{o2005video}. The dataset consists of a face speech video of approximately 5 seconds with a resolution of 360 $\times$ 540.  We extract 10 non-successive frames and mark a monochromatic block on different color channels for 10 distinct frames per color. This process results in a total of 40 color frames, including the original $10$ unmarked frames. As a monochromatic frame usually does not have a low-rank structure, we vectorize each color channel of each frame into a vector and construct a (height $\cdot$ width) $\times$ 3 $\times$ frames tensor. The targeted Tucker rank is set as $\bm{r} = (3,3,3)$ for all videos.  For those matrix algorithms, including AAP and IRCUR, we unfold the tensor to a (height $\cdot$ width) $\times$ ($3~\cdot$ frames) matrix and rank $3$ is used. We set \alg\ parameters $\upsilon=2$, $\zeta^{(0)}=255$, $\gamma=0.7$ in this experiment. 
	
	\Cref{fig: face2} presents the test examples and visual results; \Cref{tab:face time} summarizes the runtime for each method applied on this task. One can see that the matrix-based methods fail to detect the monochromatic outlier blocks since they lose the structural connection between color channels after matricization, albeit they spend less time on this task. In contrast, all variants of \alg\ successfully detect the outlier blocks. The other two TPRCA methods, ADMM and RGD, partially detect the outlier blocks. Since ADMM is based on tubal decomposition, it is not a surprise to see different performances in this experiment.
	The empirical results in this section verify our claim in \Cref{rmk:uncertainty principle}. 
	
	\begin{table}[h]
		\caption{Runtime comparison (in seconds) for face modeling task. The matrix RPCA approaches (AAP and IRCUR) meet the termination condition earlier with the unfolded tensor, but they failed to detect the artificial noise in this task (see \Cref{fig: face2}).}
		\label{tab:face time}
		\centering
		\vspace{-0.1in}
		\begin{small}
			\begin{tabular}{c|c||c|c}\toprule
				\textsc{Method} & \textsc{Runtime}& \textsc{Method} & \textsc{Runtime} \\ \midrule
				\alg-FF & 2.247 & ADMM &  30.61 \\
				\alg-RF & 2.289 & AAP &  1.754 \\
				\alg-FC & 2.319 & IRCUR &  1.307 \\
				\alg-RC & 2.701 &  RGD &  1430.8 \\\bottomrule
			\end{tabular}
		\end{small}
	\end{table}

	\subsection{Color Video Background Subtraction} \label{sec:video}
	\begin{table}[t]
		\caption{Video information and runtime comparison (in seconds) for color video background subtraction task.  
		}\label{table:video}
		\centering
		\vspace{-0.1in}
		\resizebox{\textwidth}{!}{
			\begin{tabular}{ c|c|ccccccc}
				\toprule
				~     &\textsc{Video size} & {RTCUR-FF} &{RTCUR-RF} &{RTCUR-FC} &{RTCUR-RC}  &{ADMM} &{AAP} & {IRCUR}               \cr
				\midrule
				\textit{Shoppingmall}
				&
				$256\times 320\times 3\times  1250$ &\textbf{3.53} & 5.83 & 10.68 & 10.75 & 783.67 & 50.38 & 15.71   \cr
				\textit{Highway}
				&
				$240\times 320\times 3\times440$  & \textbf{3.15} & 5.47 & 6.80 & 7.55 & 168.55  & 18.10 & 3.87 \cr
				\textit{Crossroad}
				&
				$350\times640\times 3\times600$  & \textbf{6.15} & 13.33 & 8.46 & 12.01 & 1099.3 & 97.85 & 35.47 \cr
				\textit{Port}
				&
				$480\times640\times 3\times1000$  & \textbf{11.04} & 18.34 & 26.63 & 27.93  & 2934.3 & 154.30 & 71.64 \cr
				\textit{Parking-lot}
				&
				$360\times640\times 3\times400$  & \textbf{3.79} & 4.52 & 6.62 & 8.14 & 854.50 & 34.70 & 17.38 \cr
				\bottomrule
			\end{tabular}
		}
	\end{table}

	\begin{figure}[ht]
		\centering
		\vspace{-0.17in}
		\subfloat[Original]{\includegraphics[width=0.12\linewidth]{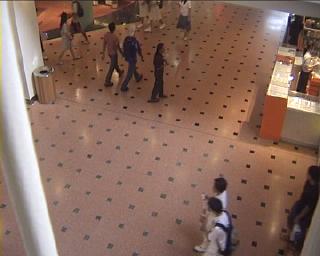}}\hfill
		\subfloat[RTCUR-FF]{\includegraphics[width=0.12\linewidth]{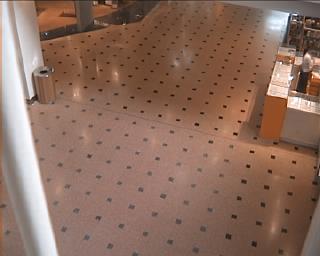}}\hfill
		\subfloat[RTCUR-RF]{\includegraphics[width=0.12\linewidth]{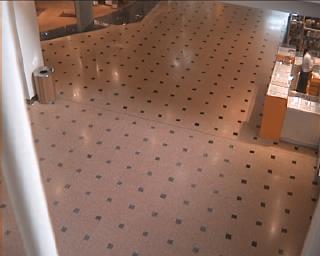}}\hfill
		\subfloat[RTCUR-FC]{\includegraphics[width=0.12\linewidth]{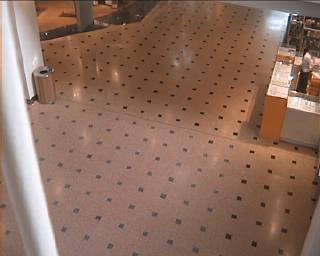}}\hfill
		\subfloat[RTCUR-RC]{\includegraphics[width=0.12\linewidth]{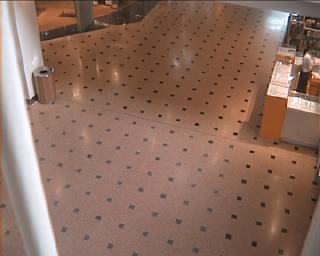}}\hfill
		\subfloat[ADMM]{\includegraphics[width=0.12\linewidth]{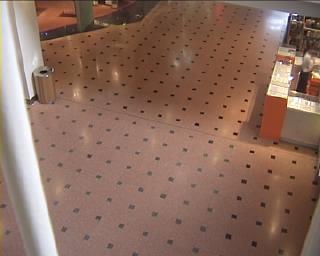}}\hfill
		\subfloat[AAP]{\includegraphics[width=0.12\linewidth]{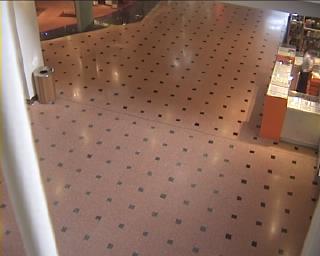}}\hfill
		\subfloat[IRCUR]{\includegraphics[width=0.12\linewidth]{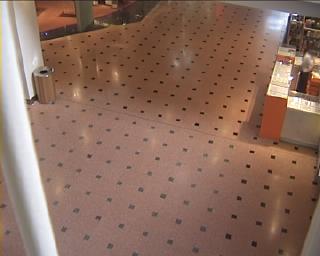}}
		\vspace{-0.12in}
		\\
		\subfloat{\includegraphics[width=0.12\linewidth]{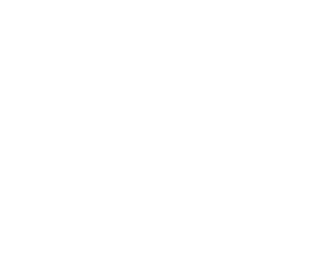}}\hfill
		\subfloat{\includegraphics[width=0.12\linewidth]{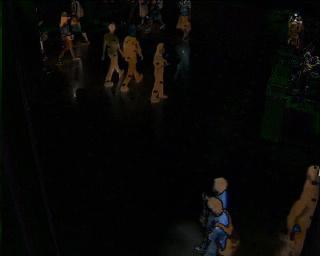}}\hfill
		\subfloat{\includegraphics[width=0.12\linewidth]{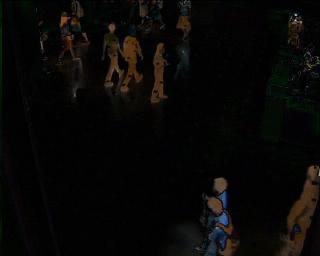}}\hfill
		\subfloat{\includegraphics[width=0.12\linewidth]{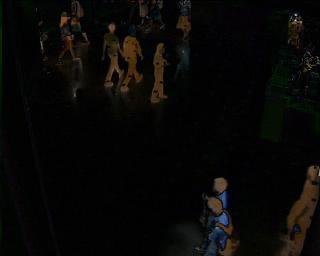}}\hfill
		\subfloat{\includegraphics[width=0.12\linewidth]{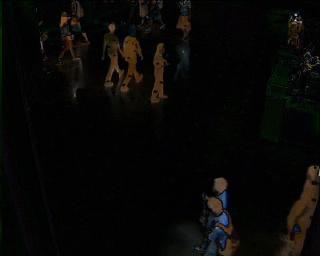}}\hfill
		\subfloat{\includegraphics[width=0.12\linewidth]{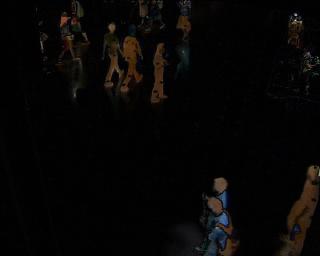}}\hfill
		\subfloat{\includegraphics[width=0.12\linewidth]{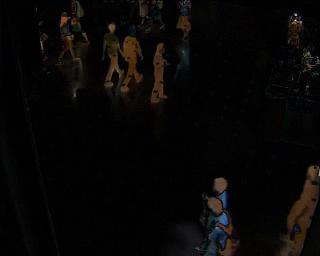}}\hfill
		\subfloat{\includegraphics[width=0.12\linewidth]{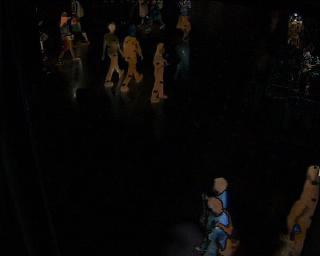}}
		\\
		\vspace{-0.12in}
		\subfloat{\includegraphics[width=0.12\linewidth]{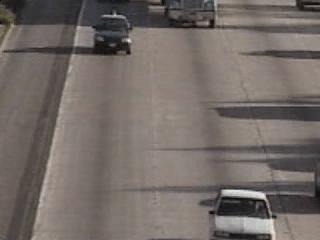}}\hfill
		\subfloat{\includegraphics[width=0.12\linewidth]{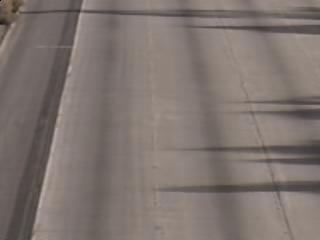}}\hfill
		\subfloat{\includegraphics[width=0.12\linewidth]{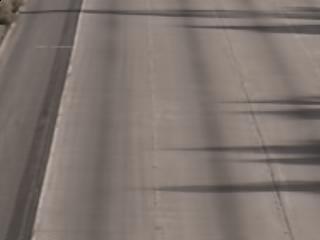}}\hfill
		\subfloat{\includegraphics[width=0.12\linewidth]{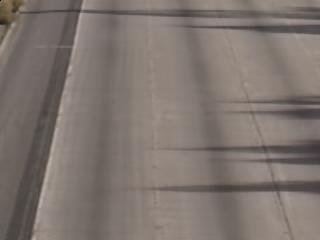}}\hfill
		\subfloat{\includegraphics[width=0.12\linewidth]{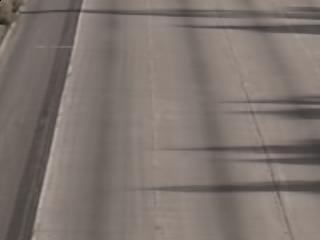}}\hfill
		\subfloat{\includegraphics[width=0.12\linewidth]{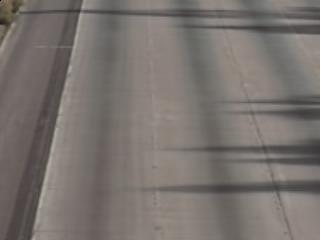}}\hfill
		\subfloat{\includegraphics[width=0.12\linewidth]{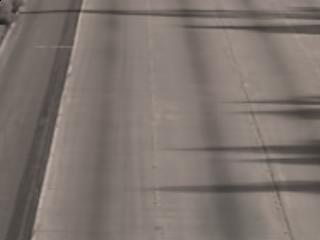}}\hfill
		\subfloat{\includegraphics[width=0.12\linewidth]{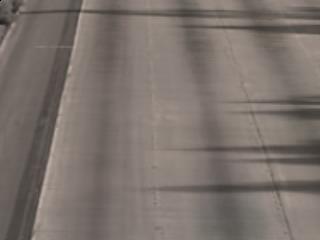}}
		\\
		\vspace{-0.12in}
		\subfloat{\includegraphics[width=0.12\linewidth]{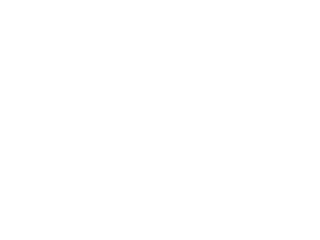}}\hfill
		\subfloat{\includegraphics[width=0.12\linewidth]{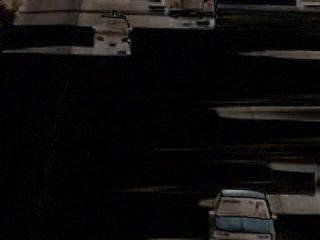}}\hfill
		\subfloat{\includegraphics[width=0.12\linewidth]{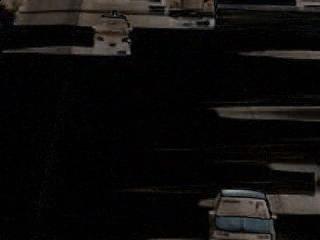}}\hfill
		\subfloat{\includegraphics[width=0.12\linewidth]{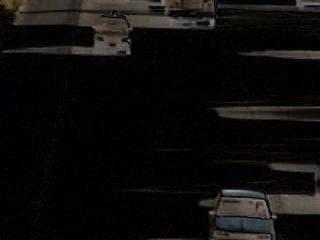}}\hfill
		\subfloat{\includegraphics[width=0.12\linewidth]{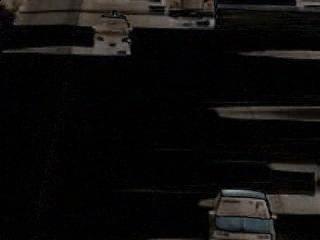}}\hfill
		\subfloat{\includegraphics[width=0.12\linewidth]{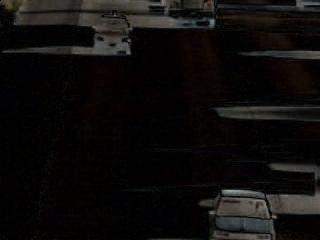}}\hfill
		\subfloat{\includegraphics[width=0.12\linewidth]{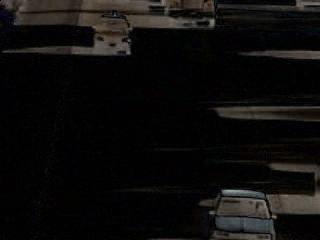}}\hfill
		\subfloat{\includegraphics[width=0.12\linewidth]{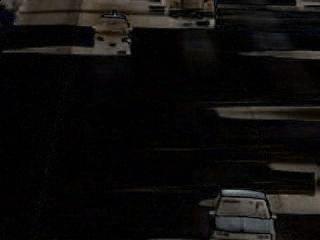}} 
		\\
		\vspace{-0.12in}
		\subfloat{\includegraphics[width=0.12\linewidth]{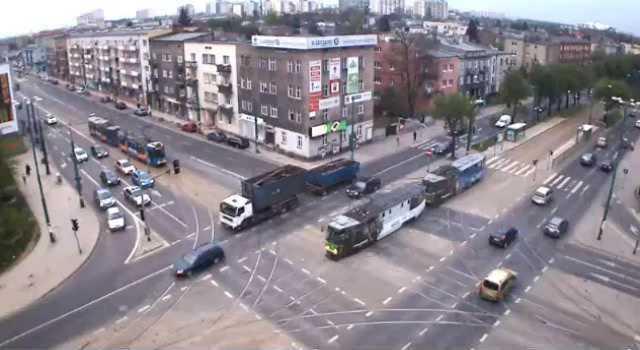}}\hfill
		\subfloat{\includegraphics[width=0.12\linewidth]{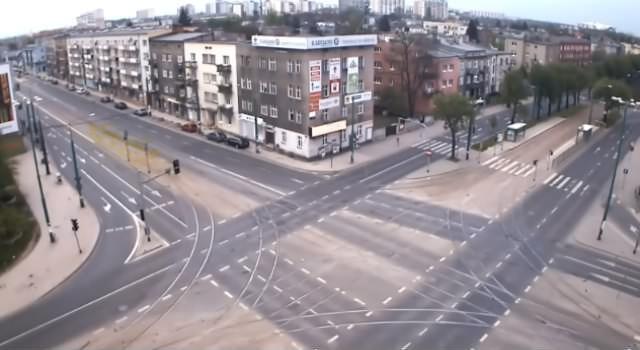}}\hfill
		\subfloat{\includegraphics[width=0.12\linewidth]{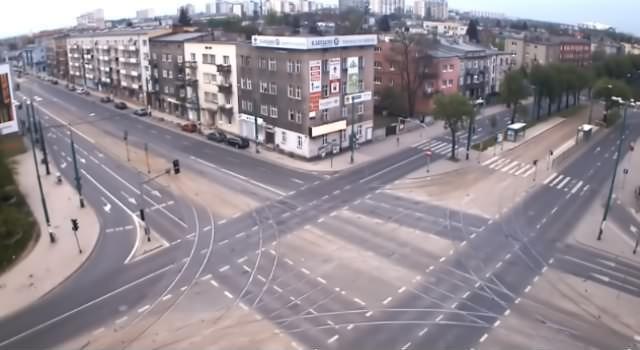}}\hfill
		\subfloat{\includegraphics[width=0.12\linewidth]{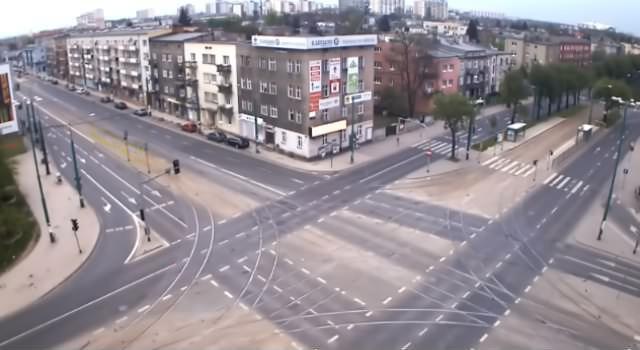}}\hfill
		\subfloat{\includegraphics[width=0.12\linewidth]{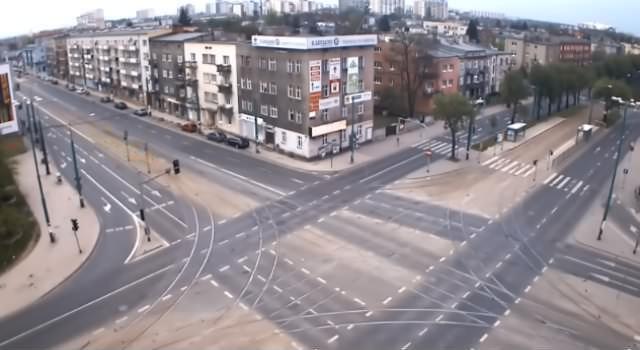}}\hfill
		\subfloat{\includegraphics[width=0.12\linewidth]{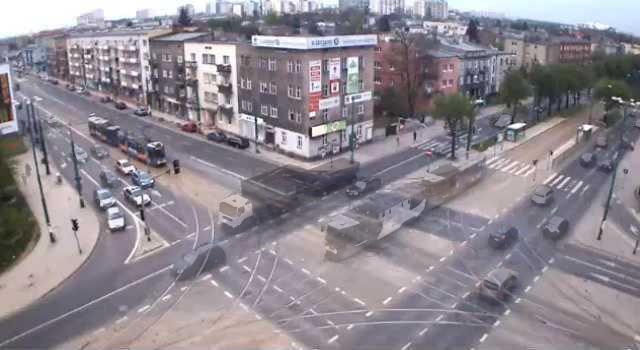}}\hfill
		\subfloat{\includegraphics[width=0.12\linewidth]{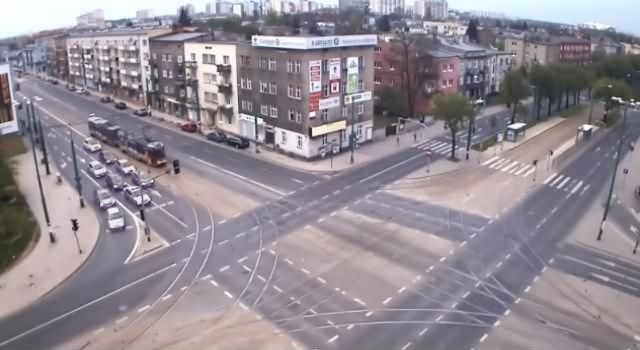}}\hfill
		\subfloat{\includegraphics[width=0.12\linewidth]{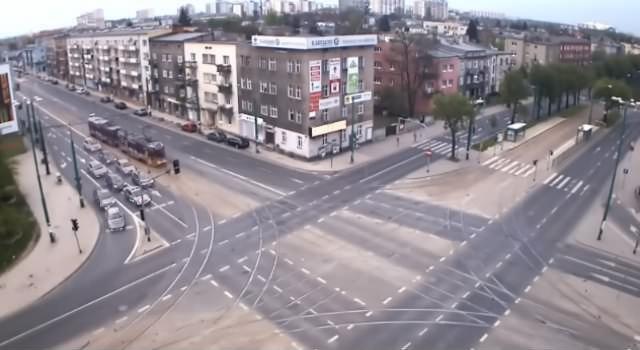}}
		\vspace{-0.12in}
		\\
		\subfloat{\includegraphics[width=0.12\linewidth]{video/video2_blank.jpg}}\hfill
		\subfloat{\includegraphics[width=0.12\linewidth]{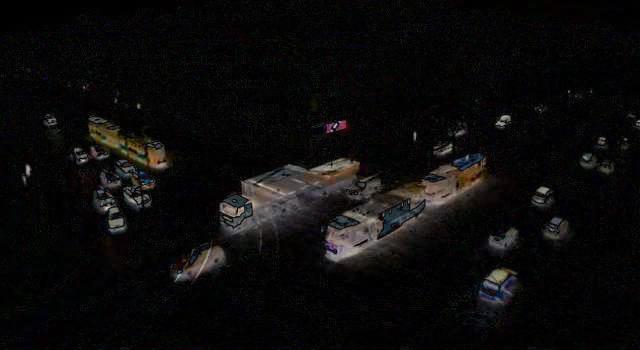}}\hfill
		\subfloat{\includegraphics[width=0.12\linewidth]{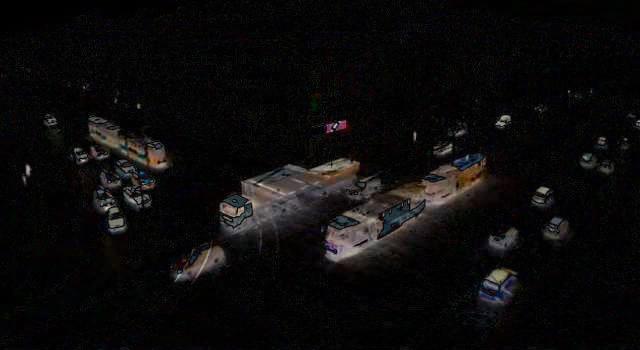}}\hfill
		\subfloat{\includegraphics[width=0.12\linewidth]{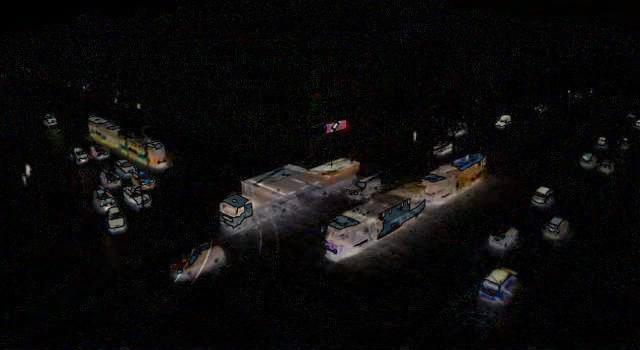}}\hfill
		\subfloat{\includegraphics[width=0.12\linewidth]{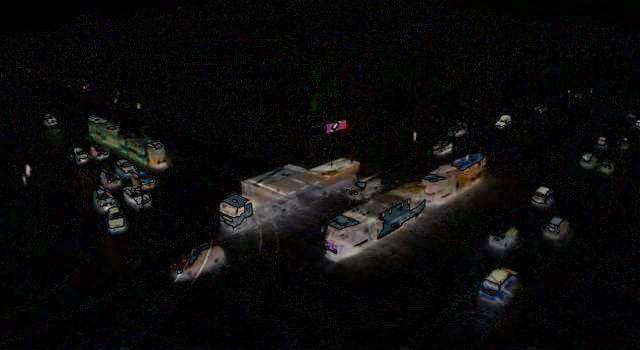}}\hfill
		\subfloat{\includegraphics[width=0.12\linewidth]{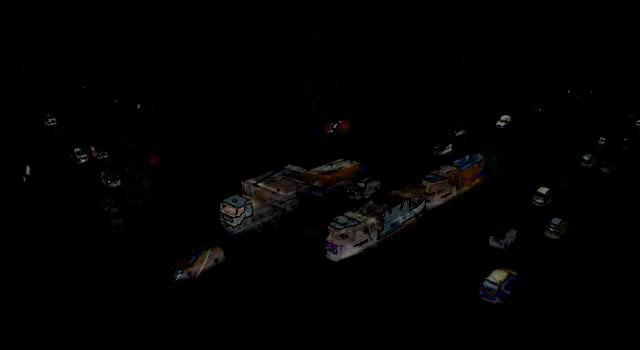}}\hfill
		\subfloat{\includegraphics[width=0.12\linewidth]{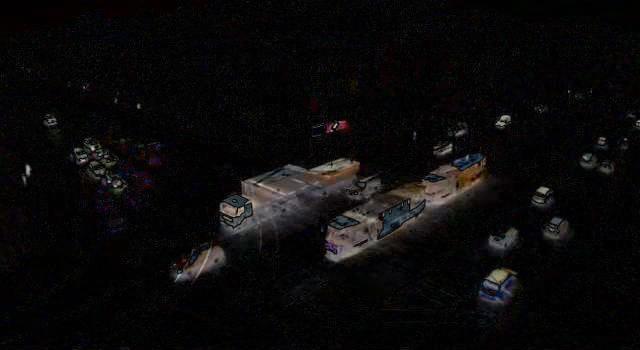}}\hfill
		\subfloat{\includegraphics[width=0.12\linewidth]{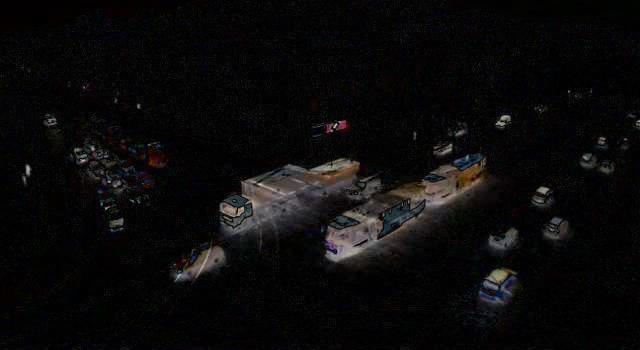}} 
		\\
		\vspace{-0.30in}
		\subfloat{\includegraphics[width=0.12\linewidth]{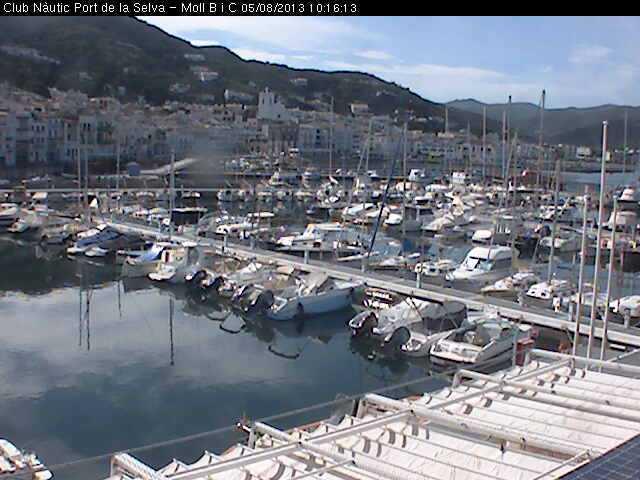}}\hfill
		\subfloat{\includegraphics[width=0.12\linewidth]{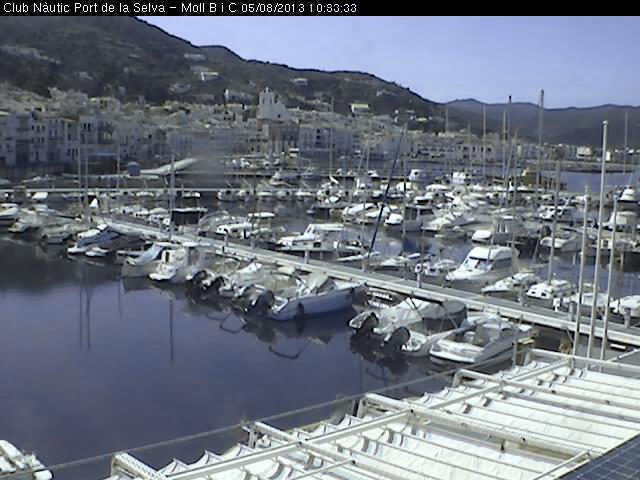}}\hfill
		\subfloat{\includegraphics[width=0.12\linewidth]{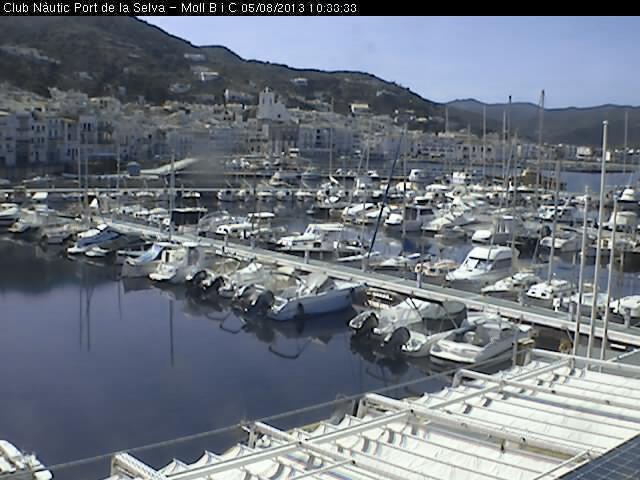}}\hfill
		\subfloat{\includegraphics[width=0.12\linewidth]{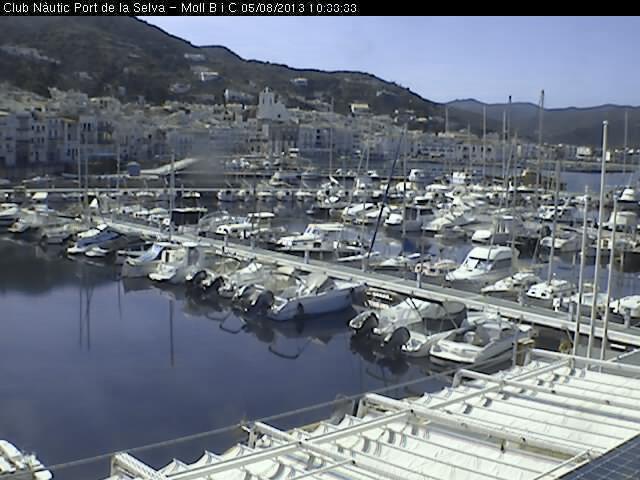}}\hfill
		\subfloat{\includegraphics[width=0.12\linewidth]{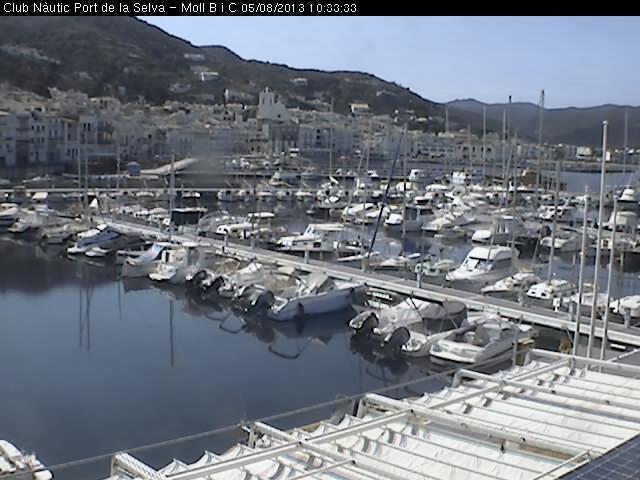}}\hfill
		\subfloat{\includegraphics[width=0.12\linewidth]{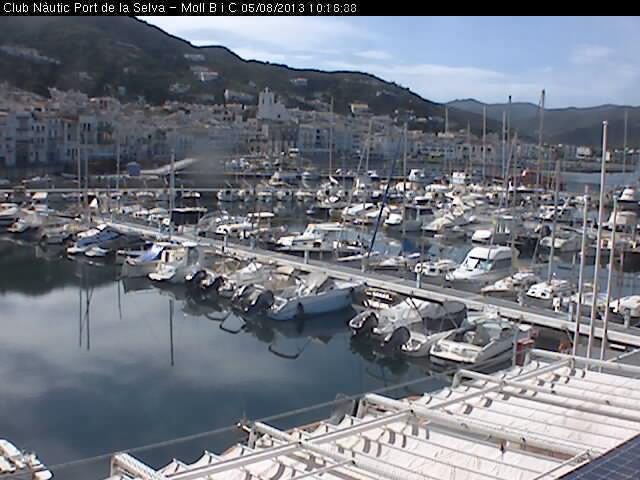}}\hfill
		\subfloat{\includegraphics[width=0.12\linewidth]{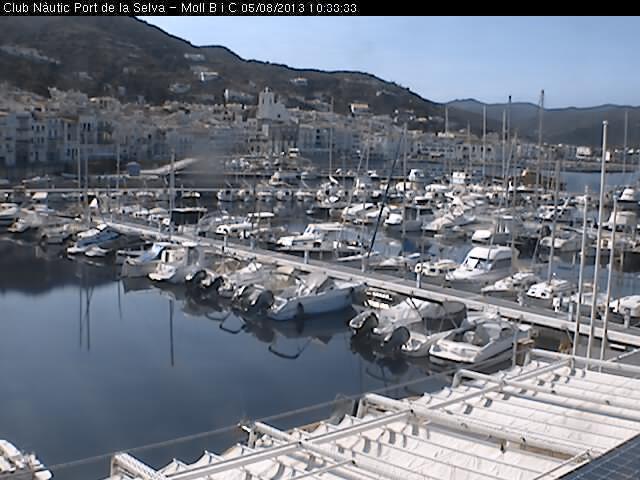}}\hfill
		\subfloat{\includegraphics[width=0.12\linewidth]{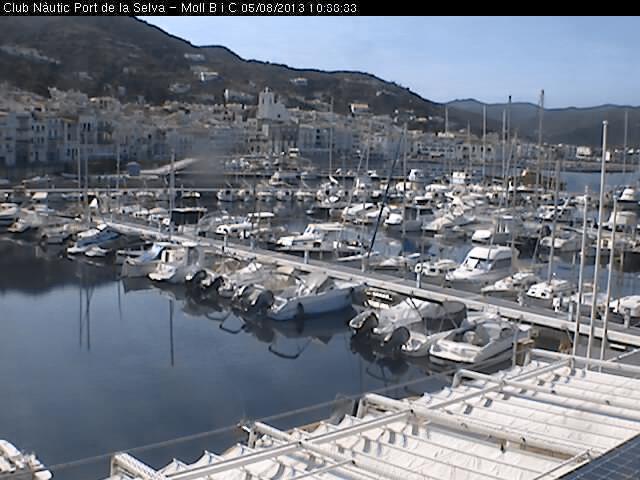}}
		\vspace{-0.12in}
		\\
		\subfloat{\includegraphics[width=0.12\linewidth]{video/video2_blank.jpg}}\hfill
		\subfloat{\includegraphics[width=0.12\linewidth]{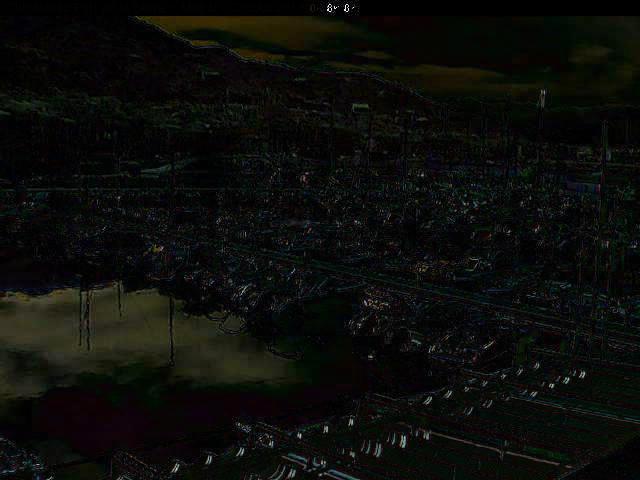}}\hfill
		\subfloat{\includegraphics[width=0.12\linewidth]{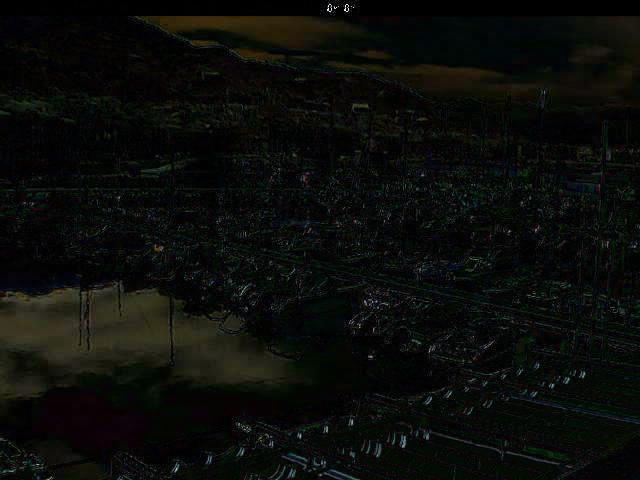}}\hfill
		\subfloat{\includegraphics[width=0.12\linewidth]{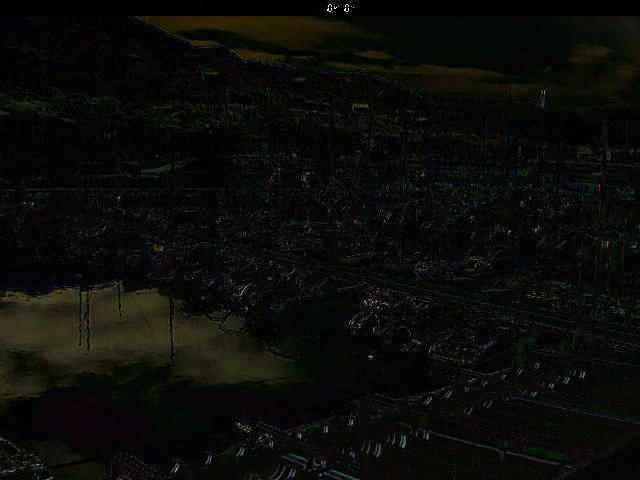}}\hfill
		\subfloat{\includegraphics[width=0.12\linewidth]{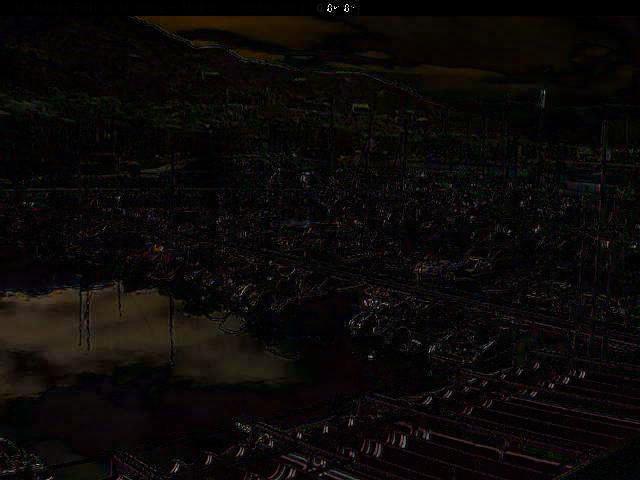}}\hfill
		\subfloat{\includegraphics[width=0.12\linewidth]{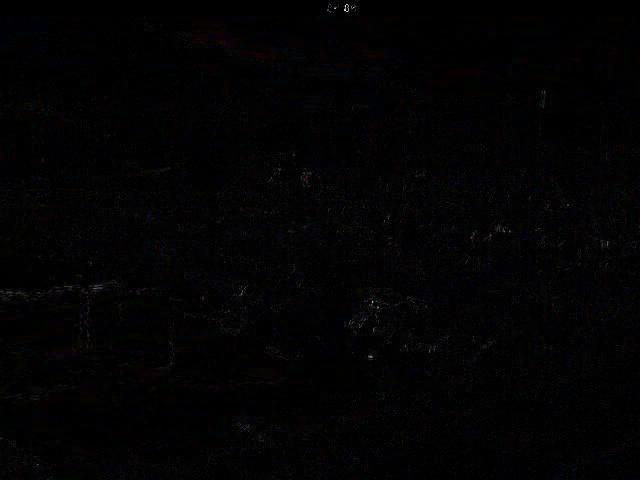}}\hfill
		\subfloat{\includegraphics[width=0.12\linewidth]{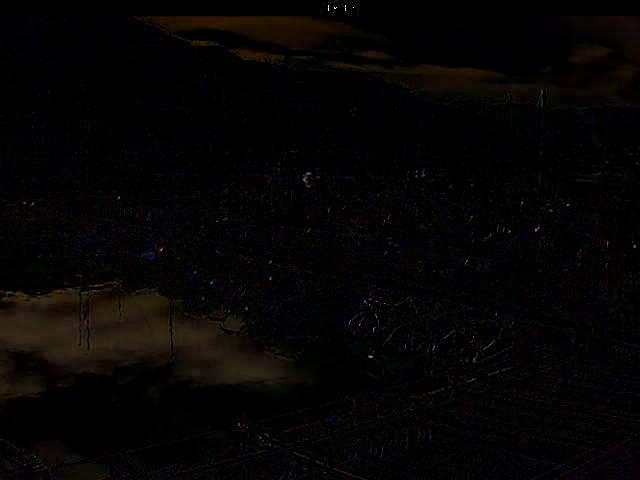}}\hfill
		\subfloat{\includegraphics[width=0.12\linewidth]{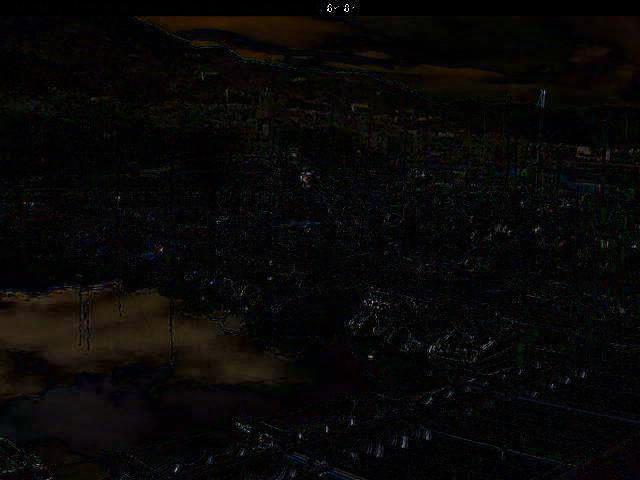}} 
		\vspace{-0.12in}
		
		\subfloat{\includegraphics[width=0.12\linewidth]{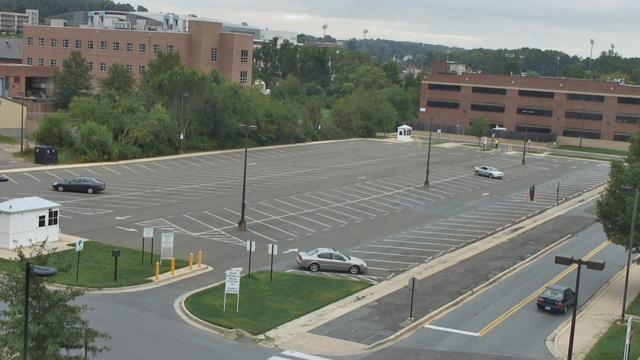}}\hfill
		\subfloat{\includegraphics[width=0.12\linewidth]{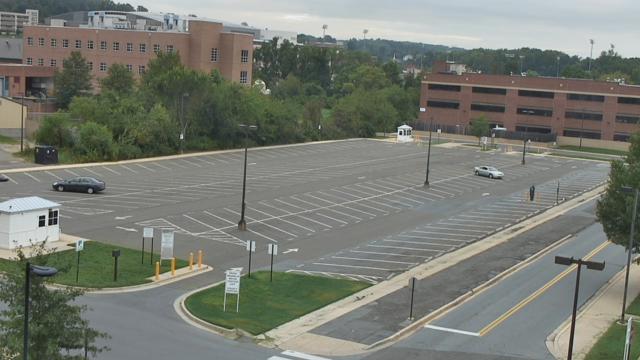}}\hfill
		\subfloat{\includegraphics[width=0.12\linewidth]{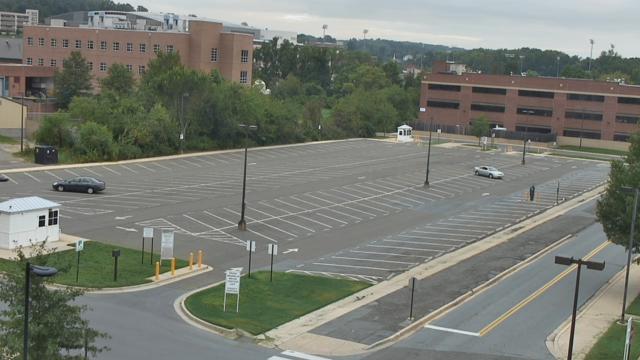}}\hfill
		\subfloat{\includegraphics[width=0.12\linewidth]{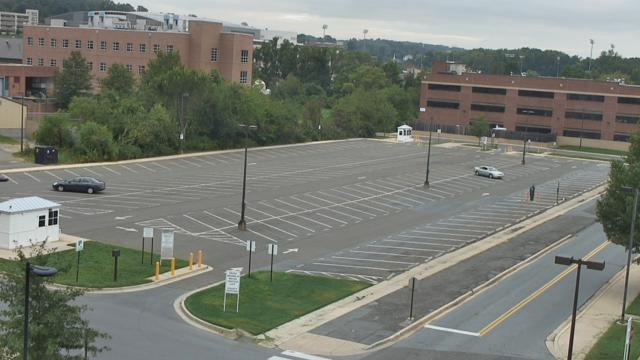}}\hfill
		\subfloat{\includegraphics[width=0.12\linewidth]{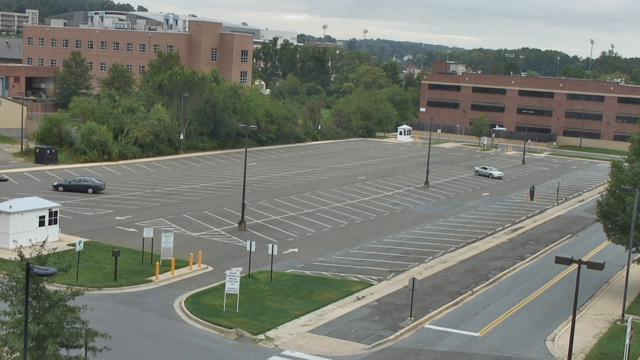}}\hfill
		\subfloat{\includegraphics[width=0.12\linewidth]{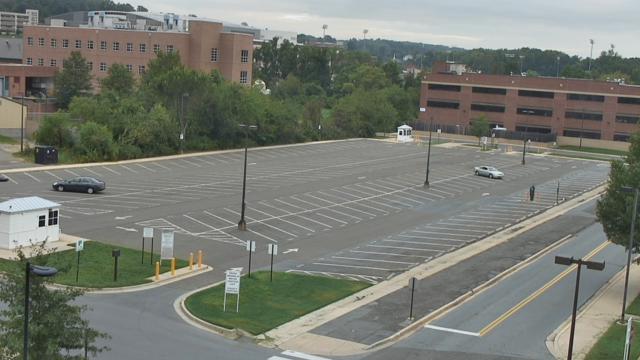}}\hfill
		\subfloat{\includegraphics[width=0.12\linewidth]{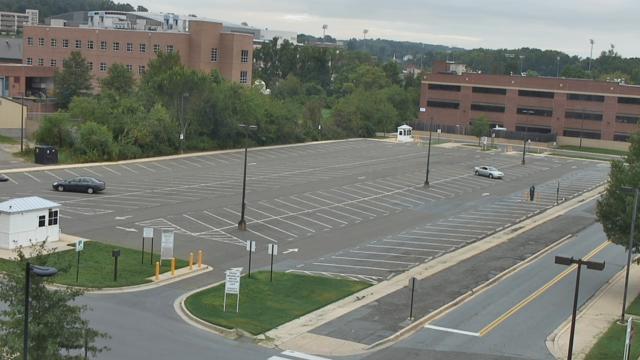}}\hfill
		\subfloat{\includegraphics[width=0.12\linewidth]{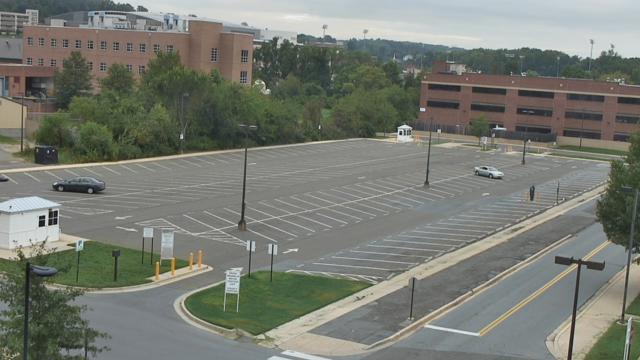}}
		\vspace{-0.12in}
		\\
		\subfloat{\includegraphics[width=0.12\linewidth]{video/video2_blank.jpg}}\hfill
		\subfloat{\includegraphics[width=0.12\linewidth]{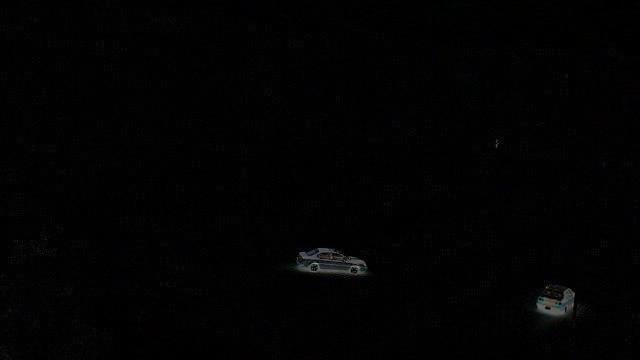}}\hfill
		\subfloat{\includegraphics[width=0.12\linewidth]{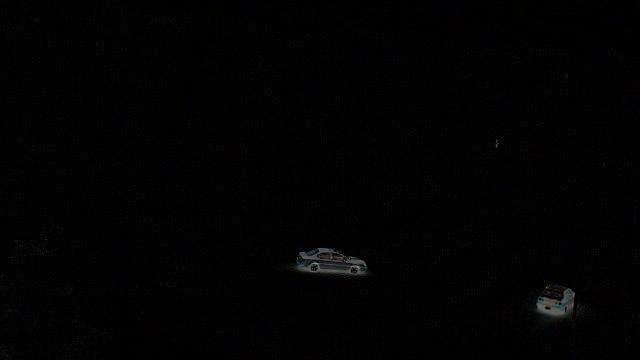}}\hfill
		\subfloat{\includegraphics[width=0.12\linewidth]{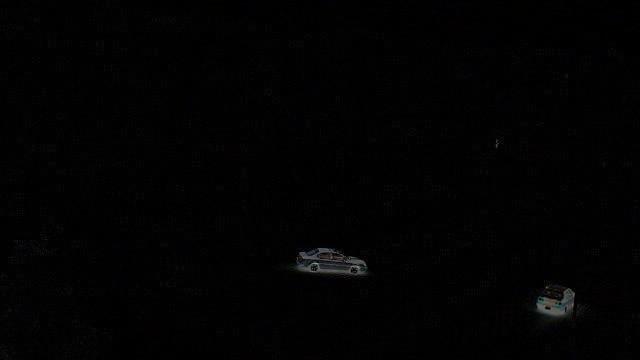}}\hfill
		\subfloat{\includegraphics[width=0.12\linewidth]{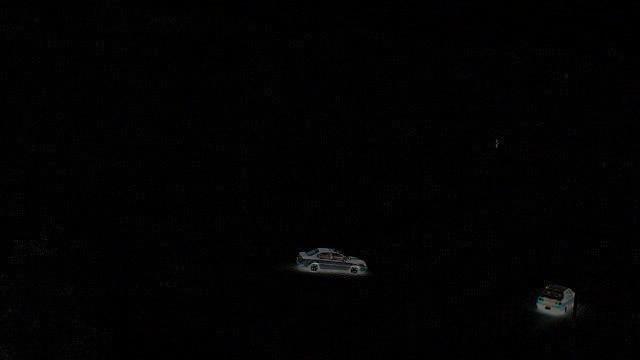}}\hfill
		\subfloat{\includegraphics[width=0.12\linewidth]{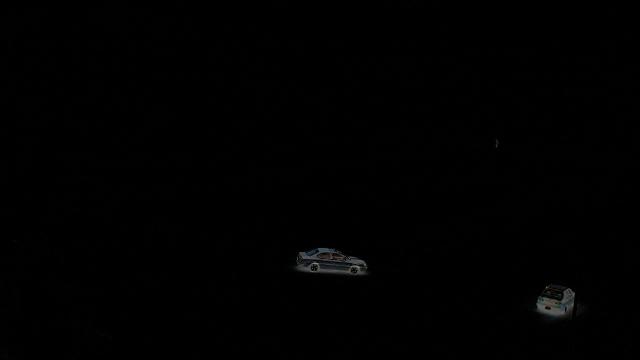}}\hfill
		\subfloat{\includegraphics[width=0.12\linewidth]{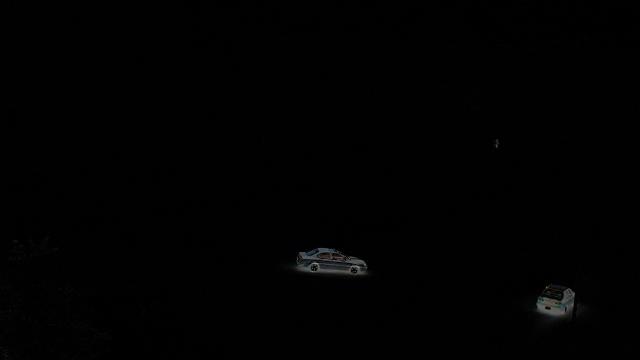}}\hfill
		\subfloat{\includegraphics[width=0.12\linewidth]{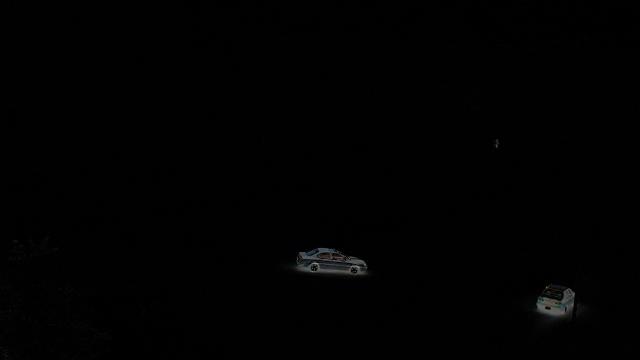}} 
		\vspace{-0.25in}
		\caption{Visual results for color video background subtraction. The \textbf{first two rows} are separated   backgrounds and foregrounds corresponding to a frame from \textit{Shoppingmall}, the \textbf{3rd and 4th rows}  are separated  backgrounds and foregrounds corresponding to a frame from \textit{Highway}, the \textbf{5th and 6th} are separated  backgrounds and foregrounds corresponding to a frame from \textit{Crossroad}, \textbf{7th and 8th} correspond to a frame from \textit{Port}, and the \textbf{last two rows} correspond to a frame from \textit{Parking lot}, except the first column which is the original frame.
		}\label{FIG: video}
		\vspace{-0.1in}
	\end{figure}
	
	We apply the four variants of \alg\ and aforementioned tensor/matrix RPCA algorithms on the color video background subtraction task. We obtain 5 color video datasets from various sources: \textit{Shoppingmall} \cite{l}, \textit{Highway} \cite{bouwmans2017scene}, \textit{Crossroad} \cite{bouwmans2017scene}, \textit{Port} \cite{bouwmans2017scene}, and \textit{Parking-lot} \cite{oh2011large}. 
	Similar to \Cref{sec:face}, we vectorize each frame and construct a $(\text{height} \cdot \text{width})\times 3\times \text{frame}$ data tensor. 
	The tensor is unfolded to a (height $\cdot$ width) $\times$ ($3~\cdot$ frame) matrix for matrix RPCA methods. We use Tucker rank $(3,3,3)$ for tensors methods and rank 3 for matrix methods. 
	We exclude RGD from this experiment because the disk space required for RGD exceeds our server limit.
	Among the tested videos, \textit{Shoppingmall}, \textit{Highway}, and \textit{Parking-lot} are normal speed videos with almost static backgrounds. \textit{Crossroad} and \textit{Port} are outdoor time-lapse videos, hence their background colors change slightly between different frames. We observe that all tested algorithms perform very similarly for videos with static backgrounds and produce visually desirable output. On the other hand, the color of the extracted background varies slightly among different algorithms on time-lapse videos. Since the color of the background keeps changing slightly for the time-lapse videos, we cannot decide the ground-truth color of the background, hence we do not rank the performance of different algorithms. The runtime for results along with video size information are summarized in \Cref{table:video}. By comparing the runtime of four variants of \alg, we can observe that the experiment result generally agrees with the analysis on computational efficiency in \Cref{sec:4vars}. All \alg\ variants accomplish the background subtraction task faster than the guideline methods. In addition, we provide some selected visual results in \Cref{FIG: video}.

	\subsection{Network Clustering} \label{sec:network}
	In this section, we apply our \alg\ algorithm, and the TRPCA algorithm RGD from \cite{cai2022generalized}  for the community detection task on the co-authorship network data from \cite{ji2016coauthorship} and compare their results and efficiency. This dataset contains 3248 papers with a list of authors for each paper (3607 authors in total); hence could naturally serve as the adjacency matrix for the weighted co-authorship graph. The original paper for this dataset tests a number of community detection algorithms, including network spectral clustering, profile likelihood, pseudo-likelihood approach, etc., on a selected subset with 236 authors and 542 papers. To obtain this subset, we generate an undirected graph for the 3607 authors and put an edge between two authors if   they have co-authored two or more papers. We then take the largest connected component, and all the nodes in this component are the subset of authors we are interested in.  
	Since we aim to explore the higher-order interactions among this co-authorship network, we convert this connected component into a 3-mode adjacency tensor $\tens{T}$ and apply TRPCA algorithms to obtain the low-rank component from $\tens{L}$. We construct the adjacency tensor $\tens{T}$ with the following rules:
	\begin{itemize}
		\item For any two connected authors $(i,j)$, which means author $i$ and $j$ have worked together for at least two papers, we set $\tens{T}_{\mathfrak{G}(i,i,j)} = \tens{T}_{\mathfrak{G}(i,j,j)} = 1$ 
		\item For any three pairwisely connected authors $(i,j,k)$, we set $\tens{T}_{\mathfrak{G}(i,j,k)} = 1$. Notice that these three authors may not appear in one paper at the same time, but each pair of them have worked together for at least two papers.
	\end{itemize}
	
\begin{figure}[ht]
		\centering
		\subfloat{\includegraphics[trim=100 40 100 50,clip,width =
			0.49\linewidth]{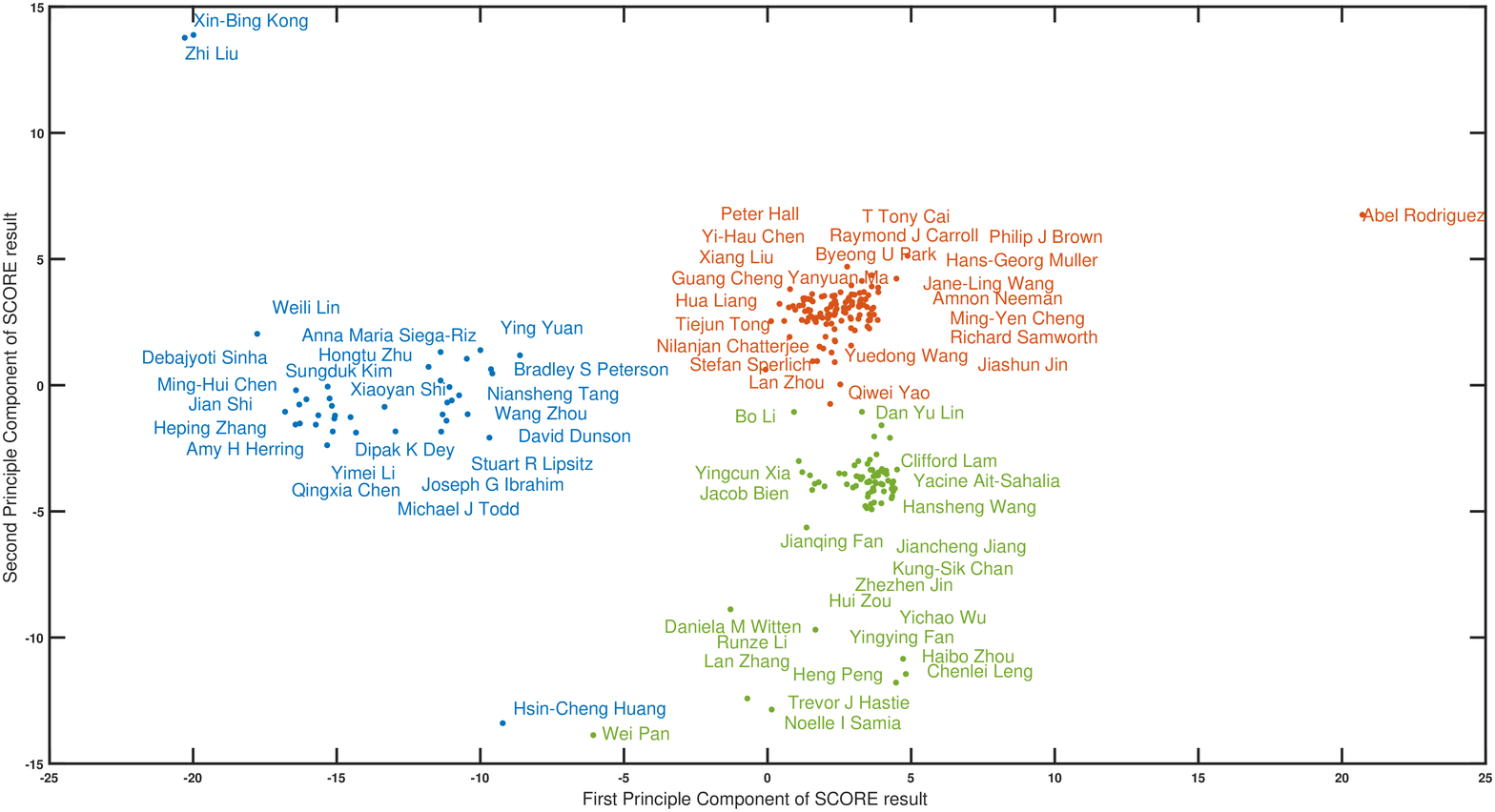}}\hfill
		\subfloat{\includegraphics[trim=100 40 100 50,clip,width = 0.49\linewidth]{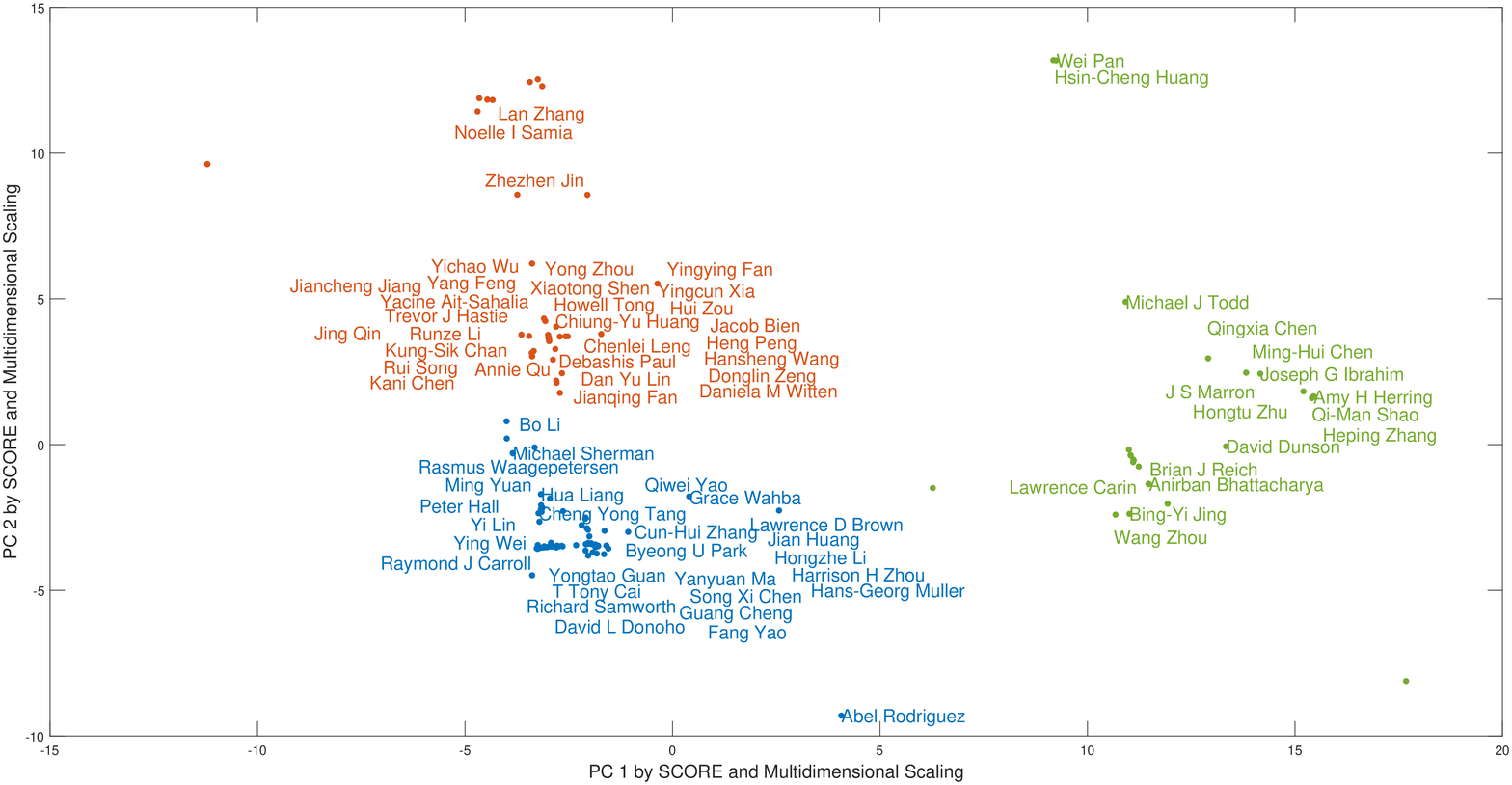}}\\
		\subfloat{\includegraphics[trim=100 40 100 50,clip,width = 0.49\linewidth]{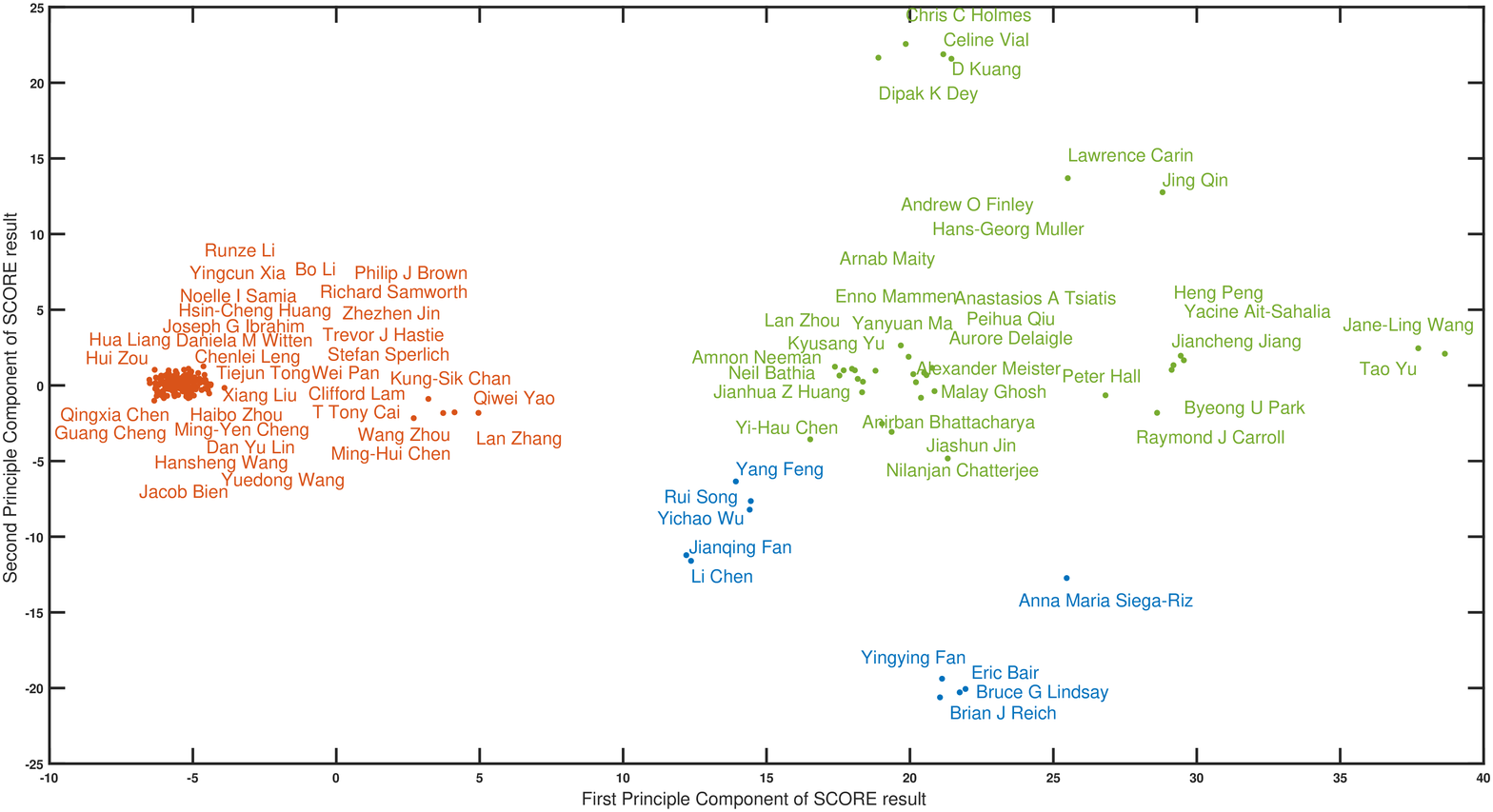}}\hfill
		\subfloat{\includegraphics[trim=100 40 100 50,clip,width = 0.49\linewidth]{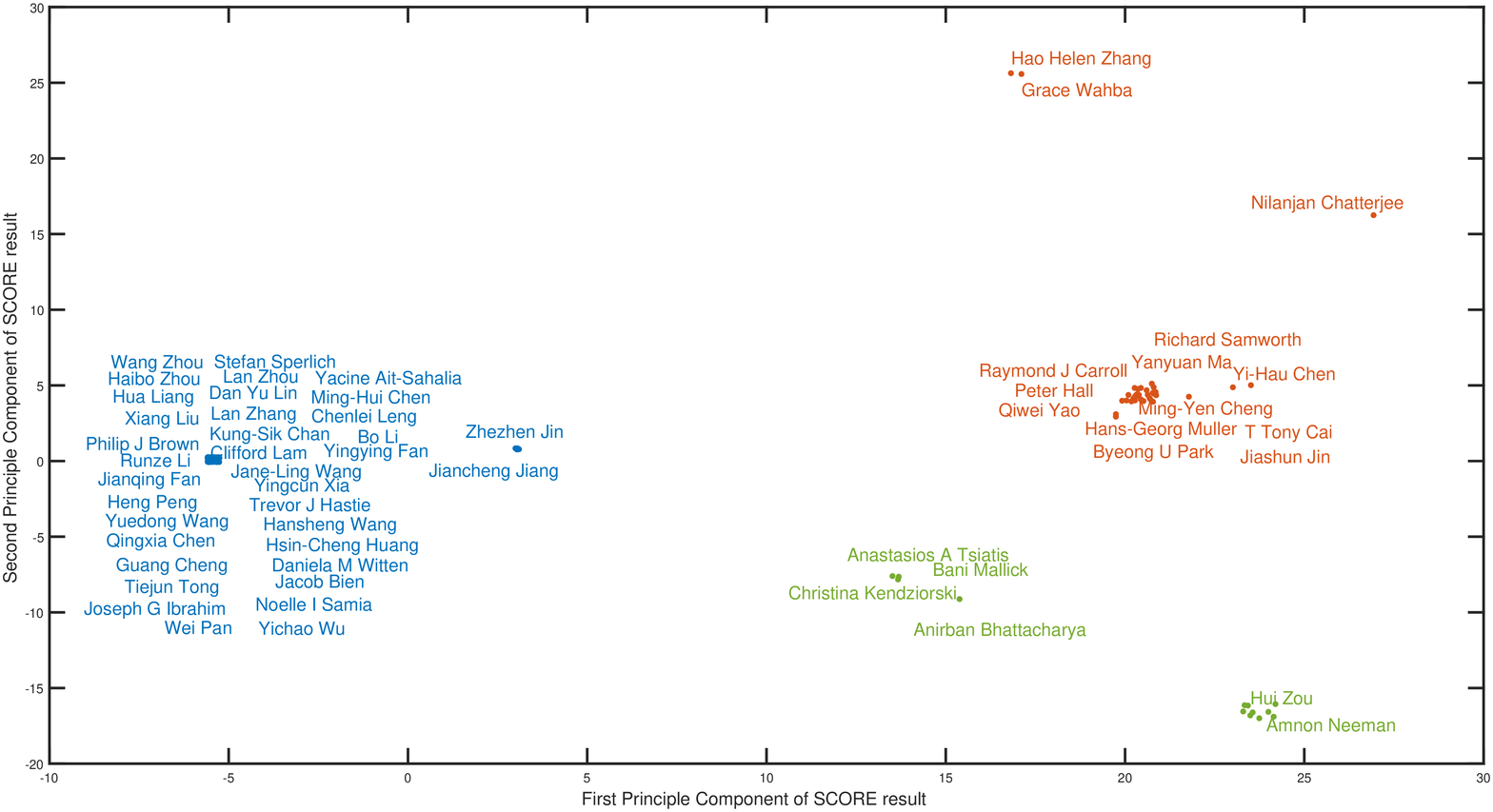}}\\
		\subfloat{\includegraphics[trim=100 40 100 50,clip,width = 0.49\linewidth]{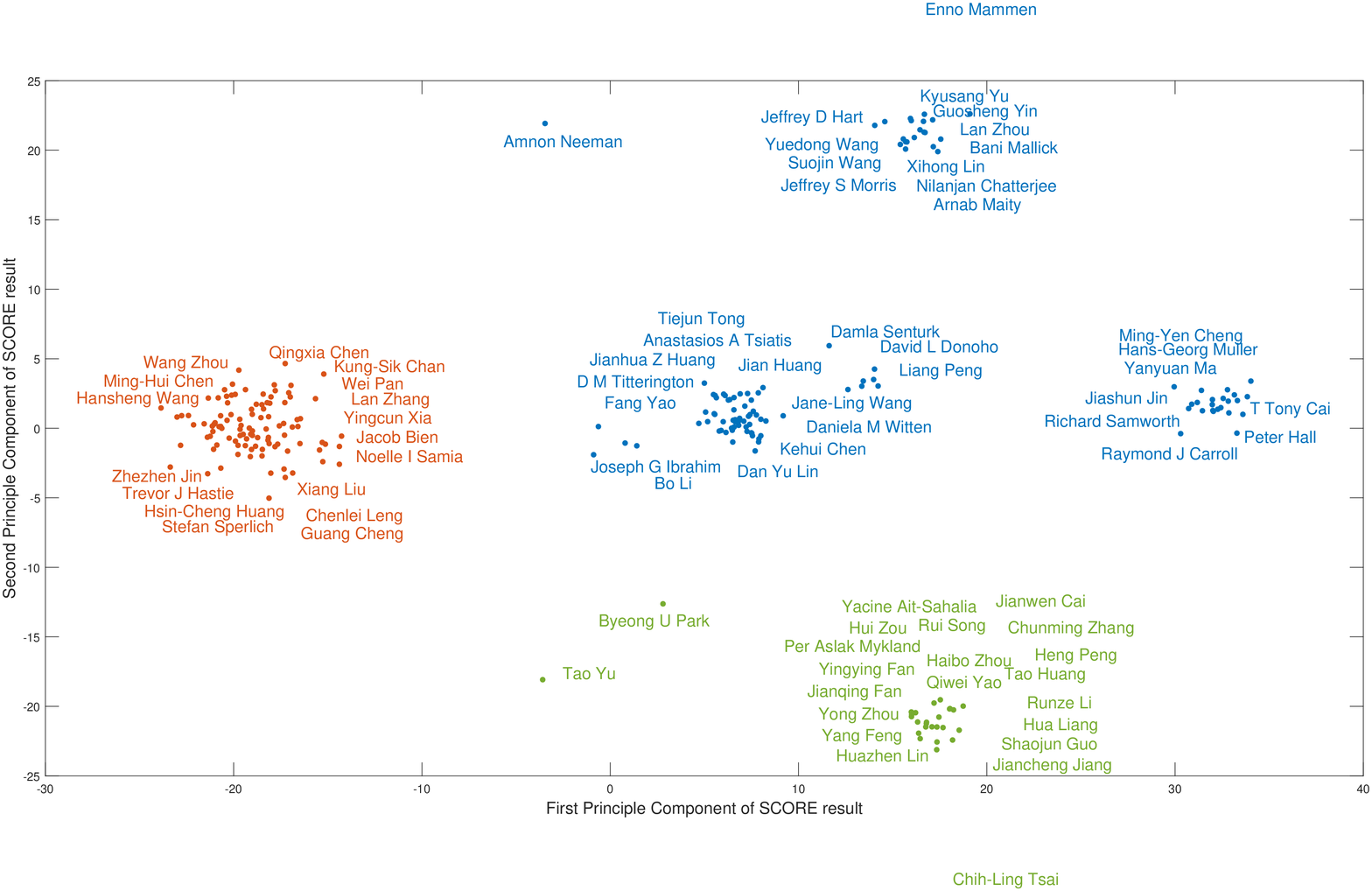}}\hfill
		\subfloat{\includegraphics[trim=100 40 100 50,clip,width = 0.49\linewidth]{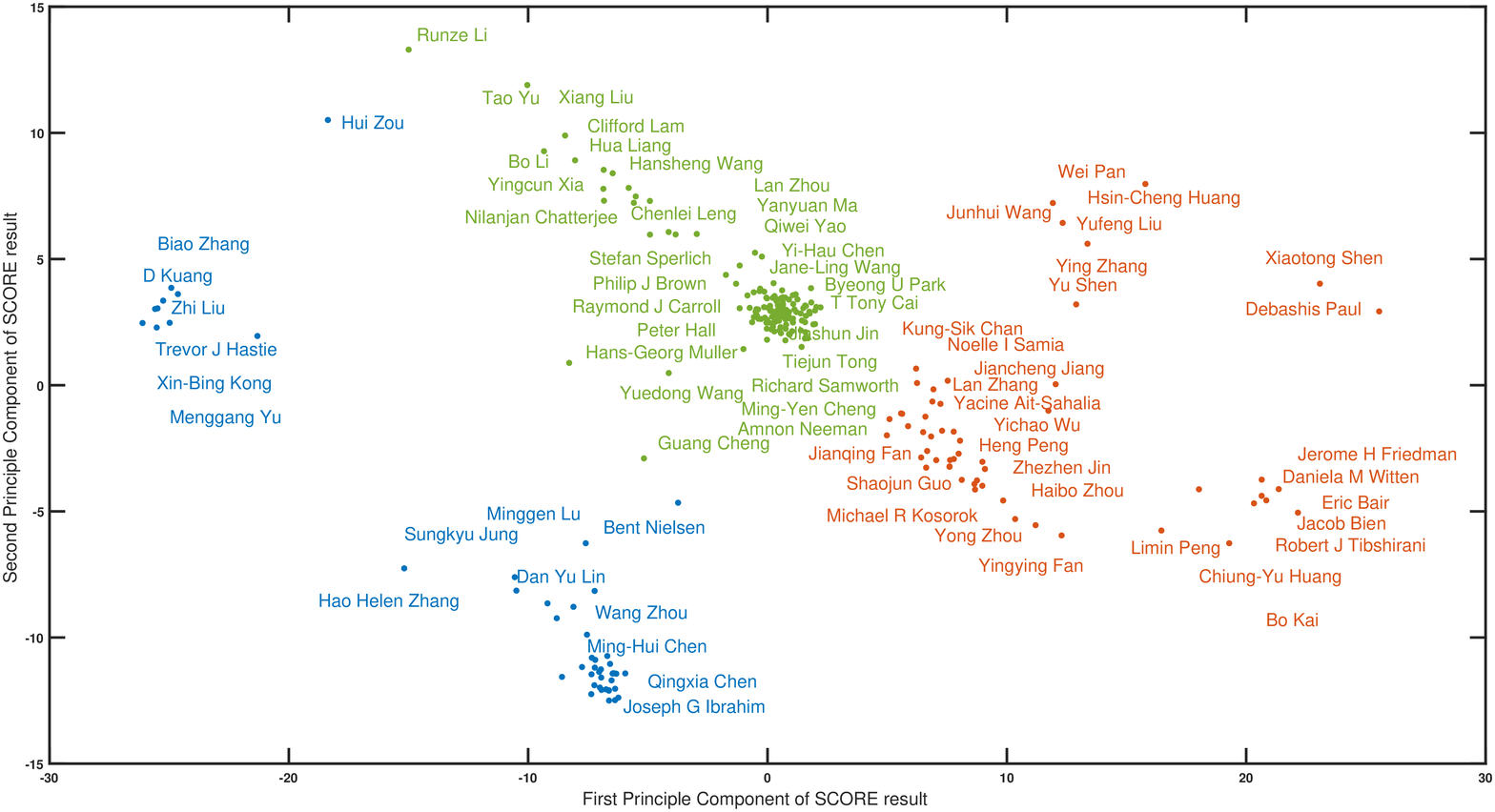}}
		\caption{ 
			Three communities detected in the “High-Dimensional Data Analysis” co-authorship network with SCORE\cite{ji2016coauthorship}, RGD \cite{cai2022generalized}, and \alg\ . \textbf{Top left}: Result from SCORE on original tensor; \textbf{top right}: Result from RGD and SCORE; \textbf{middle left}: Result from \alg-FF and SCORE, with $\upsilon$ = 6; \textbf{middle right}: Result from \alg-FC and SCORE, with $\upsilon$ = 2; \textbf{bottom left}: Result from \alg-FC and SCORE, with $\upsilon$ = 6; \textbf{bottom right}: Result from \alg-FC and SCORE, with $\upsilon$ = 11. All TRPCA methods are applied on the adjacency tensor $\tens{T}$ with Tucker rank $= (4,4,4)$.} 
		\label{fig:clusters}
		\vspace{-0.1in}
	\end{figure}

	Here $\mathfrak{G}(S)$ denotes all permutations of the set $S$. Therefore the adjacency tensor $\tens{T}$ is symmetric. Now we apply \alg\ with different sampling constants as well as the TRPCA algorithm RGD \cite{cai2022generalized} to learn the low-rank component $\tens{L}$ with Tucker rank $(4,4,4)$, which is used to infer the communities in this network. Then we apply the SCORE algorithm \cite{ji2016coauthorship} as the clustering algorithm on the low-rank tensor $\tens{L}$. We use SCORE instead of other traditional clustering algorithms such as spectral clustering because SCORE could mitigate the influence of node heterogeneity \cite{ji2016coauthorship}. We plot the results from each TRPCA algorithm in \Cref{fig:clusters}.

	The clustering of the “High-Dimensional Data Analysis” co-authorship network is an unsupervised task, which means the ground truth of labeling an author with a certain community does not exist. Therefore, we do not focus on qualitatively evaluating each result, but we present the new findings from higher-order interactions among co-authors and analyze the results from different choices of parameters. Previous studies on this co-authorship network generally provide three clusters with names: “Carroll–Hall” group, “North Carolina” community, and “Fan and others” group \cite{ji2016coauthorship, cai2022generalized}. Among them, “Carroll-Hall” group generally includes researchers in nonparametric and semi-parametric statistics, functional estimation, and high-dimensional statistics; “North Carolina” group generally includes researchers from Duke University, University of North Carolina, and North Carolina State University; “Fan and Others” group includes primarily the researchers collaborating closely with Jianqing Fan or his co-authors, and other researchers who do not obviously belong to the first two groups \cite{cai2022generalized}. For conciseness, we will make use of the same name of each group as in previous studies on this co-authorship network.
	
	The top two plots of \Cref{fig:clusters} are existing results from \cite{cai2022generalized}. With SCORE as the clustering method, the original tensor and the RGD output 
	both successfully reveal the two groups: the “Carroll–Hall” group and a “North Carolina” community, with slightly different clustering result for researchers who do not obviously belong to one group, such as Debajyoti Sinha, Michael J Todd, and Abel Rodriguez. One can observe that Fiber \alg\ detected the “Carroll–Hall” group. However, Fiber \alg\ labels most authors not having a strong connection with Peter Hall, Raymond Carroll, and Jianqing Fan as the “North Carolina” community. Similarly, Chidori \alg\ with $\upsilon = 2$ generates the center of the “Carroll–Hall” group as one cluster and categorizes most authors with a lower number of co-authorships and not co-authored with kernel members into the “Fan and others” group. We infer that the tendency to cluster most members into one group is due to insufficient sampling. The co-authorship tensor is very sparse, with only about 2\% entries being non-zero, so the feature of each node may not be sufficiently extracted from Fiber sampling or Chidori sampling with small $\upsilon$. From the middle two plots, we can observe that most authors in the largest group have very close first and second principal components in the two-dimension embedding, providing the evidence that the algorithm ignored some non-zeros entries for nodes with fewer numbers of connections during the sampling process.
	
	Note that the sampling constant of Fiber sampling should be $r\log d$  
	times the constant of Chidori sampling in order to access the same amount of data from the original tensor, where $d$ denotes the number of authors in this experiment). So we only test the Chidori sampling with a larger sampling constant on the co-authorship tensor $\tens{L}$ for efficiency. The result is shown in the bottom: $\upsilon = 6$ for bottom-left and $\upsilon = 11$ for bottom-right of \Cref{fig:clusters}. Both settings generate the “Carroll–Hall” group with authors having strong ties to Peter Hall and Raymond Carroll, such as Richard Samworth, Hans-Georg Muller, Anastasios Tsiatis, Yuanyuan Ma, Yuedong Wang, Lan Zhou, etc. The “Fan and others” is also successfully detected, including co-authors of high-degree node Jianqing Fan such as Hua Liang, Yingying Fan, Haibo Zhou, Yong Zhou, Jiancheng Jiang, Qiwei Yao, etc. The sizes of the three clusters generated from these two settings are more balanced than the result from \alg\ with smaller $\upsilon$. Therefore we can conclude that, in this real-world network clustering task, different choices of sampling constant provide the same core members of each group, and the group size is more balanced with a larger amount of sampling data at the cost of computation efficiency.  \Cref{tab:network time} shows the runtime for each algorithm and sampling method.
	
	\vspace{-5mm}
	\begin{table}[h]
		\caption{Runtime comparison (in seconds) of TRPCA algorithms: RGD and \alg.}
		\label{tab:network time}
		\centering
		\vspace{-0.1in}
		\begin{small}
			\begin{tabular}{c|c}\toprule
				\textsc{Method} & \textsc{Runtime}\\\midrule
				RGD \cite{cai2022generalized}& 120.5 \\
				\alg-FF with $\upsilon = 6$~\,  & 3.526 \\
				\alg-FC with $\upsilon = 2$~\,  & 0.571 \\
				\alg-FC with $\upsilon = 6$~\,  & 4.319 \\
				\alg-FC with $\upsilon = 11$  & 7.177 \\\bottomrule
			\end{tabular}
		\end{small}
	\end{table}
	
	\vspace{-5mm}
	
	\section{Conclusion} 
	{This paper presents a highly efficient algorithm \alg\ for large-scale TRPCA problems. \alg\ is developed by introducing a novel inexact low-Tucker-rank tensor approximation via  mode-wise tensor decomposition. The structure of this approximation significantly reduces the computational complexity, resulting in a computational complexity for each iteration of $\cO(n^2d(r\log d)^2)$ with Fiber CUR and $\cO(nd(r\log d)^n)$ with Chidori CUR. This is much lower than  the minimum of $\cO(rd^n)$ required by HOSVD-based algorithms. Numerical experiments on synthetic and real-world datasets also demonstrate the efficiency advantage of \alg\ against other state-of-the-art tensor/matrix RPCA algorithms. Additionally, the fixed sampling variants of \alg\ only require partial information from the input tensor for the TRPCA task. }
	
	\bibliographystyle{siamplain}
	\bibliography{reference}

\begin{thebibliography}{10}

\bibitem{abdi2010principal}
{\sc H.~Abdi and L.~J. Williams}, {\em Principal component analysis}, Wiley
  interdisciplinary reviews: computational statistics, 2 (2010), pp.~433--459.

\bibitem{bergqvist2010higher}
{\sc G.~Bergqvist and E.~G. Larsson}, {\em The higher-order singular value
  decomposition: Theory and an application [lecture notes]}, IEEE signal
  processing magazine, 27 (2010), pp.~151--154.

\bibitem{bouwmans2018applications}
{\sc T.~Bouwmans, S.~Javed, H.~Zhang, Z.~Lin, and R.~Otazo}, {\em On the
  applications of robust {PCA} in image and video processing}, Proceedings of
  the IEEE, 106 (2018), pp.~1427--1457.

\bibitem{bouwmans2017scene}
{\sc T.~Bouwmans, L.~Maddalena, and A.~Petrosino}, {\em Scene background
  initialization: A taxonomy}, Pattern Recognition Letters, 96 (2017),
  pp.~3--11.

\bibitem{cai2021accelerated}
{\sc H.~Cai, J.-F. Cai, T.~Wang, and G.~Yin}, {\em Accelerated structured
  alternating projections for robust spectrally sparse signal recovery}, IEEE
  Transactions on Signal Processing, 69 (2021), pp.~809--821.

\bibitem{cai2019accelerated}
{\sc H.~Cai, J.-F. Cai, and K.~Wei}, {\em Accelerated alternating projections
  for robust principal component analysis}, The Journal of Machine Learning
  Research, 20 (2019), pp.~685--717.

\bibitem{cai2022structured}
{\sc H.~Cai, J.-F. Cai, and J.~You}, {\em Structured gradient descent for fast
  robust low-rank {H}ankel matrix completion}, SIAM Journal on Scientific
  Computing, 45 (2023), pp.~1172--1198.

\bibitem{cai2021fast}
{\sc H.~Cai, Z.~Chao, L.~Huang, and D.~Needell}, {\em Fast robust tensor
  principal component analysis via {F}iber {CUR} decomposition}, in Proceedings
  of the IEEE/CVF International Conference on Computer Vision Workshops, 2021,
  pp.~189--197.

\bibitem{cai2020rapid}
{\sc H.~Cai, K.~Hamm, L.~Huang, J.~Li, and T.~Wang}, {\em Rapid robust
  principal component analysis: {CUR} accelerated inexact low rank estimation},
  IEEE Signal Processing Letters, 28 (2020), pp.~116--120.

\bibitem{cai2021mode}
{\sc H.~Cai, K.~Hamm, L.~Huang, and D.~Needell}, {\em Mode-wise tensor
  decompositions: Multi-dimensional generalizations of {CUR} decompositions},
  The Journal of Machine Learning Research, 22 (2021), pp.~1--36.

\bibitem{cai2021robust}
{\sc H.~Cai, K.~Hamm, L.~Huang, and D.~Needell}, {\em Robust {CUR}
  decomposition: Theory and imaging applications}, SIAM Journal on Imaging
  Sciences, 14 (2021), pp.~1472--1503.

\bibitem{cai2021learned}
{\sc H.~Cai, J.~Liu, and W.~Yin}, {\em Learned robust {PCA}: A scalable deep
  unfolding approach for high-dimensional outlier detection}, Advances in
  Neural Information Processing Systems, 34 (2021).

\bibitem{cai2022generalized}
{\sc J.-F. Cai, J.~Li, and D.~Xia}, {\em Generalized low-rank plus sparse
  tensor estimation by fast riemannian optimization}, Journal of the American
  Statistical Association,  (2022), pp.~1--17.

\bibitem{caiafa2010generalizing}
{\sc C.~F. Caiafa and A.~Cichocki}, {\em Generalizing the column--row matrix
  decomposition to multi-way arrays}, Linear Algebra and its Applications, 433
  (2010), pp.~557--573.

\bibitem{candes2007sparsity}
{\sc E.~Candes and J.~Romberg}, {\em Sparsity and incoherence in compressive
  sampling}, Inverse {P}roblems, 23 (2007), p.~969.

\bibitem{candes2011robust}
{\sc E.~J. Cand{\`e}s, X.~Li, Y.~Ma, and J.~Wright}, {\em Robust principal
  component analysis?}, Journal of the ACM (JACM), 58 (2011), pp.~1--37.

\bibitem{carroll1970analysis}
{\sc J.~D. Carroll and J.-J. Chang}, {\em Analysis of individual differences in
  multidimensional scaling via an n-way generalization of “eckart-young”
  decomposition}, Psychometrika, 35 (1970), pp.~283--319.

\bibitem{chambua2018tensor}
{\sc J.~Chambua, Z.~Niu, A.~Yousif, and J.~Mbelwa}, {\em Tensor factorization
  method based on review text semantic similarity for rating prediction},
  Expert Systems with Applications, 114 (2018), pp.~629--638.

\bibitem{chandrasekaran2011rank}
{\sc V.~Chandrasekaran, S.~Sanghavi, P.~A. Parrilo, and A.~S. Willsky}, {\em
  Rank-sparsity incoherence for matrix decomposition}, SIAM Journal on
  Optimization, 21 (2011), pp.~572--596.

\bibitem{chao2021hosvd}
{\sc Z.~Chao, L.~Huang, and D.~Needell}, {\em {HOSVD}-based algorithm for
  weighted tensor completion}, Journal of Imaging, 7 (2021), p.~110.

\bibitem{chen2009tensor}
{\sc J.~Chen and Y.~Saad}, {\em On the tensor {SVD} and the optimal low rank
  orthogonal approximation of tensors}, SIAM journal on Matrix Analysis and
  Applications, 30 (2009), pp.~1709--1734.

\bibitem{de2000best}
{\sc L.~De~Lathauwer, B.~De~Moor, and J.~Vandewalle}, {\em On the best rank-1
  and rank-($r_1,..., r_n$) approximation of higher-order tensors}, SIAM
  Journal on Matrix Analysis and Applications, 21 (2000), pp.~1324--1342.

\bibitem{dong2022fast}
{\sc H.~Dong, T.~Tong, C.~Ma, and Y.~Chi}, {\em Fast and provable tensor robust
  principal component analysis via scaled gradient descent}, Information and
  Inference: A Journal of the IMA, 12 (2023), p.~iaad019.

\bibitem{Goreinov}
{\sc S.~A. Gore\u\i{}nov, N.~L. Zamarashkin, and E.~E. Tyrtyshnikov}, {\em
  Pseudo-skeleton approximations}, Doklay Akdemii Nauk, 343 (1995),
  pp.~151--152.

\bibitem{gu2014robust}
{\sc Q.~Gu, H.~Gui, and J.~Han}, {\em Robust tensor decomposition with gross
  corruption}, Advances in Neural Information Processing Systems, 27 (2014),
  pp.~1422--1430.

\bibitem{HH2020}
{\sc K.~Hamm and L.~Huang}, {\em Perspectives on {CUR} decompositions}, Applied
  and Computational Harmonic Analysis, 48 (2020), pp.~1088--1099.

\bibitem{hamm2020stability}
{\sc K.~Hamm and L.~Huang}, {\em Stability of sampling for {CUR}
  decompositions}, Foundations of Data Science, 2 (2020), p.~83.

\bibitem{hamm2022RieCUR}
{\sc K.~Hamm, M.~Meskini, and H.~Cai}, {\em Riemannian {CUR} decompositions for
  robust principal component analysis}, in Topological, Algebraic and Geometric
  Learning Workshops 2022, PMLR, 2022, pp.~152--160.

\bibitem{HF1928}
{\sc F.~L. Hitchcock}, {\em Multiple invariants and generalized rank of a
  $p$-way matrix or tensor}, Journal of Mathematical Physics, 7 (1928),
  pp.~39--79.

\bibitem{hu2019dstpca}
{\sc Y.~Hu, J.-X. Liu, Y.-L. Gao, and J.~Shang}, {\em {DSTPCA}: Double-sparse
  constrained tensor principal component analysis method for feature
  selection}, IEEE/ACM Transactions on Computational Biology and
  Bioinformatics,  (2019).

\bibitem{hu2020robust}
{\sc Y.~Hu and D.~B. Work}, {\em Robust tensor recovery with fiber outliers for
  traffic events}, ACM Transactions on Knowledge Discovery from Data (TKDD), 15
  (2020), pp.~1--27.

\bibitem{huang2014provable}
{\sc B.~Huang, C.~Mu, D.~Goldfarb, and J.~Wright}, {\em Provable low-rank
  tensor recovery}, Optimization-Online, 4252 (2014), pp.~455--500.

\bibitem{ji2016coauthorship}
{\sc P.~Ji and J.~Jin}, {\em Coauthorship and citation networks for
  statisticians}, The Annals of Applied Statistics, 10 (2016), pp.~1779--1812.

\bibitem{l}
{\sc E.~L}, {\em Foreground detection in video sequences with probabilistic
  self-organizing maps}, \url{http://www.lcc.uma.es/~ezeqlr/fsom/fsom.html}.

\bibitem{lee1999learning}
{\sc D.~D. Lee and H.~S. Seung}, {\em Learning the parts of objects by
  non-negative matrix factorization}, Nature, 401 (1999), pp.~788--791.

\bibitem{li2004statistical}
{\sc L.~Li, W.~Huang, I.~Y.-H. Gu, and Q.~Tian}, {\em Statistical modeling of
  complex backgrounds for foreground object detection}, IEEE Transactions on
  Image Processing, 13 (2004), pp.~1459--1472.

\bibitem{li2018tucker}
{\sc X.~Li, D.~Xu, H.~Zhou, and L.~Li}, {\em Tucker tensor regression and
  neuroimaging analysis}, Statistics in Biosciences, 10 (2018), pp.~520--545.

\bibitem{Lin2017}
{\sc Z.~Lin and H.~Zhang}, in Low-Rank Models in Visual Analysis, Elsevier,
  2017, pp.~1--2.

\bibitem{liu2012tensor}
{\sc J.~Liu, P.~Musialski, P.~Wonka, and J.~Ye}, {\em Tensor completion for
  estimating missing values in visual data}, IEEE transactions on pattern
  analysis and machine intelligence, 35 (2012), pp.~208--220.

\bibitem{liu2022characterizing}
{\sc T.~Liu, M.~Yuan, and H.~Zhao}, {\em Characterizing spatiotemporal
  transcriptome of the human brain via low-rank tensor decomposition},
  Statistics in Biosciences,  (2022), pp.~1--29.

\bibitem{liu2018improved}
{\sc Y.~Liu, L.~Chen, and C.~Zhu}, {\em Improved robust tensor principal
  component analysis via low-rank core matrix}, IEEE Journal of Selected Topics
  in Signal Processing, 12 (2018), pp.~1378--1389.

\bibitem{lu2016tensor}
{\sc C.~Lu, J.~Feng, Y.~Chen, W.~Liu, Z.~Lin, and S.~Yan}, {\em Tensor robust
  principal component analysis: Exact recovery of corrupted low-rank tensors
  via convex optimization}, in Proceedings of the IEEE Conference on Computer
  Vision and Pattern Recognition, 2016, pp.~5249--5257.

\bibitem{lu2019tensor}
{\sc C.~Lu, J.~Feng, Y.~Chen, W.~Liu, Z.~Lin, and S.~Yan}, {\em Tensor robust
  principal component analysis with a new tensor nuclear norm}, IEEE
  Transactions on Pattern Analysis and Machine Intelligence, 42 (2019),
  pp.~925--938.

\bibitem{LU2011survey}
{\sc H.~Lu, K.~N. Plataniotis, and A.~N. Venetsanopoulos}, {\em A survey of
  multilinear subspace learning for tensor data}, Pattern Recognition, 44
  (2011), pp.~1540--1551.

\bibitem{luo2015subgraph}
{\sc Y.~Luo, Y.~Xin, E.~Hochberg, R.~Joshi, O.~Uzuner, and P.~Szolovits}, {\em
  Subgraph augmented non-negative tensor factorization ({SANTF}) for modeling
  clinical narrative text}, Journal of the American Medical Informatics
  Association, 22 (2015), pp.~1009--1019.

\bibitem{MMD2008}
{\sc M.~W. Mahoney, M.~Maggioni, and P.~Drineas}, {\em Tensor-{CUR}
  decompositions for tensor-based data}, SIAM Journal on Matrix Analysis and
  Applications, 30 (2008), pp.~957--987.

\bibitem{mohammadpour2022randomized}
{\sc M.~Mohammadpour~Salut and D.~Anderson}, {\em Randomized tensor robust
  principal component analysis},  (2022).

\bibitem{netrapalli2014non-convex}
{\sc P.~Netrapalli, N.~U~N, S.~Sanghavi, A.~Anandkumar, and P.~Jain}, {\em
  Non-convex robust {PCA}}, in Advances in Neural Information Processing
  Systems, vol.~27, 2014.

\bibitem{oh2011large}
{\sc S.~Oh, A.~Hoogs, A.~Perera, N.~Cuntoor, C.-C. Chen, J.~T. Lee,
  S.~Mukherjee, J.~Aggarwal, H.~Lee, L.~Davis, et~al.}, {\em A large-scale
  benchmark dataset for event recognition in surveillance video}, in CVPR 2011,
  IEEE, 2011, pp.~3153--3160.

\bibitem{o2005video}
{\sc A.~J. O'Toole, J.~Harms, S.~L. Snow, D.~R. Hurst, M.~R. Pappas, J.~H.
  Ayyad, and H.~Abdi}, {\em A video database of moving faces and people}, IEEE
  Transactions on Pattern Analysis and Machine Intelligence, 27 (2005),
  pp.~812--816.

\bibitem{rabanser2017introduction}
{\sc S.~Rabanser, O.~Shchur, and S.~G{\"u}nnemann}, {\em Introduction to tensor
  decompositions and their applications in machine learning}, arXiv preprint
  arXiv:1711.10781,  (2017).

\bibitem{rajwade2012image}
{\sc A.~Rajwade, A.~Rangarajan, and A.~Banerjee}, {\em Image denoising using
  the higher order singular value decomposition}, IEEE Transactions on Pattern
  Analysis and Machine Intelligence, 35 (2012), pp.~849--862.

\bibitem{sofuoglu2018two}
{\sc S.~E. Sofuoglu and S.~Aviyente}, {\em A two-stage approach to robust
  tensor decomposition}, in 2018 IEEE Statistical Signal Processing Workshop
  (SSP), IEEE, 2018, pp.~831--835.

\bibitem{song2017based}
{\sc T.~Song, Z.~Peng, S.~Wang, W.~Fu, X.~Hong, and P.~S. Yu}, {\em Based
  cross-domain recommendation through joint tensor factorization}, in
  International conference on database systems for advanced applications,
  Springer, 2017, pp.~525--540.

\bibitem{tucker1966}
{\sc L.~R. Tucker}, {\em Some mathematical notes on three-mode factor
  analysis}, Psychometrika, 31 (1966), pp.~279--311.

\bibitem{vaswani2018static}
{\sc N.~Vaswani and P.~Narayanamurthy}, {\em Static and dynamic robust {PCA}
  and matrix completion: A review}, Proceedings of the IEEE, 106 (2018),
  pp.~1359--1379.

\bibitem{wright2008robust}
{\sc J.~Wright, A.~Y. Yang, A.~Ganesh, S.~S. Sastry, and Y.~Ma}, {\em Robust
  face recognition via sparse representation}, IEEE Transactions on Pattern
  Analysis and Machine Intelligence, 31 (2008), pp.~210--227.

\bibitem{xia2021statistically}
{\sc D.~Xia, M.~Yuan, C.-H. Zhang, et~al.}, {\em Statistically optimal and
  computationally efficient low rank tensor completion from noisy entries},
  Annals of Statistics, 49 (2021), pp.~76--99.

\bibitem{zhang2014novel}
{\sc Z.~Zhang, G.~Ely, S.~Aeron, N.~Hao, and M.~Kilmer}, {\em Novel methods for
  multilinear data completion and de-noising based on tensor-{SVD}}, in
  Proceedings of the IEEE Conference on Computer Vision and Pattern
  Recognition, 2014, pp.~3842--3849.

\bibitem{zheng2016topic}
{\sc X.~Zheng, W.~Ding, Z.~Lin, and C.~Chen}, {\em Topic tensor factorization
  for recommender system}, Information Sciences, 372 (2016), pp.~276--293.

\bibitem{zhou2013tensor}
{\sc H.~Zhou, L.~Li, and H.~Zhu}, {\em Tensor regression with applications in
  neuroimaging data analysis}, Journal of the American Statistical Association,
  108 (2013), pp.~540--552.

\bibitem{zhou2017tensor}
{\sc P.~Zhou, C.~Lu, Z.~Lin, and C.~Zhang}, {\em Tensor factorization for
  low-rank tensor completion}, IEEE Transactions on Image Processing, 27
  (2017), pp.~1152--1163.

\end{thebibliography}

\end{document}